\documentclass[12pt, a4paper]{article}
\usepackage{amsmath,amssymb}
\usepackage[nottoc]{tocbibind} 
\usepackage{bm,bbm}
\usepackage{comment}
\usepackage{graphicx}
\usepackage{ascmac}
\usepackage{tikz}
\usetikzlibrary {arrows.meta}
\usepackage{xypic}
\usepackage{tikz-cd}
\usetikzlibrary{matrix, calc, arrows}

\usepackage{adjustbox}
\usepackage{amsmath,amsfonts,amssymb,mathrsfs}   
\usepackage[figurename=Fig.]{caption}
\usepackage{amsthm}
\usepackage{enumerate}
\usepackage{color}
\definecolor{tab:blue}{RGB}{31, 119, 180}
\definecolor{tab:orange}{RGB}{255, 127, 14}
\definecolor{tab:green}{RGB}{44, 160, 44}
\definecolor{tab:red}{RGB}{214, 39, 40}
\definecolor{tab:purple}{RGB}{148, 103, 189}
\definecolor{tab:brown}{RGB}{140, 86, 75}
\definecolor{tab:pink}{RGB}{227, 119, 194}
\definecolor{tab:gray}{RGB}{127, 127, 127}
\definecolor{tab:olive}{RGB}{188, 189, 34}
\definecolor{tab:cyan}{RGB}{23, 190, 207}
\usepackage[colorlinks,linkcolor=blue, citecolor=tab:red, urlcolor=cyan, 
]{hyperref}

\usepackage{geometry}
\usepackage{lineno}

\theoremstyle{plain}
\newtheorem{mainthm}{Theorem}

\newtheorem{mainprop}[mainthm]{Proposition}

\newtheorem{maincor}[mainthm]{Corollary}
\newtheorem{theorem}{Theorem}[subsection]
\newtheorem{lemma}[theorem]{Lemma}
\newtheorem{proposition}[theorem]{Proposition}
\newtheorem{corollary}[theorem]{Corollary}

\theoremstyle{definition}
\newtheorem{definition}[theorem]{Definition}
\newtheorem{example}[theorem]{Example}
\newtheorem{assumption}[theorem]{Assumption}

\theoremstyle{remark}
\newtheorem{remark}[theorem]{Remark}

\numberwithin{equation}{subsection}
\renewcommand\thesection{\arabic{section}}
\renewcommand\thesubsection{\arabic{section}.\arabic{subsection}}

\newcommand{\Ad}{\mathop{\mathrm{Ad}}\nolimits}

\newcommand{\Aut}{\mathop{\mathrm{Aut}}\nolimits}

\newcommand{\bR}{\mathbb{R}}

\newcommand{\GL}{\mathop{\mathrm{GL}}\nolimits}
\newcommand{\SL}{\mathop{\mathrm{SL}}\nolimits}
\newcommand{\SU}{\mathop{\mathrm{SU}}\nolimits}
\newcommand{\PSL}{\mathop{\mathrm{PSL}}\nolimits}

\newcommand{\Ker}{\mathop{\mathrm{Ker}}\nolimits}

\newcommand{\ord}{\mathop{\mathrm{ord}}\nolimits}

\newcommand{\sign}{\mathop{\mathrm{sign}}\nolimits}

\newcommand{\Tr}{\mathop{\mathrm{Tr}}\nolimits}

\newcommand{\Id}{\mathop{\mathrm{Id}}\nolimits}

\newcommand{\bb}{\mathfrak}

\newcommand{\Zdam}{Z_{\mathrm{O}\hspace{-.175em}-\hspace{-.175em}\mathrm{O}}}

\newcommand{\Darr}[3][0]{
\draw[->,>=stealth, rotate=#1] (#2,#3)--(#2+0.051,#3);
}

\newcommand\blfootnote[1]{
  \begingroup
  \renewcommand\thefootnote{}\footnote{#1}
  \addtocounter{footnote}{-1}
  \endgroup
}

\makeatletter
\newcommand{\dates}[1]{%
  \let\@@@oldtitle\@title%
  \gdef\@title{\@@@oldtitle\footnotetext{\emph{Date:} #1.}}%
}
\makeatother

\title{\bf Integral invariants for framed $3$-manifolds associated to trivalent graphs possibly with self-loops}
\author{Hisatoshi Kodani and Bingxiao Liu}
\date{}

\begin{document}
\maketitle

\abstract
Bott--Cattaneo's theory defines the integral invariants for a framed rational homology $3$-sphere equipped with an acyclic orthogonal local system, in terms of graph cocycles without self-loops. The $2$-loop term of their invariants is associated with the theta graph. Their definition requires a cohomological condition. Cattaneo--Shimizu removed this cohomological condition and gave a $2$-loop invariant associated with a linear combination of the theta graph and the dumbbell graph, the $2$-loop trivalent graph with self-loops. In this paper, we are concerned with an acyclic local system given by the adjoint representation of a semi-simple Lie group composed with a representation of the fundamental group of a closed $3$-manifold, and we show that through a cohomological construction eventually the integral associated with the dumbbell graph vanishes. Based on this idea, we construct a theory of graph complexes and cocycles, so that higher-loop invariants can be defined by two different but equivalent methods: the graph cocycles without self-loops as in Bott--Cattaneo's theory, and the ones with self-loops that extend Cattaneo--Shimizu's $2$-loop invariants. As a consequence, we prove that the generating series of Chern--Simons perturbation theory gives rise to topological invariants for framed $3$-manifolds in our setting, which admits a formula in terms of only trivalent graphs without self-loops.

\setcounter{tocdepth}{2} 
\tableofcontents

\section{Introduction}
The mathematical foundation of Chern--Simons perturbation theory was developed by Axelrod--Singer \cite{AS,AS2} and Kontsevich \cite{Ko} around 1990s after the breakthrough of Witten's work \cite{MR990772} on the Jones polynomials via Chern--Simons theory. The Chern--Simons perturbation theory produces a family of topological invariants, parametrized by certain linear combinations of trivalent graphs. More precisely, these invariants pertain to closed, oriented, connected smooth $3$-manifolds, denoted $M$. These invariants are associated with a homotopy class of smooth framings of $M$ (that is, smooth trivializations of the tangent bundle $TM$) and an acyclic local system $E_{\rho}$ over $M$ (that is, $H^k(M; E_{\rho}) =0$ for $0 \leq k \leq 3$). The local system concerned is $E_\rho:= \pi_{1}(M)\backslash \big(\widetilde{M}\times_{\rho} \mathfrak{g}\big)$ defined by a representation $\rho: \pi_1(M) \rightarrow G \overset{\Ad}{\rightarrow} \Aut(\mathfrak{g})$, where $\widetilde{M}$ denotes the universal covering space of $M$, $\mathfrak{g}$ is the Lie algebra of a semi-simple (connected) Lie group $G$, and $\Ad$ denotes the adjoint representation of $G$. Then the aforementioned invariants are defined as the configuration space integrals of the propagators, and the way of defining integrands is encoded by the trivalent graphs. We will always refer to such invariants as integral invariants for $M$ and $E_\rho$.

Inspired by their work, Bott and Cattaneo \cite{BC, BC2} introduced topological invariants of framed rational homology 3-spheres equipped with acyclic orthogonal local systems. Invariants (of arbitrary orders) in  \cite{BC, BC2} are defined in terms of configuration space integrals associated to appropriate graph cocycles without self-loops (self-loops here are called simple loops in \cite{Ko}).

Then Cattaneo and Shimizu in their paper \cite{CS} pointed out that there is a gap in the arguments of \cite[Lemma 1.2]{BC2} due to a possible non-vanishing cohomology group. More precisely, the invariants constructed in \cite{BC2} were based on the implicit assumption of the vanishing of a cohomology group, that is, $H^2_{-}(\Delta; \pi_1^{-1}E_{\rho} \otimes \pi_2^{-1}E_{\rho})=0$ in \cite{BC2}. In \cite{CS} and in the present paper, this cohomology group is denoted as $H^2_{-}(\Delta; E_{\rho} \otimes E_{\rho})$, where $\Delta\simeq M$ is the diagonal of $M\times M$. This motivated Cattaneo and Shimizu \cite{CS} to correct the $2$-loop term of Bott--Cattaneo invariants \cite{BC2} with $H^2_{-}(\Delta; E_{\rho} \otimes E_{\rho})\neq 0$, where the dumbbell graph appears in the construction.

The non-vanishing of the cohomology group $H^2_{-}(\Delta; E_{\rho} \otimes E_{\rho})$ implies that the boundary value of a propagator can have a factor of a non-trivial antisymmetric form $\xi$, which is analogous to the regular part of the propagator constructed in \cite{AS,AS2} via the Green kernels. This regular form $\xi$ leads to configuration space integrals associated with trivalent graphs with self-loops. Such graphs are not needed in \cite{BC, BC2} under the assumption $H^2_{-}(\Delta; E_{\rho} \otimes E_{\rho})=0$. The refined 2-loop term defined in \cite{CS} is eventually defined by a linear combination of two terms, $Z_{\Theta}$ and $\Zdam$, which are the configuration space integrals corresponding to the theta graph and the dumbbell graph, respectively. The form $\xi$ is the key factor associated with the self-loops of the dumbbell graph.  

Note that Shimizu \cite{Shi} also showed that, when $G=\SU(2)$, the cohomology group $H^2_{-}(\Delta; E_{\rho} \otimes E_{\rho})$ always vanishes (see also Proposition \ref{prop:3.5.2}) and $Z_{\Theta}$ itself becomes an invariant of closed 3-manifolds with any representation $\rho: \pi_1(M) \rightarrow \SU(2) \overset{\Ad}{\rightarrow} \Aut(\mathfrak{su}_2)$. Also note that a class of regular form $\xi$ has recently been studied by Kitano and Shimizu (see \cite{Shi3}, \cite{KS}) which is denoted by $d(M, \rho)$, motivated by its relation to Reidemeister torsion expected from the viewpoint of quantum Chern--Simons theory. Their study can be regarded as an attempt to generalize Lescop's result for closed 3-manifolds with the first Betti number 1 and an abelian representation (\cite{Les10}), where her invariant is denoted by $I_{\Delta}$ and described by the logarithmic derivative of Alexander polynomial, to those with arbitrary first Betti number and non-abelian representations.

Therefore, it is important to ask for the existence of examples of a pair of a closed 3-manifold and an acyclic local system via an adjoint representation as above, which have non-vanishing $H^2_{-}(\Delta; E_{\rho} \otimes E_{\rho})$. 
Such examples would ensure that the refinement by Cattaneo--Shimizu in \cite{CS} is meaningful, and now our first result of this paper is to give a positive answer.

\begin{mainprop}\label{mainprop1}
There exists a triple $(M, G, \rho)$, consisting of a closed oriented smooth manifold $M$, a semi-simple Lie group $G$, and an acyclic representation $\rho$ of $\pi_1(M)$ through an adjoint representation as above, such that it has $H^2_{-}(\Delta; E_{\rho} \otimes E_{\rho}) \neq 0$.
\end{mainprop}

To the best of our knowledge, we are the first to consider such examples as stated in Proposition \ref{mainprop1}. 
In fact, we state a more concrete version of Proposition \ref{mainprop1} as Proposition \ref{prop:3.-1}, where the expected examples are constructed explicitly for the tripe $(M, G, \rho)$ with the group $G=\SL(2, \mathbb{C})\times \SL(2, \mathbb{C})$.

Following this direction with possibly non-vanishing $H^2_{-}(\Delta; E_{\rho} \otimes E_{\rho})$ from \cite{CS}, we further investigate the general integral invariants of Cattaneo--Shimizu/Bott--Cattaneo defined from trivalent graphs (with higher loops); 
in particular, we aim to understand the role of trivalent graphs with self-loops that are excluded in the original work of Bott--Cattaneo \cite{BC2}.  

Our second result follows from the reexamination of the $2$-loop invariant introduced in \cite{CS}, we found that, even when $H^2_{-}(\Delta; E_{\rho} \otimes E_{\rho}) \neq 0$, there is a special choice of propagator by which the integral associated with the dumbbell graph vanishes. 
Then, roughly speaking, the essential part of this $2$-loop invariant in \cite{CS} is reduced back to the term of the theta graph, as in \cite{BC2}. 

Let us start with the definition of propagators. 
Let $C_2(M)$ denote the compactified $2$-point configuration space of $M$, we can consider it as the real blow-up of $M^2$ along the diagonal $\Delta\simeq M$. 
The manifold $C_2(M)$ is a smooth manifold with boundary, and its boundary $\partial C_2(M)$ can be identified with the sphere tangent bundle $S(TM)$ of $M$. 
We will denote by $\mathfrak{i}_\partial: \partial C_2(M)\rightarrow C_2(M)$ the inclusion.

Let $q:C_2(M)\rightarrow M\times M$ denote the blow-down map, and let $q_\partial: \partial C_2(M)\rightarrow M$ denote its restriction to the boundary, which is a smooth fibration with a fiber $\mathbb{S}^2$. We always fix an orientation $o(M)$ for $M$, and we also fix a smooth framing $f$ of $M$, that is, a smooth identification of vector bundles $TM\simeq M\times \bR^3$ over $M$. 
In this way, we identify $\partial C_2(M)\simeq M\times\mathbb{S}^2$.

We always fix a connected semi-simple Lie group $G$ with Lie algebra $\mathfrak{g}$. Note that $G$ could be noncompact. As we mentioned, we consider a representation $\rho: \pi_1(M) \rightarrow G \overset{\Ad}{\rightarrow} \Aut(\mathfrak{g})$ and the associated local system $E_\rho$ over $M$. 
Consequently, we have the induced tensor bundle $E_\rho\boxtimes E_\rho$ on $M\times M$, hence after taking its pullback bundle by the blow-down map $q$, we get a flat vector bundle $F_\rho=q^*(E_\rho\boxtimes E_\rho)$ on $C_2(M)$. 
The restriction of $F_\rho$ on $\partial C_2(M)\simeq M\times\mathbb{S}^2$ is just the pullback of vector bundle $E_\rho\otimes E_\rho\rightarrow M$. 
Meanwhile, we also define an involution $T$ acting on $M\times M$ and on $E_\rho\boxtimes E_\rho$ by swapping two factors in the product or tensor product, this action lifts to $F_\rho\rightarrow C_2(M)$ that commutes with the flat connection, then acts on the sections of $F_\rho$ and the associated de Rham cohomology groups $H^\bullet(M; E_\rho)$.

Assume $E_\rho$ to be acyclic, that is, $H^\bullet(M; E_\rho)=0$. 
In this case, a propagator is a closed differential form $\omega\in \Omega^2(C_2(M);F_\rho)$ such that 
\begin{itemize}
    \item $\omega$ is anti-symmetric, i.e., $T^\ast \omega=-\omega$ or we say $\omega\in \Omega^2_-(C_2(M);F_\rho)$, where the subscript $-$ corresponds to $(-1)$-eigenvalue of the action of $T$;
    \item there is a normalized volume form $\eta$ on $\mathbb{S}^2$ (i.e., with volume $1$) and a closed smooth form $\xi\in \Omega^2_-(M; E_\rho\otimes E_\rho)$ such that
    \begin{equation}
    \mathfrak{i}_\partial^\ast(\omega)=\eta\otimes \mathbf{1}+q_\partial^\ast(\xi),
        \label{eq:intro1}
    \end{equation}
    where $\eta$ is viewed as a fibrewise vertical volume form along the fibration $q_\partial: \partial C_2(M)\rightarrow M$ (hence depending on the framing $f$), and $\mathbf{1}$ is a flat section, called Casimir section (see Lemma \ref{lm:1.1}), of $E_\rho\otimes E_\rho$ over $M$.
\end{itemize}
Then, any chosen propagator $\omega$ gives rise to the 2-loop invariant as the linear combination of certain integrations over $C_2(M)$ associated with the theta graph and the dumbbell graph:
\begin{equation}
     \frac{1}{12} Z_{\Theta}(\omega, \xi) - \frac{1}{8} \Zdam(\omega,\xi).
\end{equation}
 Note that, in \eqref{eq:intro1}, the closed form $\xi$ defines a cohomology class $[\xi]\in H^2_{-}(M; E_{\rho} \otimes E_{\rho})$. In general, one cannot take $\xi = 0$ when $H^2_{-}(M; E_{\rho} \otimes E_{\rho})\neq 0$. This term in the boundary condition for a propagator was missing in \cite[Lemma 1.2]{BC2} and then studied by \cite{CS} to define the $2$-loop integral invariant $\Zdam(\omega,\xi)$ associated to the dumbbell graph.

As is well known, for the case of a rational homology $3$-sphere with the trivial local system, the dumbbell term $\Zdam(\omega,\xi)$ vanishes by the antisymmetric relation of Jacobi diagrams. Jacobi diagrams are a kind of abstraction of the Lie algebraic part of the perturbation theory around the trivial connection. As indicated by examples in Proposition \ref{mainprop1}, it is not the same case for general non-trivial local systems, and such antisymmetric relation cannot directly imply the vanishing of $\Zdam(\omega,\xi)$. However, in the present paper, we notice that the Lie algebraic structure (in particular, the fibre-wise semi-simplicity) of $E_{\rho}$ leads to the vanishing of the dumbbell term $\Zdam(\omega,\xi)$. More precisely, let $\mathfrak{L}: E_\rho\otimes E_\rho\rightarrow E_\rho$ be the vector bundle morphism that corresponds fibre-wisely the Lie bracket on $\mathfrak{g}$, we find that, by the semi-simplicity of $\mathfrak{g}$, there exists a propagator $\omega^{\sharp}$ whose regular part $\xi^\sharp$ in its boundary value satisfies an additional condition
\begin{equation}
    \mathfrak{L}(\xi^\sharp)=0.
    \label{eq:1.0.2nov23}
\end{equation}
Then, using such $\omega^{\sharp}$, we refine Cattaneo and Shimizu's result \cite[Theroem 2.3]{CS} as follows. 
\begin{mainthm}[see Theorem \ref{prop:3.2}]\label{thm:A}
Let $M$ be a closed, connected, orientable 
smooth $3$-manifold. Fix a homotopy class $[f]$ of smooth framings of $M$ and an orientation $o(M)$. Let $E_{\rho}$ be an acyclic local system over $M$ associated with a representation $\rho : \pi_1(M) \rightarrow G \overset{\Ad}{\longrightarrow} \Aut(\mathfrak{g})$. Then, for any  propagator $\omega^\sharp$ with the regular form $\xi^\sharp$ such that $\mathfrak{L}(\xi^\sharp)=0$, we have $\Zdam(\omega^\sharp,\xi^\sharp)=0$. In other words, the theta invariant $Z_{\Theta}(\omega^\sharp)$  gives a 2-loop invariant for a framed closed 3-manifold $M$ and $\rho$.
\end{mainthm}
Notice that the above result also recovers the corresponding result of Bott and Cattaneo \cite{BC2} since it is always possible to take $\xi^\sharp\equiv 0$ when $H^2_-(M;E_\rho\otimes E_\rho)=0$ and the condition \eqref{eq:1.0.2nov23} is automatically satisfied.  In Definition \ref{def:6.1.1}, for the acyclic $E_\rho$, such a propagator $\omega^\sharp$ is called adapted (to $E_\rho$). The condition \eqref{eq:1.0.2nov23} for adapted propagators will make the integrals vanish not only for the dumbbell graph but also for general trivalent graphs with self-loops. Based on this observation, we generalize the result in Theorem \ref{thm:A} to higher-loop invariants associated with trivalent graphs with or without self-loops, formulated in terms of certain graph complexes. 

In Section \ref{section:graph}, a graph complex (over $\mathbb{Q}$) of decorated graphs (see Definition \ref{def:6.3.1August}) possibly with self-loops, denoted by $(\mathcal{GC}_{\mathfrak{g}},\delta)$,  is defined, which is dedicated to acyclic local systems and only depends on the semi-simple Lie algebra $\mathfrak{g}$. The differential operator $\delta$ is given by contraction on each non-self-loop edge of the decorated graphs. In our convention, each vertex of the graph has at least $3$ incident half-edges, so that the degree-$0$ subspace of $\mathcal{GC}_{\mathfrak{g}}$ is exactly spanned by the decorated trivalent graphs, which is isomorphic to the linear space spanned by topological trivalent graphs (see Subsection \ref{ss6.6kk}). The order of a trivalent graph is defined as half of its total number of vertices.

Then in Subsection \ref{sec:gc_acyclic} we construct two associated graph complexes: one is the subcomplex $(\mathcal{GC}^\prime_{\mathfrak{g}},\delta)$ of the above one which consists of decorated graphs with at least one self-loop, and the second complex $(\mathcal{G}_{\mathfrak{g}},\delta^\sharp)$ is the quotient complex $\mathcal{GC}_{\mathfrak{g}}/\mathcal{GC}^\prime_{\mathfrak{g}}$, which consists of linear combinations of the decorated graphs without self-loops.

 The graph complex $\mathcal{GC}_{\mathfrak{g}}/\mathcal{GC}^\prime_{\mathfrak{g}}$ was already sketched in \cite[\S 4 Discussion]{BC2} without developing the details and inspired from the general ideas of Kontsevich \cite{MR1247289}\cite[Section 2]{Ko}. Here, different from \cite[\S 4]{BC2}, the graphs with self-loops are necessary in our graph complex $(\mathcal{GC}_{\mathfrak{g}},\delta)$ and we will give more down to earth graphical construction incorporating with weight systems given as in \cite{BN}. Note that our decorated graphs constructed in Subsection \ref{ss6.4complex} can be viewed as a variant of the graphs described by Conant and Vogtmann \cite{MR2026331} for the Lie type graph complexes, following the work of Kontsevich. However, a different grading from  \cite[\S 3]{MR2026331} is used on the graph spaces for defining our complex $(\mathcal{GC}_{\mathfrak{g}},\delta)$ in order to serve the configuration space integrals (see Remark \ref{rk:new24}).

An element in $\Ker \delta \subset \mathcal{GC}_{\mathfrak{g}}$ or in $\Ker \delta^\sharp \subset \mathcal{G}_{\mathfrak{g}}$ is called a {\it graph cocycle} in the respective graph complexes. Then we are mainly concerned with the cocycles of degree $0$ and order $n\geq 1$ (hence with $(n+1)$ loops): $H^0(\mathcal{GC}^\bullet_{\mathfrak{g}: n}, \delta)$, $H^0(\mathcal{G}^\bullet_{\mathfrak{g}: n}, \delta^\sharp)$. In particular, the $2$-loop cocycles $H^0(\mathcal{GC}^\bullet_{\mathfrak{g}: 1}, \delta)$ is $1$-dimensional and spanned by a linear combination of the theta graph and dumbbell graph, which is implicitly used in \cite{CS}.

Let $C_n(M)$ denote the Fulton--MacPherson compactification of configuration space of $n$-points in $M$ (see \ref{ss1.1corners} and \ref{app:FMAS} for detailed defintions). Associated with the graph cocycles of degree $0$ and order $n\geq 1$, we can define an integral on $C_{2n}(M)$ in terms of the propagators, see \eqref{eq:6.2.4paris}, where the integrands are encoded by the graphs. Now we can state the main result of the paper.

\begin{mainthm}[See Theorems \ref{thm:6.2.5ss} and \ref{thm:6.3.1}]\label{thm:CC}
Fix a homotopy class $[f]$ of smooth framings of $M$ and an orientation $o(M)$. Let $E_{\rho}$ be an acyclic local system on $M$ corresponding to $\rho : \pi_1(M) \rightarrow G \overset{\Ad}{\longrightarrow} \Aut(\mathfrak{g})$ as above. Fix an order $n\geq 1$. Any cocycle  $H^0(\mathcal{GC}^\bullet_{\mathfrak{g}: n}, \delta)$ (that is, $\delta \Gamma=0$), taking the associated configuration space integrals of any propagator for $E_\rho$,  gives rise to an integral invariant $Z(M,\rho,[f])(\Gamma)\in\bR$, which is independent of the choice of propagator used to define it. Similarly, a cocycle $\Gamma' \in H^0(\mathcal{G}^\bullet_{\mathfrak{g}: n}, \delta^\sharp)$ consisting of connected trivalent graphs without self-loops, using the adapted propagator to define the associated configuration space integral, gives rise to an integral invariant $Z^\sharp(M,\rho,[f])(\Gamma')\in\bR$.

Moreover, regarding $Z(M,\rho,[f])$, $Z^\sharp(M,\rho,[f])$ as linear functionals on the cocycles, the following diagram commutes:

\begin{equation}
\begin{tikzcd}
H^0(\mathcal{GC}^\bullet_{\mathfrak{g}: n}, \delta) \arrow[r, "{Z(M,\rho, [f])}"]\arrow[d] &[1.5cm] \bR \arrow[d,"="]\\
H^0(\mathcal{G}^\bullet_{\mathfrak{g}: n}, \delta^\sharp) \arrow[r, "{Z^\sharp(M,\rho, [f])}"] & \bR
\end{tikzcd}
\end{equation}
where the leftmost vertical map is given by sending graphs with self-loops to zero.
\end{mainthm}

These results imply that the computation of integral invariants associated with trivalent graphs possibly with self-loops can be reduced to those associated with trivalent graphs without self-loops. 
Again note that the reason for the non-necessity of self-loop in the graphs for invariant $Z^\sharp(M,\rho,[f])$ is different from the trivial local system case, for instance, described in \cite[\S\,The graph complex]{Ko}.

Now we explain how Chern--Simons perturbation theory fits into our results. The general idea from Chern--Simons perturbation theory is to define the topological invariants for $3$-manifolds via the configuration space integrals associated with a generating series in terms of the linear combinations of trivalent graphs (see \cite[Section 2]{Ko}, \cite{AS,AS2}). For a given order $n\geq 1$, the corresponding term in the generating series is given as
\begin{equation}
     \sum_{\substack{\text{connected trivalent } \mathfrak{G}\\ \mathrm{ord}(\mathfrak{G})=n}} \pm \frac{1}{|\Aut(\mathfrak{G})|}   \Gamma(\mathfrak{G})\in \mathcal{GC}^0_{\mathfrak{g}: n},
\label{eq:intro6new}
\end{equation}
where the sum runs over all the connected topological trivalent graph $\mathfrak{G}$ of order $n$, and $\Gamma(\mathfrak{G})$ denotes an arbitrary decorated graph with the underlying topological graph $\mathfrak{G}$ with the sign $\pm$ determined by it (see Corollary \ref{Cor:4.6.5-july}).

In Proposition \ref{prop:CSseries}, we prove that the element \eqref{eq:intro6new} is a cocycle in $\mathcal{GC}^0_{\mathfrak{g}: n}$; if we remove all the terms with self-loops from  \eqref{eq:intro6new}, we get a cocycle in $\mathcal{G}^0_{\mathfrak{g}: n}$. Combining these assertions with Theorem \ref{thm:CC}, we get the following result.

\begin{maincor}[See Corollary \ref{cor:generating_func}]
Let $M$ be a closed, connected, orientable smooth $3$-manifold. Fix a homotopy class $[f]$ of smooth framings of $M$ and an orientation $o(M)$. Let $E_{\rho}$ be an acyclic local system on $M$ associated with a representation $\rho : \pi_1(M) \rightarrow G \overset{\Ad}{\longrightarrow} \Aut(\mathfrak{g})$. Then the following formal series
\begin{equation}
   \log \mathcal{Z}_{\mathrm{CS}}(M, \rho, [f])= \sum_{\text{connected trivalent } \mathfrak{G}} \frac{\hbar^{\mathrm{ord}(\mathfrak{G})}}{|\Aut(\mathfrak{G})|}  Z(M,\rho, [f])(\pm \Gamma(\mathfrak{G})) \in \bR[[\hbar]]
\label{eq:intro7new}
\end{equation}
is a topological invariants for $(M,o(M),[f])$ and $\rho$. This series $\log \mathcal{Z}_{\mathrm{CS}}(M, \rho, [f])$ can also be written in terms of $Z^\sharp(M, \rho, [f])$ associated with the connected trivalent topological graphs without self-loops, that is,
\begin{equation}
   \log \mathcal{Z}_{\mathrm{CS}}(M, \rho, [f])= \sum_{\substack{\text{connected trivalent } \mathfrak{G}\\ \text{without self-loop}}} \frac{\hbar^{\mathrm{ord}(\mathfrak{G})}}{|\Aut(\mathfrak{G})|}  Z^\sharp(M,\rho, [f])(\pm \Gamma(\mathfrak{G})) \in \bR[[\hbar]].
\label{eq:intro8new}
\end{equation}
\end{maincor}

Finally, we remark that our construction is applicable to both real and complex semi-simple Lie groups. Thus, for simplicity of arguments, this paper mainly focuses on the case of real semi-simple Lie groups unless otherwise stated.

The organization of this paper is as follows.

In Section \ref{section3}, we investigate the de Rham cohomology groups of a local system $E_\rho$ associated with a representation $\rho : \pi_1(M) \rightarrow G \overset{\Ad}{\longrightarrow} \Aut(\mathfrak{g})$, especially, the consequences of the Lie bracket structure on $E_\rho$.

In Section \ref{section4}, we recall the existence of propagators and introduce the notion of adapted propagators for the acyclic local system $E_\rho$. Then, after revisiting the result of Cattaneo--Shimizu for their $2$-loop integral invariant, we explain how the use of an adapted propagator kills the dumbbell term.

In Section \ref{section:graph}, we introduce our version of graph complexes (which only involve the Lie algebra $\mathfrak{g}$) which is used in our theory of integral invariants. In particular, the generating series of Chern--Simons perturbation theory is a cocycle.

In Section \ref{section6}, we prove that for $E_\rho$ being acyclic, each cocycle in our complexes defines an integral invariant for the framed $3$-manifold $M$ and the representation $\rho$. The use of adapted propagators reduces the cocycles to the ones without any self-loops.

 In Appendix \ref{sectionsecond}, for the convenience of the readers, we recall the basis for the compactified configuration spaces from the viewpoint of smooth manifolds with corners.

\subsection*{Acknowledgments}
The authors would like to express their appreciation to Tatsuro Shimizu for his motivating comments, discussions, and suggestions on this project and to Tadayuki Watanabe for inspiring discussions, in particular on finding examples of local systems with specific properties. 
The authors also thank Yuji Terashima for his encouragement and helpful communication on this topic.

Many of the work in the first half of this paper was done during the authors' stay at Max Planck Institute for Mathematics (MPIM) in Bonn. The authors thank MPIM for its hospitality and financial support. B.X. Liu also gratefully acknowledges the support of DFG Priority Program 2265 ‘Random Geometric Systems’ for his stay in Universit\"{a}t zu K\"{o}ln.

\subsection*{Notation and convention}
Let $\mathbb{N}$ denote the set of natural numbers without zero. We denote the cardinality of a set $S$  by $|S|$. For a non-negative integer $r \geq 1$, let $\mathbb{S}^{r-1}$ denote the unit sphere with the induced metric from the $r$-dimensional Euclidean space $\bR^{r}$.

For a graded vector space $V=\bigoplus_{i \in \mathbb{Z}} V^i$, the degree of a homogeneous element $v \in V^i$ is denoted by $\deg(v)=i$. For two graded algebras $A$ and $B$ over a field, let $A \widehat{\otimes} B$ stand for the graded tensor product of $A$ and $B$ over the field, i.e., its underlying vector space is the tensor product $A \otimes B$ and product structure is given by the linear extension of $(a \otimes b)\cdot (a' \otimes b') = (-1)^{\deg(b)\cdot  \deg(a')} (aa' \otimes b b')$ for homogeneous elements $a, a' \in A$ and $b, b' \in B$. Depending on the context, $A^{\bullet}$ denotes a graded vector space $A^{\bullet}= \bigoplus_{i \in \mathbb{Z}} A^i$, a cochain complex $(\bigoplus_{i \in \mathbb{Z}} A^i, \partial)$, or arbitrary homogeneous component of a graded vector space $\bigoplus_{i \in \mathbb{Z}} A^i$, where $\bullet$ plays a role of the placeholder of homogeneous degree.

For an oriented smooth manifold with boundary $X$ and a local system $E$ of a real or complex vector space on $X$, we denote by $\Omega^{\bullet}(X)$ (resp. $\Omega^{\bullet}(X; E)$) the differential graded commutative algebra of (resp. $E$-valued) smooth differential forms on $X$. For a commutative ring $R$ with a unit and an $R$-module $V$, the constant local system (trivial local system) on $X$ with fiber $V$ is denoted by $\underline{V}_X$ or simply $\underline{V}$. In this paper, we often identify a flat vector bundle and its corresponding local system. 

For an oriented smooth manifold (with corners) $X$, the boundary $\partial X$ is oriented by the \textit{outword normal first} convention so that the Stokes' formula is given by $\int_X d\omega = \int_{\partial X} \omega$ for $\omega \in \Omega^{\dim X -1}(X)$.



\section{On the cohomology group $H^2(M; E_{\rho} \otimes E_{\rho})$}\label{section3}
This section studies a cohomology group $H^2(M; E_{\rho} \otimes E_{\rho})$ of $M$ associated with a local system $E_{\rho}$ defined from the adjoint action on a semi-simple Lie algebra. 
As we will see in the subsequent section, this cohomology group contains the class of regular part  $\xi$ of propagators. 
In particular, we study  $H^2(M; E_{\rho} \otimes E_{\rho})$ from the viewpoints of its Lie bracket structure in detail. 

\subsection{Preliminaries on the local system $E_{\rho}$}\label{ss:2.1-july}
This subsection recalls several basic facts on local systems $E_{\rho}$ given via adjoint representations of semi-simple Lie groups, for example, de Rham cohomology groups of such a local system, the diagonal class, the Killing form and its associated cubic trace form, and Lie structure on $E_{\rho}$. At last, we give an important result in Proposition \ref{prop:3.-1} on the nonvanishing of the cohomology group $H_{-}^2(M;E_{\rho} \otimes E_{\rho})$, this cohomology group plays an important role in our construction of (adapted) propagators in the next sections.

Let $G$ be a connected (real) semi-simple Lie group and let $\mathfrak{g}$ be its Lie algebra. Let $B:\mathfrak{g}\times\mathfrak{g}\rightarrow \bR$ denote the corresponding Killing form. Let $\rho:\pi_1(M)\rightarrow G$ be a morphism of groups, composing with the adjoint action $\mathrm{Ad}: G\rightarrow \mathrm{Aut}(\mathfrak{g})$, we get a representation of $\pi_1(M)$ with representation space $\mathfrak{g}$, which is still denoted by $\rho$. 
The representation $\rho$ defines canonically a flat vector bundle $(E_\rho,\nabla^{E_\rho})$ (equivalently, a local system) on $M$: more precisely, let $\widetilde{M}$ be the universal cover of $M$ on which $\pi_1(M)$ acts smoothly and freely, then 
$$E_{\rho}=\pi_{1}(M)\backslash \big(\widetilde{M}\times_{\rho} 
\mathfrak{g}\big),$$
and the flat connection $\nabla^{E_\rho}$ is induced from the usual differential on $\widetilde{M}$. 

Throughout the present paper, we make the following assumption unless otherwise stated.
\begin{assumption}
    We assume that the representation $\rho: \pi_1(M) \rightarrow G \overset{\Ad}{\longrightarrow} \Aut(\mathfrak{g})$ is acyclic, that is, the cohomology group $H^{i}(M; E_{\rho}) = 0$ for each $i =0, 1, 2, 3$.
\end{assumption}

Here, we recall examples of acyclic $\rho$ for the convenience of readers.

\begin{example}{(Fintushel--Stern (\cite{FS}, see also \cite[Lecture 15]{Sa}))}
Let $n \geq 3$ be a fixed integer and $a_1,\ldots, a_n$ be pairwise relatively prime integers with $a_i \geq 2$. Let $\Sigma=\Sigma(a_1, \ldots, a_n)$ be the Seifert homology sphere determined by $(a_1,\ldots, a_n)$ (cf. \cite[\S 6.3]{Sa}), which is the closed orientable $3$-manifold to be considered, and let $\varphi : \pi_1(\Sigma) \rightarrow \SU(2)$ be an irreducible representation. Recall that $\pi_1(\Sigma)$ admits a finite presentation as 
\begin{equation}
    \pi_1(\Sigma) = \langle x_1, \ldots, x_n ,h \mid [h,x_i]=1 , x_i^{a_i} h^{b_i} = 1, x_1\cdots x_n =1 \rangle
\end{equation}
where $b_1, \ldots, b_n$ are integers satisfying  the equation
\begin{equation}
    a_1 \cdots a_n \cdot \sum_{i=1}^n \frac{b_i}{a_i} = 1.
\end{equation}
Suppose that, for an integer $m \geq 3$, $\varphi(x_k) \neq \pm 1$ for $k=1,\ldots, m$ and $\varphi(x_k) = \pm 1$ for $k=m+1,\ldots, n$. Then, it turns out that $H^1(\Sigma; E_{\rho})= \mathbb{R}^{2m-6}$ with $\rho = \Ad \circ \varphi.$ In particular, when $m=3$, one obtains $H^1(\Sigma; E_{\rho})=0$. Since $\varphi$ is irreducible, we have $H^0(\Sigma; E_{\rho})=0$. Therefore, using Poincar\'e duality, we conclude $H^{\bullet}(\Sigma; E_{\rho})=0$ and $\rho$ with $m=3$ gives an acyclic local system on $\Sigma$.
\end{example}

\subsubsection{Lie bracket on {$E_\rho$}}
We denote the Lie bracket operator by $\mathfrak{L}=[\cdot, \cdot]:\mathfrak{g}\otimes\mathfrak{g}\rightarrow \mathfrak{g}$, $a\otimes 
b\mapsto [a,b]$. Recall that $\pi_{1}(M)$ acts on $\mathfrak{g}\otimes\mathfrak{g}$ by 
the diagonal action via $\rho$. The Lie bracket operator
$\mathfrak{L}$ is equivariant with respect to the actions of $\pi_{1}(M)$.

Moreover, the Lie bracket operator $\mathfrak{L}$
defines canonically  a morphism of vector bundles on $M$,
\begin{equation}
	\mathfrak{L}:E_{\rho}\otimes E_{\rho}\rightarrow E_{\rho}.
	\label{eq:1.5.1}
\end{equation}
By construction, it preserves the flat connections, that is, when acting 
on smooth sections,
$\nabla^{E_{\rho}}\circ \mathfrak{L}=\mathfrak{L}\circ \nabla^{E_{\rho}\otimes E_{\rho}}$.

We also extend it on $\Omega^{\bullet}(M,E_{\rho})$ such that if 
$\alpha,\beta\in \Omega^{\bullet}(M)$, $s_{1}, s_{2}\in 
C^{\infty}(M,E_{\rho})$,
\begin{equation}
	\mathfrak{L}(\alpha s_{1} \otimes \beta s_{2})=\alpha\wedge \beta 
	\mathfrak{L}(s_{1}\otimes s_{2}).
	\label{eq:1.5.3}
\end{equation}
Then $\mathfrak{L}$ induces the morphism of de Rham cohomology groups,
\begin{equation}
	\mathfrak{L}: H^{\bullet}(M; E_{\rho}\otimes E_{\rho})\rightarrow 
	H^{\bullet}(M; E_{\rho}).
	\label{eq:1.5.4}
\end{equation}

\subsubsection{Killing form and cubic trace form}\label{ss3.3Sept}
Recall that the Killing form $B \in \mathfrak{g}^{\ast} \otimes \mathfrak{g}^{\ast}$ is a 
non-degenerate bilinear form, we have the corresponding Casimir 
element $\mathbf{1} \in \mathfrak{g} \otimes \mathfrak{g}$. Let $e_1, 
\ldots, e_{\dim \mathfrak{g}}$ be a basis of $\mathfrak{g}$, and let 
$e_1^{\ast}, \ldots, e_{\dim\mathfrak{g}}^{\ast}\in\mathfrak{g}$ be 
dual basis of $\{e_i \}$ with respect to $B$, i.e., 
$B(e_{i},e_{j}^{\ast})=\delta_{ij}$. Then $\mathbf{1}$ can be explicitly written as
\begin{equation}
	\mathbf{1} = \sum_{i=1}^{\dim \mathfrak{g}} e_i \otimes e_i^{\ast}.
	\label{eq:1.3.1}
\end{equation}

Moreover, the element $\mathbf{1}$ is fixed by the diagonal action of $\pi_1(M)$ on 
$\mathfrak{g}\otimes\mathfrak{g}$, hence we can view $\mathbf{1}$ 
as a smooth section of $E_{\rho}\otimes E_{\rho}$ on $M$. Moreover, we get a well-defined map
\begin{equation}
	I : \underline{\mathbb{R}} \rightarrow E_{\rho} \otimes E_{\rho}, \quad 1 \mapsto \mathbf{1}.
	\label{eq:1.3.3}
\end{equation}
\begin{lemma}\label{lm:1.1}
	The section $\mathbf{1}$ is a flat section.
\end{lemma}

The invariant bilinear form $B$ also induces a fiberwise non-degenerate bilinear form on $E_{\rho}\rightarrow M$. By abuse of notation, we use the same $B$ to denote it, i.e., we have
\begin{equation}
	B : E_{\rho} \otimes E_{\rho} \rightarrow \underline{\mathbb{R}},
	\label{eq:1.3.2}
\end{equation}
where $\underline{\mathbb{R}}$ stands for the trivial local system on $M$.

We now introduce the cubic trace form.  For $a,b,c\in\mathfrak{g}$, the cubic trace form
\begin{equation}
    \Tr: \mathfrak{g}^{\otimes 3} \rightarrow \bR
\end{equation}
is defined  as
\begin{equation}
	\mathrm{Tr}[a\otimes b\otimes c]=B([a,b],c) = B(a, [b,c]).
\end{equation}
Then $\mathrm{Tr}\in\Lambda^{3}(\mathfrak{g}^{\ast})$, which is  $\mathrm{Ad}(G)$-invariant. 
Hence it extends to 
$\Lambda^{3}(E_{\rho}^{\ast})$, which is flat with respect to 
$\nabla^{E_{\rho}}$ and is also denoted by
\begin{equation}
    \Tr: E_{\rho}^{\otimes 3} \rightarrow \underline{\bR}.
\end{equation}
Therefore, for smooth sections $s_1, s_2, s_3$ of $E_{\rho}$ we have
	\begin{equation}
			\mathrm{Tr}[s_{1}\otimes s_{2} \otimes 
		s_{3}]=B(\mathfrak{L}(s_{1}\otimes s_{2})\otimes s_{3})\in C^{\infty}(M).
		\label{eq:1.5.2bis}
	\end{equation}

\subsection{An isomorphism of {$H^2_{-}(M;E_{\rho} \otimes E_{\rho})$}}

Let $\mathfrak{h}\subset \mathfrak{g}\otimes\mathfrak{g}$ be the kernel space of $\mathfrak{L}$. Since $\mathfrak{g}$ is semi-simple, then $\mathfrak{L}$ is surjective. Then 
\begin{equation}
     \dim_{\bR} {\mathfrak{h}} =(\dim_{\bR} {\mathfrak{g}})^2-\dim_{\bR} {\mathfrak{g}}.
\end{equation}
In particular, the symmetric tensor space $S^2\mathfrak{g}$ is a subspace of $\mathfrak{h}$, i.e, $S^2 \mathfrak{g} \subset \mathfrak{h}$.

Set the vector bundle
\begin{equation}
	H_{\rho}=\ker(\mathfrak{L}:E_{\rho}\otimes E_{\rho}\rightarrow E_{\rho}),
	\label{eq:1.5.5}
\end{equation}
which is a subbundle of $E_\rho\otimes E_\rho$ equipped with the induced flat connection $\nabla^{H_\rho}$. In fact, the vector bundle $H_{\rho}$ can also be defined as
\begin{equation}
    H_\rho=\pi_1(M)\backslash(\widetilde{M}\times_\rho \mathfrak{h}).
\end{equation}

By our construction, we have the short exact sequence of flat vector bundles: 
\begin{equation}
0\rightarrow H_{\rho}\hookrightarrow 
E_{\rho}\otimes E_{\rho}\rightarrow E_{\rho}\rightarrow 0.
\label{eq:2.2.4-july-24}
\end{equation}
Then we obtain the following long exact sequence of de Rham cohomology groups,
\begin{equation}
	\begin{split}
		&\cdots\rightarrow H^{0}(M; E_{\rho})\rightarrow 
		H^{1}(M; H_{\rho})\rightarrow 
		H^{1}(M; E_{\rho}\otimes E_{\rho}) \overset{\mathfrak{L}}{\longrightarrow}
		H^{1}(M; E_{\rho})\\
		&\qquad\qquad\rightarrow 
		H^{2}(M; H_{\rho}) \rightarrow 
		H^{2}(M; E_{\rho}\otimes E_{\rho}) \overset{\mathfrak{L}}{\longrightarrow}
		H^{2}(M; E_{\rho})\rightarrow \cdots  
	\end{split}
	\label{eq:1.5.6}
\end{equation}

\begin{example}\label{ex:3.5.1}
Consider the flat section $\mathbf{1}$ of $E_\rho\otimes E_\rho$ over $M$, that is, $\mathbf{1}\in H^0(M; E_\rho\otimes E_\rho)$. A straightforward computation shows that $\mathfrak{L}(\mathbf{1})=0$, hence we have $\mathbf{1}\in H^0(M; H_\rho)$. Since $B$ is a symmetric bilinear form, we also conclude $T^*\mathbf{1}=\mathbf{1}$, i.e., $\mathbf{1}\in H^0_+(M,H_\rho)$.
\end{example}

Note that $T$ acts on $M$ as identity but exchanges the factors of the tensor $E_\rho\otimes E_\rho$. Induces involutions $T^\ast$ acting on the de Rham complexes and on the cohomology groups $H^\bullet(M; H_\rho)$, $H^\bullet(M; E_\rho\otimes E_\rho)$. Let $\Omega^\bullet_{\pm}(M; H_\rho)$, $\Omega^\bullet_{\pm}(M; E_\rho\otimes E_\rho)$, $H^\bullet_{\pm}(M; H_\rho)$, $H^\bullet_{\pm}(M; E_\rho\otimes E_\rho)$ denote the eigenspaces of $T^\ast$ corresponding to the eigenvalues $\pm 1$. We also make $T$ act on $E_\rho$ by $-\mathrm{Id}_{E_\rho}$ point-wisely, so that the action of $T$ preserves the short exact sequence \eqref{eq:2.2.4-july-24}. Note that $\Omega^\bullet_{\pm}(M;\bullet)=\Omega^\bullet(M;\bullet_\pm)$, then by \eqref{eq:1.5.6}, we have the following result.

\begin{proposition}\label{prop:3.5.2paris}
    \begin{enumerate}[(1)]
    \item For $\rho$ which may not be acyclic, we have the isomorphism
    \begin{equation}
        H^\bullet_+(M; H_\rho)\simeq H^\bullet_+(M;E_\rho\otimes E_\rho).
        \label{eq:3.5.9paris}
    \end{equation}
    We also have the long exact sequence as follows
    \begin{equation}
	\begin{split}
		&\cdots\rightarrow H^{0}(M; E_{\rho})\rightarrow 
		H^{1}_{-}(M; H_{\rho})\rightarrow 
		H^{1}_{-}(M; E_{\rho}\otimes E_{\rho}) \overset{\mathfrak{L}}{\longrightarrow}
		H^{1}(M; E_{\rho})\\
		&\qquad\qquad\rightarrow 
		H^{2}_{-}(M; H_{\rho}) \rightarrow 
		H^{2}_{-}(M; E_{\rho}\otimes E_{\rho}) \overset{\mathfrak{L}}{\longrightarrow} 
		H^{2}(M; E_{\rho})\rightarrow \cdots  
	\end{split}
	\label{eq:3.5.10paris}
\end{equation}
\item When $\rho$ is acyclic, we also have
\begin{equation}
	H^{\bullet}_{-}(M; H_{\rho})\simeq H^{\bullet}_{-}(M; E_{\rho}\otimes 
	E_{\rho}).
	\label{eq:1.5.7}
\end{equation}
\end{enumerate}
\end{proposition}

Therefore, for each cohomology class $[\xi]\in 
H^{\bullet}_{\pm}(M; E_{\rho}\otimes E_{\rho})$, there exists a closed form 
$\xi_{0}\in\Omega^{\bullet}_{\pm}(M; H_{\rho})$, such that $[\xi_{0}]=[\xi]$ 
and 
\begin{equation}
	\mathfrak{L}(\xi_{0})=0.
	\label{eq:1.5.8}
\end{equation}

\begin{remark}
    By Example \ref{ex:3.5.1}, we see that $\mathbf{1}\in H^0_+(M;E_\rho\otimes E_\rho)$, which means that $E_\rho\otimes E_\rho$ can never be acyclic.
\end{remark}

 The following proposition can be viewed as an extension of \cite[Lemma 4.6]{Shi}.
\begin{proposition}\label{prop:3.5.2}
If $G$ is a real $3$-dimensional simple Lie group, then we have
\begin{equation}
    H^\bullet_{-}(M; H_\rho)=0,\; H^\bullet_{-}(M;E_\rho\otimes E_\rho)\simeq H^\bullet(M;E_\rho).
\end{equation}
In particular, if in addition $\rho$ is acyclic, then 
$$H^\bullet_{-}(M; E_\rho\otimes E_\rho)=0.$$
\end{proposition}
We need to point out that if $G$ is semi-simple and real $3$-dimensional, then it has to be a simple Lie group. In fact, for such a linear Lie group $G$, if $G$ is compact, then $G=\mathrm{SU}(2)$ or $\mathrm{SO}(3)$; if $G$ is noncompact, then $G=\mathrm{SL}_2(\bR)$ or $\mathrm{SO}(2,1)$.

\begin{proof}[Proof of Proposition \ref{prop:3.5.2}]
When $\mathfrak{g}$ is simple with $\dim_{\bR} \mathfrak{g}=3$, we directly conclude $\mathfrak{h}=S^2\mathfrak{g}$. In fact, the simplicity of $\mathfrak{g}$ implies an exact sequence $0  \rightarrow \mathfrak{h}\rightarrow S^2 \mathfrak{g} \oplus \Lambda^2 \mathfrak{g} \rightarrow \mathfrak{g} \rightarrow 0$, and $\dim_{\bR} \mathfrak{g}=3$ leads to $\dim_{\bR} S^2 \mathfrak{g} = 6$ and $\dim_{\bR} \Lambda^2 \mathfrak{g} = 3$. Since $S^2\mathfrak{g} \subset \mathfrak{h}$, one obtains $\mathfrak{h}=S^2 \mathfrak{g}$ for the dimensional reason. Since $T^*$ acts on $H^\bullet(M;H_\rho)$ as identity, we get
\begin{equation}
    H^\bullet_+(M;H_\rho)=H^\bullet(M;H_\rho),\;\; H^\bullet_-(M;H_\rho)=0.
\end{equation}
Then this proposition follows from Proposition \ref{prop:3.5.2paris}.
\end{proof}

\subsection{Construction of examples for Proposition \ref{mainprop1}}

The proof of Proposition \ref{mainprop1} is built on the explicit construction of certain examples of the pair $(M, \rho)$ such that $H^{\bullet}(M; E_{\rho}) =0$, but $H^1_{-}(M; E_{\rho} \otimes E_{\rho}) \simeq H^2_{-}(M; E_{\rho} \otimes E_{\rho}) \neq 0$ when $G=\SL(2, \mathbb{C})\times \SL(2, \mathbb{C})$. So, Proposition \ref{mainprop1} is a consequence of the following proposition.

\begin{proposition}\label{prop:3.-1}
    Let $M$ be an oriented closed hyperbolic 3-manifold that contains a totally geodesic surface. Consider the representation
    \begin{equation}
        \rho:\pi_1(M) \overset{\widetilde{\mathrm{hol}}}{ \xrightarrow{\hspace*{0.5cm}}} \SL(2, \mathbb{C})  \overset{\Id \times \overline{\,\Id\,}}{ \xrightarrow{\hspace*{1cm}}}\SL(2,\mathbb{C}) \times \SL(2,\mathbb{C}) \overset{\Ad}{\longrightarrow} \Aut_{\mathbb{C}} (\mathfrak{sl}_2(\mathbb{C}) \oplus \mathfrak{sl}_2(\mathbb{C})),
    \end{equation}
    where $\pi_1(M) \overset{\widetilde{\mathrm{hol}}}{\longrightarrow} \SL(2, \mathbb{C})$ is a lift of the holonomy representation $\mathrm{hol}: \pi_1(M) \rightarrow \PSL(2, \mathbb{C})$  corresponding to the complete hyperbolic structure of $M$ and $\overline{\,\Id\,}: \SL(2, \mathbb{C}) \rightarrow \SL(2,\mathbb{C})$ denotes the complex conjugation of matrix.  Then, we have
    \begin{equation}
        H^{\bullet}(M; E_{\rho}) = 0,\quad H^1_{-}(M; H_{\rho}) \simeq H^1(M; \Lambda^2 E_{\rho}) \neq 0.
    \end{equation}
\end{proposition}

\begin{proof}
  To prove the statement, we use several facts about hyperbolic 3-manifolds, for example, summarized in Porti's paper \cite{P}. Let $V_{2,0}$ denote the space of complex homogeneous polynomials with two variables and of degree $2$. Then, $\SL(2,\mathbb{C})$ acts on  $V_{2,0}$ by $(A,P) \mapsto P \circ A^T$ for $A \in \SL(2, \mathbb{C})$ and $P \in V_{2,0}$, where $A^T$ denotes the transposition of $A$. We set $V_{0,2}:= \overline{V_{2,0}}$ the complex conjugate representation to $V_{2,0}$. Then, it is known that  $V_{2,0}, V_{0,2}$ and $V_{2,2}:= V_{2,0} \otimes V_{0,2}$ are irreducible representations of $\SL(2, \mathbb{C})$ (cf. \cite[Chapter 2.\S3]{MR0855239}).
  
  Let $\mathfrak{sl}_2(\mathbb{C})_{\Ad}$ and $\mathfrak{sl}_2( \mathbb{C})_{\overline{\Ad}}$ denote the $\mathfrak{sl}_2(\mathbb{C})$ as $\SL(2, \mathbb{C})$-modules via adjoint representation $\Ad:\SL(2, \mathbb{C})\rightarrow \Aut(\mathfrak{sl}_2(\mathbb{C}))$ and its complex conjugate $\overline{\Ad}$ respectively. Then, as $\SL(2, \mathbb{C})$-modules we have $\mathfrak{sl}_2(\mathbb{C})_{\Ad} \simeq V_{2,0}$ and $\mathfrak{sl}_2(\mathbb{C})_{\overline{\Ad}} \simeq V_{0,2}$. Hence,by Raghunathan vanishing theorem (\cite[Theorem 5.1]{P}, \cite{R1}), the irreducibility of the representations $V_{2,0}$ and $V_{0,2}$, and the Poincar\'{e} duality, we conclude that $H^{\bullet}(M; E_{\rho}) =H^{\bullet}(M; V_{2,0}) \oplus H^{\bullet}(M; V_{0,2})= 0$. 
  
  Here, we denote the local systems associated with $\pi_1(M) \overset{\widetilde{\mathrm{hol}}}{ \longrightarrow} \SL(2, \mathbb{C}) \rightarrow \GL(V)$ for $V=V_{2,0}, V_{0,2}$ by the same symbols $V_{2,0}$ and $V_{0,2}$ respectively. We have  isomorphisms as $\SL(2, \mathbb{C})$-modules 
\begin{equation}
\begin{split}
    & \Lambda^2(\mathfrak{sl}_2(\mathbb{C})_{\Ad} \oplus \mathfrak{sl}_2(\mathbb{C})_{\overline{\Ad}}) \\
    & \simeq \Lambda^2( V_{2,0} \oplus V_{0,2}) \\
    & \simeq \Lambda^2(V_{2,0}) \oplus (V_{2,0} \otimes V_{0,2}) \oplus \Lambda^2 (V_{0,2})\\
    & \simeq V_{2,0} \oplus V_{2,2} \oplus V_{0,2},
\end{split}
\end{equation}
where we refer, for example, \cite[Excercise 11.35]{FH} for the last isomorphism. Again by Raghunathan vanishing theorem, we get
\begin{equation}
    H^1(M; \Lambda^2 E_{\rho}) \simeq H^1(M; V_{2,2}).
\end{equation}
Then, by Millson's theorem (\cite[Proposition 5.4]{P}, \cite{Mi}), under the assumption that $M$ contains a totally geodesic surface,  we conclude $H^1(M; V_{2,2} ) \neq 0$. The claim is proved.
\end{proof}

\begin{remark}
    It is known that there are infinitely many hyperbolic rational homology $3$-spheres containing closed embedded totally geodesic surfaces \cite[Theorem 2]{DB}. Therefore, by combining this fact with Proposition \ref{prop:3.-1}, one sees that there are infinitely many examples of a pair $(M, \rho)$ which satisfies the condition  $H^{\bullet}(M; E_{\rho}) =0$ and $H^1_{-}(M; E_{\rho} \otimes E_{\rho}) = H^2_{-}(M; E_{\rho} \otimes E_{\rho}) \neq 0$.
\end{remark}

\section{Propagators on {$C_2(M)$} for acyclic local systems}\label{section4}
This section first recalls the basics of propagators, from which one can define configuration space integrals, for acyclic local systems following \cite{BC2} and \cite{CS}. Then, we introduce a class of propagators, called adapted propagators.  They have a distinguished feature that is crucial for our results in subsequent sections.

Recall that $p_i : M \times M \rightarrow 
M$ denotes the $i$-th projection map ($i=1,2$). 
In the sequel, we also denote the induced smooth map by the same notation $p_i : C_2(M) \rightarrow M$. Moreover, for $n\geq 2$, the projection ($i\neq j$)
     \begin{equation}
        p_{ij}:\mathrm{Conf}_n(M)\ni (x_1,\ldots, x_n)\mapsto (x_i, x_j)\in \mathrm{Conf}_2(M)
     \end{equation}
     induces a smooth map of manifolds with corners
       \begin{equation}
      C_n(M)\rightarrow C_2(M)
     \end{equation}
     which is still denoted by $p_{ij}$. 

We define the pull-back vector bundle $F_{\rho}:= q^{\ast}(E_{\rho} 
\boxtimes E_{\rho}) \rightarrow C_2(M)$. Note that 
$F_{\rho}|_{\partial C_2(M)}$ is just the pull-back of $E_{\rho} 
\otimes E_{\rho}\rightarrow M$ by the projection $q_\partial: \partial 
C_2(M)\rightarrow\Delta\simeq M$, which we still denote by $E_{\rho} 
\otimes E_{\rho}$. Moreover, we have the corresponding induced flat connection 
$\nabla^{F_{\rho}}$ on $F_\rho$.

\subsection{An element in $H^2_{-}(\partial C_2(M); E_\rho\otimes E_\rho)$}

The involution $T$ on $M\times M$ given by $(x_1, x_2) \mapsto (x_2, x_1)$ extends to an involution on $C_{2}(M)$, which preserves the boundary $\partial C_2(M)$. It also lifts to the bundle $F_\rho$ by exchanging the factors of the tensor product. Let $\Omega^\bullet_{\pm}(C_2(M); F_\rho)$ (resp. $\Omega^\bullet_{\pm}(\partial C_2(M); E_\rho\otimes E_\rho)$) denote the $(\pm 1)$-eigenspaces of the action of $T^\ast$, and we also use similar convention for the cohomology groups.

We consider the unit oriented sphere $\mathbb{S}^2$ in $\bR^3$. Let $T_{\mathbb{S}^2}$ be the involution on $\mathbb{S}^2$ given by $T_{\mathbb{S}^2}(v)=-v, v\in \mathbb{S}^2$. Let $\eta$ denote a smooth normalized volume form on the unit sphere $\mathbb{S}^2$ such that $T^*_{\mathbb{S}^2}\eta=-\eta$. 
Consider the obvious projection $\pi: M\times\mathbb{S}^2\rightarrow \mathbb{S}^2$, then $\pi^* \eta$ is a closed $2$-form on $M\times\mathbb{S}^2$.

Note that the sphere normal bundle $S\nu_{\Delta}$ is identified with $S(TM)$ by
\begin{equation}
S\nu_{\Delta} \overset{\simeq}{\rightarrow} S(TM), \quad ((x,x), 
(-v,v)) \mapsto (x, v).
\label{eq:2.1.1}
\end{equation}
The involution $T$ on $C_2(M)$ restricting to the boundary corresponds to the 
involution on $S(TM)$: $(x, v) \mapsto (x, -v)$. We always use $T$ to 
denote all these involution operators.

The given framing $f$ of $M$ induces a canonical identification $M\times\mathbb{S}^2\simeq S(TM)\simeq \partial C_2(M)$ (see also \ref{A.3-july}). In this way, we view $\pi^* \eta$ as a closed $2$-form on $\partial C_2(M)$, such that $T^* (\pi^*\eta)=\pi^* (T^*_{\mathbb{S}^2} \eta)=-\pi^*  \eta$. Under the above consideration, we have an identification of vector spaces
\begin{equation}
    H^2_{-}(\partial C_2(M); E_\rho\otimes E_\rho)\simeq H^0_{+}(S^2; \mathbb{R}) \otimes H^2_{-}(M; E_\rho\otimes E_\rho)\oplus\big(H^2(\mathbb{S}^2;\bR)\otimes H^0_{+}(M; E_\rho\otimes E_\rho)\big),
    \label{eq:4.1.2paris}
\end{equation}
where $H^2(\mathbb{S}^2; \bR)=\bR[\eta]$ is $1$-dimensional.

Recall that the element $\mathbf{1}$ defined in \eqref{eq:1.3.1} can also be regarded as a flat section of $E_{\rho}\otimes E_{\rho}$ 
over $M$, $S(TM)$, and $\partial C_2(M)$.
We define a 2-form on $\Omega^2(\partial C_2(M); E_{\rho} \otimes E_{\rho})$ by 
\begin{equation}
	I(\eta) = \pi^*\eta \otimes \mathbf{1},
	\label{eq:2.1.2}
\end{equation}
where the notation $I(\cdot)$ is compatible with the definition given in \eqref{eq:1.3.3}, and the form $\pi^*\eta$ should be viewed as the $2$-form on $\partial C_2(M)$ given by the pull-back of $f$, sometimes we also denote it by $f^*\eta$ to emphasize the role of the framing $f$.

The following lemma is an analog of \cite[Proposition 3.1]{BC}.
\begin{lemma}\label{lm:2.1}
The 2-form $I(\eta) \in \Omega^2(\partial C_2(M); E_{\rho} \otimes E_{\rho})$ in \eqref{eq:2.1.2} satisfies the following properties:
\begin{enumerate}[(i)]
\item $I(\eta)$ is a closed form, its $($$E_{\rho} \otimes E_{\rho}$-valued$)$ fiber integration along $q_\partial$ is $\mathbf{1}$;
\item $T^{\ast}(I(\eta)) = - I(\eta)$.
\end{enumerate}
\end{lemma}

The following result was already implied in \cite[\S 5.2 Proof of Theorem 5.1]{CS}
\begin{lemma}\label{lem:4.1.2paris}
    If $f'$ is another smooth framing of $M$ which induces the same orientation as $f$ does, then by taking an oriented normalized volume form $\eta'$ on $\mathbb{S}^2$, the corresponding closed $2$-form $I(\eta')$ lies in the same de Rham cohomology class $[I(\eta')]$ as of $I(\eta)$ in $H^2_{-}(\partial C_2(M);E_\rho\otimes E_\rho)$.
\end{lemma}
\begin{proof}
Note that $H^2_{-}(\partial C_2(M);\mathbb{R})\simeq H^2_{-}(M\times \mathbb{S}^2; \bR)\simeq H^2(\mathbb{S}^2;\bR)$ is $1$-dimensional vector space spanned by $[\eta]$. Hence $[f^*\eta]=[(f')^*\eta']\in H^2_{-}(\partial C_2(M);\bR)$. Then there exists $\beta'\in \Omega^1_{-}(\partial C_2(M); \bR)$ such that
\begin{equation}
    f^*\eta - (f')^*\eta' = d\beta'.
    \label{eq:4.1.4paris}
\end{equation}
As a consequence, we conclude the identity in $\Omega^2_{-}(\partial C_2(M); E_\rho\otimes E_\rho)$,
\begin{equation}
    I(\eta) - I(\eta') = (d\beta')\otimes \mathbf{1}=d(\beta'\otimes \mathbf{1}).
\end{equation}
This way, we conclude this lemma.
\end{proof}

In the above lemma, the framing $f'$ is not necessary to be homotopic to $f$. When $f'$ is homotopic to $f$, the form $\beta'$ in \eqref{eq:4.1.4paris} can be constructed more explicitly as follows.
\begin{lemma}\label{lm:4.1.3July}
      Fix an oriented normalized volume form $\eta$ on $\mathbb{S}^2$. If $f'$ is another smooth framing of $M$ which is homotopic to $f$, let $I'(\eta)$ be the closed $2$-form in $\Omega^2_{-}(\partial C_2(M); E_\rho\otimes E_{\rho})$ defined by $f'$. Then there is a $f$-vertical $1$-form $\beta'\in \Omega^1_{-}(\partial C_2(M);\bR)$ (see \eqref{eq:A.3.1}) such that
      \begin{equation}
    I'(\eta) - I(\eta) =d(\beta'\otimes \mathbf{1}).
\end{equation}
\end{lemma}
\begin{proof}
    Note that in this case, $f'\circ f^{-1}$ is connected to the identity section by a smooth path in $\mathscr{C}^\infty(M,\mathrm{Diff}(\mathbb{S}^2))$. Let $\psi_\cdot: [0,1]\ni t\mapsto \psi_t \in\mathscr{C}^\infty(M,\mathrm{Diff}(\mathbb{S}^2))$ denote such a path with $\psi_0(x)=\mathrm{Id}_{\mathbb{S}^2}$ and $\psi_1(x)=(f'\circ f^{-1})_x$. In particular, we view $\psi_t$ as a diffeomorphism of $M \times \mathbb{S}^2$. Let $\Gamma(\mathbb{S}^2)$ denote the space of smooth vector fields on $\mathbb{S}^2$. Set $X_t=\frac{\partial}{\partial t}\psi_t \in \mathscr{C}^\infty(M,\Gamma(\mathbb{S}^2))$. Fix an oriented normalized volume form $\eta$ on $\mathbb{S}^2$, viewed as a constant form on $M\times \mathbb{S}^2$, then
\begin{equation}
    \frac{\partial}{\partial t}\psi_t^\ast \eta = d\iota_{X_t}\psi_t^\ast \eta,
\end{equation}
where $\iota_{X_t}$ denotes the contraction of vector fields $X_t$.
A direct computation shows that
\begin{equation}
\psi_1^\ast\eta-\eta=d\int_0^1 \iota_{X_t}\psi_t^\ast \eta dt=:d\beta.
\end{equation}
Then
\begin{equation}
(f')^\ast\eta-f^\ast\eta=d \big(f^*\beta\big).
\end{equation}

Note that $\beta$ is a vertical $1$-form on $M\times\mathbb{S}^2$, hence $\beta':=f^*\beta$ is $f$-vertical as desired. This way, we complete the proof.
\end{proof}

Now we need to introduce a commutative diagram to better understand all the different cohomology groups $H^\bullet_{-}(\cdots;\cdots)$ that we have seen. Recall that $\mathfrak{i}_{\partial}:\partial C_{2}(M)\hookrightarrow 
C_{2}(M)$ denotes the inclusion, and $\mathfrak{i}:\Delta\hookrightarrow M\times M$ denotes the inclusion of the diagonal.

Note that we have the following short exact sequence associated to the relative de Rham complex
\begin{equation}
    \begin{split}
        &0\rightarrow \big(\Omega^\bullet(C_2(M),\partial C_2(M); F_\rho),d\big) \hookrightarrow \big(\Omega^\bullet(C_2(M); F_\rho),d\big)\\
&\hspace{7cm}\overset{\mathfrak{i}_\partial^\ast}{\longrightarrow} \big(\Omega^\bullet(\partial C_2(M); E_\rho\otimes E_\rho),d\big)\rightarrow 0.
    \end{split}
\end{equation}
Together with the excision theorem for the pairs $(C_2(M),\partial C_2(M))$ and $(M\times M, \Delta)$, we get the following commutative diagram where the horizontal lines are exact sequences ($j=0,1, \ldots, 6$)
\begin{equation}
\adjustbox{scale=0.8,center}{
    \begin{tikzcd}
H^j_-(C_2(M); \mathbf{F}_\rho) \arrow[r,"\mathfrak{i}_\partial^\ast"]  &
  H^j_-(\partial C_2(M); \mathbf{F}_\rho)  \arrow[r,"\delta^*_{C_2(M)}"] & H^{j+1}_-(C_2(M), \partial C_2(M); \mathbf{F}_\rho) \arrow[r, "\mathrm{incl}"]  &
  H^{j+1}_-(C_2(M); \mathbf{F}_\rho)  \\
H^j_-(M\times M; \mathbf{F}_\rho) \arrow[r,"\mathfrak{i}^\ast"] \arrow[u,"q^\ast"] & H^j_-(\Delta; \mathbf{F}_\rho) \arrow[u,"q_\partial^\ast"]  \arrow[r,"\delta_{M^2}^\ast"] & H^{j+1}_-(M\times M, \Delta; \mathbf{F}_\rho) \arrow[r,"\mathrm{incl}"] \arrow[u,"q^\ast", "\simeq"'] & H^{j+1}_-(M\times M; \mathbf{F}_\rho) \arrow[u,"q^\ast"]
\end{tikzcd}
}
\label{eq:4.1.7paris}
\end{equation}
where the bundle $\mathbf{F}_\rho$ represents $F_\rho$ on $C_2(M)$, $E_\rho \otimes E_\rho$ on $\Delta$ or $\partial C_2(M)$, $E_\rho\boxtimes E_\rho$ on $M\times M$ respectively. The maps $\delta^*_{C_2(M)}$ and $\delta^\ast_{M^2}$ are connecting homomorphisms. Note that, when $E_\rho$ is assumed to be acyclic,  $\delta^\ast_{M^2}$ is an isomorphism.

\subsection{Propagators associated with acyclic local systems}\label{ss4.3sept}
Now we recall the definition of propagators and their existence result from \cite{CS}.

\begin{definition}[Propagator]\label{defn:prop1}
    Assume $\rho$ to be acyclic. A smooth 2-form
\begin{equation}
	\omega \in \Omega^2(C_2(M); F_{\rho})
\label{eq:6.1.1}
\end{equation}
is called a \textit{propagator} if $\omega$ satisfies the following three properties:
\begin{enumerate}[(i)]
\item $d \omega=0$;
 \item the restriction of $\omega$ on the boundary $\omega|_{\partial C_2(M)}$ has the form
 \begin{equation}
 	\mathfrak{i}_\partial^{*}(\omega) = I(\eta) + q_\partial^{\ast}(\xi)
  \label{eq:6.1.2}
 \end{equation}
 where $\eta$ is a normalized (oriented) volume $2$-form on $\mathbb{S}^2$, and $\xi \in \Omega^2_{-}(\Delta, E_{\rho} \otimes E_{\rho})$;
\item the form $\omega$ on $C_2(M)$ is antisymmetric under the action of $T^{\ast}$, that is, $T^{\ast}(\omega) = - \omega$.
\end{enumerate}
\end{definition}

\begin{proposition}\label{cor:2.2.1}
	Assume $\rho$ to be acyclic. Then there always exists a propagator $\omega \in \Omega^2(C_2(M); F_{\rho})$ and a form $\xi\in\Omega^{2}(\Delta; E_{\rho}\otimes E_{\rho})$, such that 
	$T^{\ast}(\xi)=-\xi$ and
	\begin{equation}
		\mathfrak{i}_{\partial}^{\ast}(\omega)=I(\eta)+q_\partial^{\ast}(\xi).
		\label{eq:2.57}
	\end{equation}
 Moreover, $\xi$ is closed and the class $[\xi]\in H^2_-(\Delta; E_\rho\otimes E_\rho)$ is independent of the choice of $\xi$ or the oriented framing $f$ (which is compatible with the given $o(M)$).
\end{proposition}
\begin{proof}
    The first part was done in the proof of \cite[Proposition 2.1]{CS}. For the second part follows from an easy modification of same proof. By the definition $\delta^\ast_{C_2(M)}$ in \eqref{eq:4.1.7paris}, we have
    \begin{equation}
       \delta^\ast_{C_2(M)} [ I(\eta)+q_\partial^{\ast}(\xi)]=0.
    \end{equation}
    In this case, $\delta^\ast_{M^2}$ is an isomorphism, set
$$\Phi:=(\delta^\ast_{M^2})^{-1}\circ (q^\ast)^{-1}\circ\delta^\ast_{C_2(M)}: H^2_-(\partial C_2(M); E_\rho\otimes E_\rho)\rightarrow H^2_{-}(\Delta; E_\rho\otimes E_\rho).$$
Then we have $\Phi\circ q^\ast_\partial=\mathrm{Id}_{H^2_{-}(\Delta; E_\rho\otimes E_\rho)}$, and
\begin{equation}
[\xi]=-\Phi[I(\eta)]\in H^2_{-}(\Delta; E_\rho\otimes E_\rho).
\end{equation}
Finally, by Lemma \ref{lem:4.1.2paris}, $[I(\eta)]$ is independent of the choice of oriented framing $f$.
\end{proof}
In order to emphasize the boundary condition \eqref{eq:2.57}, we use the pair $(\omega, \xi)$ or the triplet $(\omega,\eta,\xi)$ to denote our propagator. When we also want to emphasize the role of the framing $f$, we sometimes use $(\omega,f,\eta,\xi)$ to denote a propagator.

\begin{remark}
    The cohomology class $[\xi] \in H_{-}^2(\Delta; E_{\rho} \otimes E_{\rho})$ in Proposition \ref{cor:2.2.1} is the Poincar\'{e} dual of (the anti-symmetric part of) the invariant $d(M, \rho)$ studied by Kitano--Shimizu in \cite{KS}. 
\end{remark}
Note that the propagators for an acyclic $E_\rho$ are generally not unique, but the cohomology class of propagators is unique. 
\begin{proposition}[Uniqueness of propagators for acyclic $E_\rho$]\label{prop:uniqueness}
   Fix a homotopy class of the smooth framings $[f]$ of $M$ and an orientation $o(M)$, we also assume $E_\rho$ to be acyclic, and let $\omega\in \Omega^2_-(C_2(M); F_\rho)$ be a propagator in Definition \ref{defn:prop1}. Then the de Rham cohomology class $[\omega]\in H^2_-(C_2(M); F_\rho)$ is unique (which is independent of the choice of a framing $f\in [f]$ but depends on the homotopy class $[f]$).
\end{proposition}
\begin{proof}
    We consider the diagram \eqref{eq:4.1.7paris} but for the cohomology groups of degrees $1$ and $2$. Note that the map $\delta^*_{M^2}: H^1_-(\Delta;E_\rho\otimes E_\rho)\rightarrow H^2_-(M\times M,\Delta; E_\rho\boxtimes E_\rho)$ is an isomorphism. 
    
    Meanwhile, we have the isomorphism: 
    $$q^*: H^2_-(M\times M,\Delta; E_\rho\boxtimes E_\rho)\rightarrow H^2_-(C_2(M),\partial C_2(M); F_\rho).$$ 
    As a consequence, $\delta^*_{C_2(M)}: H^1_-(\partial C_2(M);E_\rho\otimes E_\rho)\rightarrow H^2_-(C_2(M),\partial C_2(M); F_\rho)$ is surjective. Therefore, we conclude that the restriction map 
    \begin{equation}
    \mathfrak{i}_\partial^\ast: H^2_-(C_2(M); F_\rho) \rightarrow
  H^2_-(\partial C_2(M); E_\rho\otimes E_\rho)
  \label{eq:4.3.30kk}
  \end{equation}
  is injective.

  Note that $H^{2}(I \times \mathbb{S}^2; \mathbb{R}) \simeq H^0(I; \mathbb{R}) \otimes H^2(\mathbb{S}^2; \mathbb{R}) \simeq H^2(\mathbb{S}^2; \mathbb{R})$, then by Proposition \ref{cor:2.2.1}, for any propagator $\omega$ defined with a framing $f\in [f]$, the cohomology class $\mathfrak{i}_\partial^\ast[\omega]$ is uniquely determined by $(M,o(M),[f],\rho)$. Finally, the uniqueness of $[\omega]\in H^2_-(C_2(M); F_\rho)$ follows from the injectivity of $\mathfrak{i}_\partial^\ast$ in \eqref{eq:4.3.30kk}. 
\end{proof}

\subsection{Adapted propagators for acyclic local systems}\label{subsection4.3}
Based on our consideration in \eqref{eq:1.5.7}, together with Proposition \ref{cor:2.2.1}, for an acyclic local system $E_\rho$, we can define a propagator $\omega$ which has an extra property with respect to the Lie bracket operator $\mathfrak{L}$, which we call an adapted propagator.

One motivation for definition is to construct the integral invariants for the triplet $(M,f,\rho)$ associated to trivalent graphs without self-loops. Note that this definition is cohomologically canonical and in the spirit of Bott--Cattaneo \cite{BC2}.

\begin{definition}[Adapted propagator]\label{def:6.1.1}
	For an acyclic local system $E_{\rho}$ associated with a representation $\rho : \pi_1(M) \rightarrow G \overset{\Ad}{\longrightarrow} \Aut(\mathfrak{g})$, 
a propagator $\omega$ is called an {\it adapted propagator} if it satisfies the following condition (iv) in addition to the above (i), (ii), (iii) of Definition \ref{defn:prop1}:
\begin{enumerate}[(iv)]
    \item $\mathfrak{L}(\xi) = 0$, or equivalently, $\xi\in\Omega^2_-(\Delta; H_\rho)$. It is also equivalent to $\mathfrak{L}(\mathfrak{i}_\partial^*(\omega))=0$. 
\end{enumerate}
\end{definition}

 Our main result for this subsection is as follows.
\begin{theorem}\label{thm:enhanced}
    Given a framing $f$, an oriented normalized volume form $\eta$ on $\mathbb{S}^2$ and an acyclic local system $E_{\rho}$ via a representation $\rho : \pi_1(M) \rightarrow G \overset{\Ad}{\rightarrow} \Aut(\mathfrak{g})$, the adapted propagator $\omega\in \Omega^2(C_2(M); F_\rho)$ always exists with the boundary condition \eqref{eq:6.1.2}.
\end{theorem}
\begin{proof}
    Let $\omega'$ be a propagator constructed as in Proposition \ref{cor:2.2.1} which also satisfies \eqref{eq:2.57} with a closed $2$-form $\xi'\in \Omega^2_{-}(\Delta; E_\rho\otimes E_\rho)$, i.e.,
    $$\mathfrak{i}_\partial^\ast (\omega')=I(\eta)+q^\ast_\partial(\xi').$$

 By Proposition \ref{cor:2.2.1}, we have $[\xi']\in H^2_{-}(\Delta; E_\rho\otimes E_\rho)$. In the same time, by \eqref{eq:1.5.7}, there eixists $\xi_0\in \Omega^2_-(\Delta; H_\rho)$ (i.e., $\mathfrak{L}(\xi_0)=0$) such that 
 \begin{equation}
     \xi'-\xi_0=d\psi,
 \end{equation}
 where $\psi\in \Omega^1_-(\Delta; E_\rho\otimes E_\rho)$.

We now extend $\psi$ to a smooth 1-form $\widetilde{\psi} \in \Omega^1(M \times M; E_{\rho} \boxtimes E_{\rho})$ as follows:  Let $N \supset \Delta$ be a tubular neighbourhood of the diagonal $\Delta$ in $M \times M$  which is invariant by $T$ with a projection $p : N \rightarrow \Delta$. Note that, by definition, we have 
\begin{equation}
	\begin{split}
p^{\ast}(E_{\rho}\boxtimes E_{\rho}|_{\Delta})_{(x,y)} &= (E_{\rho}\boxtimes E_{\rho}|_{\Delta})_{p(x,y)}\\
&= (E_{\rho}\boxtimes E_{\rho}|_{\Delta})_{(x_0, x_0)} \quad (\text{here, we  set}\ p(x,y)=(x_0, x_0))\\
& = E_{\rho, x_0} \otimes E_{\rho, x_0}.
	\end{split}
	\label{eq:2.2.4}
\end{equation}
Since $E_{\rho}$ is equipped with a flat connection $\nabla^{E_\rho}$, we have a linear isomorphism
\begin{equation}
	\tau_{(x,y)} : (E_{\rho}\boxtimes E_{\rho}|_N)_{(x,y)} \overset{\sim}{\longrightarrow} E_{\rho, x_0} \otimes E_{\rho, x_0} = p^{\ast}(E_{\rho}\boxtimes E_{\rho}|_{\Delta})_{(x,y)}
	\label{eq:2.2.5}
\end{equation}
given by parallel transport along the projection $p$ with respect to 
$\nabla^{E_\rho}$ . This way, we get an isomorphism of (flat) vector 
bundles on $N$, 
\begin{equation}
	\tau : E_{\rho}\boxtimes E_{\rho}|_N \overset{\sim}{\rightarrow} 
	p^{\ast}(E_{\rho}\boxtimes E_{\rho}|_{\Delta}).
	\label{eq:2.2.6}
\end{equation}
Thus, we get a 1-form $\tau^{-1} \psi \in \Omega^2(N, E_{\rho}\boxtimes E_{\rho}|_N)$. Let us consider a smooth cutoff function $\chi : C_2(M) \rightarrow \mathbb{R}$ such that $\chi|_{C_2(M) \setminus N} \equiv 0$, $\chi|_U \equiv 1$,and $T^{\ast} \chi = \chi$ where $U$ is an open set such that $\partial C_2(M) \subset U 
\subsetneq 
\mathrm{Inner}(N) \subset C_2(M)$ and $U$ is $T$-invariant. Now, we can define a 1-form $\widetilde{\psi}$ on $M \times M$ as follows: 
\begin{equation}
    \widetilde{\psi}=\chi \cdot \tau^{-1} \psi \in \Omega^1_{-}(M \times M; E_{\rho}\boxtimes E_{\rho}).
\end{equation}
Note that the 1-form $\widetilde{\psi}\in \Omega^{1}_{-}(M\times M, 
	E_{\rho}\boxtimes E_{\rho})$ which is supported near $\Delta$ and 
	satisfies
	\begin{equation}
		\mathfrak{i}^{\ast}(\widetilde{\psi})=\psi.
		\label{eq:3.0.10}
	\end{equation}

 Set $\omega=\omega'-dq^\ast\widetilde{\psi}\in \Omega^2_-(C_2(M);F_{\rho})$. Then
 \begin{equation}
 \mathfrak{i}^\ast_\partial \omega= \mathfrak{i}^\ast_\partial \omega'-q^\ast_\partial d\psi=I(\eta)+q^\ast_\partial (\xi'-d\psi)=I(\eta)+q^\ast_\partial (\xi_0).
 \end{equation}
 This closed $2$-form $\omega$ is an adapted propagator as we defined.
\end{proof}

A modification of the above proof gives the following statement, which corresponds to the main framework in \cite{BC,BC2}.
\begin{proposition}
     Given a framing $f$ and an oriented normalized volume form $\eta$ on $\mathbb{S}^2$. Assume that $E_\rho$ is acyclic and
     \begin{equation}
     H^1_-(M; E_\rho\otimes E_\rho)=0,
         \label{eq:4.3.6July}
     \end{equation}
     then there is an adapted propagator $\omega\in \Omega^2_-(C_2(M);F_\rho)$ such that $d\omega=0$ and
     \begin{equation}
         \mathfrak{i}^\ast_\partial (\omega)=I(\eta).
         \label{eq:4.3.6paris}
     \end{equation}
\end{proposition}

\begin{remark}
Together with Proposition \ref{prop:3.5.2}, if $G$ is a real $3$-dimensional simple Lie group and $E_\rho$ is acyclic, then the conditions in the above proposition are always satisfied. In particular, the above results apply to the cases $G=\mathrm{SU}(2)$ or $\mathrm{SL}_2(\bR)$. In general, as we saw in Proposition \ref{prop:3.-1}, there are examples of triples of $(M, G, \rho)$ with $E_\rho$ being acyclic but $H^1_{-}(M; E_{\rho} \otimes E_{\rho}) \neq 0$.
\end{remark}

\subsection{Cattaneo--Shimizu's result on 2-loop invariant}

In this subsection, we revisit the 2-loop invariants of framed closed 3-manifolds equipped with acyclic local systems introduced in \cite{CS}. After recalling its definition and the result of \cite{CS}, we observe that the choice of adapted propagators gives a vanishing of integration associated with the dumbbell graph.

On $C_{2}(M)$, we have the canonical identification of flat vector 
bundles
\begin{equation}
	F^{\otimes 3}_{\rho}=q^{\ast}(E_{\rho}^{\otimes 3}\boxtimes 
	E_{\rho}^{\otimes 3}).
	\label{eq:3.0.1}
\end{equation}
Then the cubic trace form $\Tr: E_{\rho}^{\otimes 3} \rightarrow \underline{\bR}$ induces the linear morphism $\mathrm{Tr}^{\boxtimes 2}: F^{\otimes 3}_{\rho}\rightarrow \underline{\mathbb{R}}(= q^{\ast}(\underline{\mathbb{R}} \boxtimes \underline{\mathbb{R}}))$.

As in \cite{BC, BC2} and in \cite{CS}, we now consider the integral invariants for the $2$-loops terms in Chern--Simons perturbation theory. The theta graph and the dumbbell graph are the only two connected topological trivalent graphs with $2$-loop, see Fig. \ref{fig:6.1.theta}. For each graph, we can define a configuration space integral as our potential invariant. An integral invariant $Z_1(M,\rho)$ introduced by Cattaneo--Shimizu in \cite{CS} is given as a linear combination of theta-invariant and dumbbell invariant. We always fix a framing $f$ and an orientation $o(M)$ of $M$.
\begin{definition}\label{defn:3.4.1-24}
Fix an acyclic local system $E_{\rho}$ on $M$ associated with a representation $\rho : \pi_1(M) \rightarrow G \overset{\Ad}{\longrightarrow} \Aut(\mathfrak{g})$. Provided a propagator $\omega$ as in Definition \ref{defn:prop1}: in particular, there exits $\xi\in\Omega^{2}_{-}(\Delta,E_{\rho}\otimes E_{\rho})$ such that
	\begin{equation}
		\mathfrak{i}_{\partial}^{\ast}(\omega)=I(\eta)+q^{\ast}_\partial(\xi).
		\label{eq:3.0.2}
	\end{equation}
	We define the following integrals,
	\begin{equation}
		Z_{\Theta}(\omega)=\int_{C_{2}(M)}\mathrm{Tr}^{\boxtimes 
		2}[\omega^{3}],\quad  \Zdam(\omega,\xi)=\int_{C_{2}(M)}\mathrm{Tr}^{\boxtimes 
		2}[(p_{1}^{\ast}\xi)(p_{2}^{\ast}\xi)\omega],
  \label{eq:3.0.3}
	\end{equation}
	and set
	\begin{equation}
		Z_{1}(\rho;\omega,\xi)= Z_{\Theta}(\omega) - \frac{3}{2} 
			\Zdam(\omega,\xi).
			\label{eq:3.0.4}
	\end{equation}
\end{definition}

Note that if we take an arbitrary propagator $\omega$ given as in Proposition 
	\ref{cor:2.2.1}, the theta term $Z_\Theta(\omega)$ will depend on the choice of $\omega$. Cattaneo and Shimizu \cite{CS} introduced a correction term---the dumbbell term $\Zdam(\omega,\xi)$---to finally obtain a $2$-loop integral invariant $Z_1(M,\rho)= Z_{1}(\rho;\omega,\xi)$. The precise statement is as follows.

\begin{theorem}[{\cite[Theorem 2.3]{CS}}]\label{prop:3.1}
	Given the homotopy class of framing $[f]$ and the orientation $o(M)$ for $M$. If $E_\rho$ is acyclic, then $Z_1(\rho;\omega,\xi)$ is independent of the choice of the 
	triplet 
	$(\omega,\eta,\xi)$, so that it is an invariant for $(M, E_{\rho}, [f])$, which is denoted by $Z_{1}(M,\rho)$. 
 
\end{theorem}

 \begin{remark}
     Note that we put the coefficient $\frac{3}{2}$ instead of $3$ (the coefficient originally given in \cite{CS}) in front of the dumbbell term, this difference follows from our convention of the computations (comparing \eqref{eq:self-loop_jacobi} with \cite[\S 4.2. Proof of Proposition 4.2]{CS}), more details are referred to Examples \ref{ex:6.2.1} \& \ref{exm:6.6.7sss}.
 \end{remark}

Using our construction of an adapted propagator, we can refine Cattaneo--Shimizu's result (Theorem \ref{prop:3.1}) as follows.
\begin{theorem}\label{prop:3.2}
	Assume $E_\rho$ to be acyclic. Let $\omega^\sharp\in \Omega^2_-(C_2(M); F_\rho)$ be an adapted propagator with $	\mathfrak{i}_\partial^{*}(\omega^\sharp) = I(\eta) + q_\partial^{\ast}(\xi^\sharp)$, then
	\begin{equation}
		\Zdam(\omega^\sharp,\xi^\sharp)=0,
	\end{equation}
 and therefore
\begin{equation}
    Z_{1}(M,\rho)=Z_{\Theta}(\omega^\sharp).
\end{equation}
Equivalently, $Z_{\Theta}(\omega^\sharp)$ itself  gives the 2-loop invariant for $(M,\rho, [f])$.
\end{theorem}
\begin{proof}
We only need to prove that
	\begin{equation}
		\Zdam(\omega^\sharp,\xi^\sharp)=0.
		\label{eq:3.0.13}
	\end{equation}
 Note that for an adapted propagator, we have $\mathfrak{L}(\xi^\sharp)=0$.

By \eqref{eq:1.5.2bis}, \eqref{eq:3.0.1}, we have
	\begin{equation}
		\begin{split}
			&\mathrm{Tr}^{\boxtimes 
		2}[(p_{1}^{\ast}\xi^\sharp)(p_{2}^{\ast}\xi^\sharp)\omega^\sharp]\\
		&=B_{1,2}(\mathfrak{L}(\xi^\sharp)\boxtimes 
		\mathfrak{L}(\xi^\sharp),\omega^\sharp)\\
		&=0.
		\end{split}
		\label{eq:3.0.14}
	\end{equation}
	This implies exactly \eqref{eq:3.0.13}. The proof is then completed.	
\end{proof}

The result of Theorem \ref{prop:3.2} shows that, for an acyclic local system $E_\rho$, the use of an adapted propagator defined in Definition \ref{def:6.1.1} can reduce the computation of $Z_1(M,\rho)$ to compute only the theta term, hence $Z_1(M,\rho)$ is essentially the theta-invariant. Note that the dumbbell term corresponds to the dumbbell graph, which is the only connected $2$-loop trivalent graph with self-loops, the proof of the above proposition indicates that the extra condition $\mathfrak{L}(\xi^\sharp)=0$ in Definition \ref{def:6.1.1} is the key point to vanish the self-loops. This idea will be exploited further for the integral invariants associated with higher-loop terms in subsequent sections.


\section{Graph complex associated to acyclic adjoint local systems}\label{section:graph}
This section introduces a graph complex associated with an acyclic local system which corresponds to $\rho: \pi_1(M) \rightarrow G \overset{\Ad}{\longrightarrow} \Aut(\mathfrak{g})$. The construction is an analogous version of the one defined by Bott--Cattaneo \cite{BC2} specialized so that $\rho$ is given as above and the equivariant homomorphisms associated with (internal) vertices are defined from $\Tr$. Unlike Bott--Cattaneo \cite{BC2}, we include the graphs with self-loops in our graph complex. In this section, only $\mathfrak{g}$ is involved, information from $M$ or $\rho$ is not needed.

\subsection{Preliminary on graphs}
Here we always consider the finite graph (that is, with finite number of vertices and edges).
\begin{definition}
\begin{enumerate}[(1)]
\item A {\it self-loop} of a graph is an edge that connects the same vertex.
\item If two distinct vertices of a graph are connected by exactly one edge, then this edge is said to be {\it regular}. A graph is said to be connected if it is path connected (every two vertices can be connected by a path of edges). 
\item Let $\Gamma$ be a graph whose edges are directed. Let $v(\Gamma)$ and $e(\Gamma)$ denote the sets of vertices and edges of $\Gamma$ respectively. For an directed edge $e$ of $\Gamma$ connecting the vertex $i$ to $j$, we define a map $s : e(\Gamma) \rightarrow v(\Gamma)$ and $t : e(\Gamma) \rightarrow v(\Gamma)$ by $s(e) = i$ and $t(e) = j$. Then, a {\it half edge} of a graph $\Gamma$ is defined as an element of the form
\begin{equation*}
	(s(e), e, +1) \mathrm{\;or\;} (t(e), e, -1) \in v(\Gamma) \times e(\Gamma) \times \{\pm 1\}
\end{equation*}
for $e \in e(\Gamma)$. We call the number of half-edges at a vertex $i$ {\it valency} of the vertex $i$. Usually, we use $h(\Gamma)$ to denote the set of all half-edges of $\Gamma$.
\item A graph $\Gamma$ is said to be {\it trivalent} (resp. {\it uni-trivalent}) if the valencies for vertices all are $3$ (resp. $1$ or $3$).
\item A univalent vertex of a graph $\Gamma$ is  called {\it external vertex} and a vertex with valency $\geq 2$ of $\Gamma$ is called {\it internal vertex}. Similarly, an edge of $\Gamma$ which connects two internal vertices is called {\it internal edge} and called {\it external edge} otherwise.
\end{enumerate}
\end{definition}

\begin{figure}[h]
 \centering
 \begin{tikzpicture}
    \draw[densely dotted] (0,0) circle (0.5cm);
    \draw[densely dotted] (-1,-0.5) -- (-0.5,0);
    \draw[densely dotted] (-1,0.5) -- (-0.5,0);
     \node at (-0.5, 0) {$\bullet$};
      \node at (-1.5, 0.1) {$\vdots$};
\end{tikzpicture}
 \caption[A self-loop edge]{A self-loop of a graph, vertex with valency $4$ as displayed}
 \label{fig:6.2.1}
\end{figure}
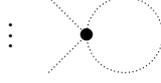

In the sequel, a connected graph always means a connected graph whose internal vertices have valency $\geq 3$. Next, we define several orientations of a connected graph.

\begin{definition}
Let $\Gamma$ be a connected uni-trivalent graph.  A vertex-wise orientation of $\Gamma$ is a collection of cyclic order of half-edges attached to each internal vertex of $\Gamma$.
\end{definition}


\begin{definition}[Vertex orientation and edge orientation of half-edges]
    Let $\Gamma$ be a connected (uni-)trivalent graph $\Gamma$. Let  $h^\mathrm{int}(\Gamma)$ be the set of internal half-edges (i.e., half-edges attached to internal vertices).  An orientation of $h^\mathrm{int}(\Gamma)$ is a numbering on $h^\mathrm{int}(\Gamma)$ up to even permutations, i.e., a bijection $h^\mathrm{int}(\Gamma) \simeq \{1,2,\ldots, |h^\mathrm{int}(\Gamma)|\}$ where two such bijections are identified if they are related by even permutations. We introduce the following two orientations for $h^\mathrm{int}(\Gamma)$:
    \begin{enumerate}[(1)]
    \item Suppose that $\Gamma$ is vertex-wise oriented and the set $v^{\mathrm{int}}(\Gamma)$ is ordered. Then, a {\it vertex orientation} of $h^\mathrm{int}(\Gamma)$ is defined as follows: according to the order of  $v^{\mathrm{int}}(\Gamma)$, take an internal vertex $v$ and order the set of half-edges at $v$.
    \item Suppose that all of edges of $\Gamma$ are directed and the set $e(\Gamma)$ is ordered. Then, an {\it edge orientation} of $h^\mathrm{int}(\Gamma)$ is defined as the induced orientation from that of $e(\Gamma)$ and directions of edges. Here, for an oriented self-loop $e$ connecting the vertex $v$, the order of two half-edges $(v, e, +1)$ and $(v, e, -1)$ is defined so that $(v, e, +1)$ is putted just before $(v, e, -1)$.
    \end{enumerate}
\end{definition}

\subsection{Weight systems associated with uni-trivalent trees}
Now let $G$ be a connected semi-simple Lie group with Lie algebra $\mathfrak{g}$ and the Killing form $B$ as considered in Section \ref{section3}. 
Here, we describe a way to obtain $\Ad(G)$-invariant multilinear map $\mathfrak{g}^{\otimes n}  \rightarrow \bR$ associated with $B$ for uni-trivalent trees, as in \cite{BN}.

Let $T$ be a vertex-wise oriented uni-trivalent tree diagram with $n$ external vertices (hence it has no loops). Suppose that the set of $n$ external vertices is ordered. Then, associated with $T$, the {\it weight system} 
\begin{equation}
    W_{T}:  \mathfrak{g}^{\otimes n} \rightarrow \bR
    \label{eq:6.2.1August}
\end{equation}
is defined as follows. For each external vertex, we associate it with $n$ inputs of elements of $\mathfrak{g}$ according to the order on the set of external vertices. For each trivalent vertex, we assign the cubic trace form $\Tr$, defined in Subsection \ref{ss:2.1-july}, according to the cyclic order of half-edges at the vertex, and for each internal edge, we assign the Casimir element $\mathbf{1}$. Then, taking contraction with respect to $B$ along internal edges, we obtain the desired multilinear form $W_{T}$. By construction, $W_T$ is $\Ad(G)$-invariant and independent of the order of contractions. 

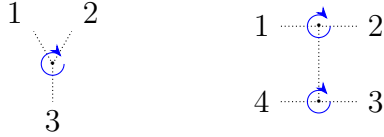
\begin{figure}[h]
\captionsetup{margin=2cm}
 \centering
 \begin{tikzpicture}
\begin{scope}
\node (1) at (-1/4-0.25,1.73/4+0.25) {1};
\node (2) at (1/4+0.25,1.73/4+0.25) {2};
\node (3) at (0,-0.75) {3};
\draw[densely dotted] (0,0) -- (1/4,1.73/4);
\draw[densely dotted] (0.0,0) -- (-1/4,1.73/4);
\draw[densely dotted] (0,0) -- (0, -0.5);
\node at (0,0) {$\cdot$};
\draw[-Stealth, color=blue]  (0.0+0.15,0.0) arc(0:-320:0.15);
\end{scope}
\begin{scope}[xshift=3cm]
\node (1) at (-0.25,0.5) {1};
\node (2) at (1.25,0.5) {2};
\node (3) at (1.25,-0.5) {3};
\node (4) at (-0.25,-0.5) {4};
\draw[densely dotted] (1.0,0.5) -- (-0,0.5);
\draw[densely dotted] (-0.0,-0.5) -- (1.0,-0.5);
\draw[densely dotted] (0.5,0.5) -- (0.5,-0.5);
\node at (0.5,0.5) {$\cdot$};
\node at (0.5,-0.5) {$\cdot$};
\draw[-Stealth, color=blue]  (0.15+0.5,-0.5) arc(0:-320:0.15);
\draw[-Stealth, color=blue]  (0.15+0.5,0.5) arc(0:-320:0.15);
\end{scope}
\end{tikzpicture}
 \caption[Examples of uni-trivalent trees]{Examples of uni-trivalent trees whose univalent vertices are ordered and trivalent vertices are cyclically ordered.}
 \label{fig:W_T}
\end{figure}

\begin{example}
We give some elementary examples of the maps $W_T$ for the case that $\mathfrak{g}=\mathfrak{sl}_2$. As in \cite[\S 10. 6]{FH}, the Lie algebra $\mathfrak{sl}_2$ is generated by 
\begin{equation}
    h = \begin{pmatrix}
        1 & 0 \\
        0 & -1
    \end{pmatrix},
    \quad e = \begin{pmatrix}
        0 & 1\\
        0 & 0
    \end{pmatrix},
    \quad f = \begin{pmatrix}
        0 & 0 \\
        1 & 0
    \end{pmatrix}
\end{equation}
subject to the relations
\begin{equation}
    [h, e] = 2 e, \quad [h,f] = -2 f, \quad [e,f] = h.
\end{equation}
For $X, Y \in \mathfrak{sl}_2$, the Killing form $B(X,Y) = 4 \Tr(XY)$, then $B(h,h) = 8, B(e,f)= 4, B(h,e)= B(h,f)=0$. Therefore, the Casimir element $\mathbf{1}$ is written as
\begin{equation}
    \mathbf{1} = \frac{1}{8} h \otimes h + \frac{1}{4} e\otimes f + \frac{1}{4} f \otimes e.
    \label{eq:casimir-4.2.4}
\end{equation}
    \begin{enumerate}[(1)]
    \item Consider the left Y-shaped uni-trivalent tree in Fig. \ref{fig:W_T}, denoted by $Y$. We assign $\Tr$ to the unique trivalent vertex of $Y$. Then the associated weight system $W_Y$ is $W_Y[v_1, v_2, v_3] = \Tr[v_1 \otimes v_2 \otimes v_3]
            = B([v_1, v_2], v_3)$, where each vector $v_j\in \mathfrak{sl}_2$ is attached to the external vectex labelled by $j$. If we write $v_j = a_{jh} h + a_{je} e + a_{jf} f \in \mathfrak{sl}_2$ $(j \in \{1,2,3\})$, then 
    \begin{equation}
        \begin{split}
            &W_Y[v_1, v_2, v_3]\\
            = &8\big ((a_{1e} a_{2f} - a_{1f} a_{2e}) a_{3h} + (a_{1f} a_{2h} - a_{1h} a_{2f})a_{3e} + (a_{1h} a_{2e} - a_{1e} a_{2h}) a_{3f}\big).
        \end{split}
    \end{equation}
    \item Consider the right $H$-shaped uni-trivalent tree in Fig. \ref{fig:W_T}. We assign $\Tr$ to two trivalent vertices and $\mathbf{1}$ in \eqref{eq:casimir-4.2.4} to the unique internal edge of $H$. Then the associated weight system $W_H$ is given by
    \begin{equation}
        \begin{split}
             W_H[v_1, v_2, v_3, v_3] 
             &= 8(a_{1e} a_{2f} - a_{1f} a_{2e})(a_{3e} a_{4f} - a_{3f} a_{4e}) \\
             &\quad + 16 (a_{1f} a_{2h} - a_{1h} a_{2f}) (a_{3h} a_{4e} - a_{3e} a_{4h})\\
             &\quad + 16(a_{1h} a_{2e} - a_{1e} a_{2h})  (a_{3f} a_{4h} - a_{3h} a_{4f}),
        \end{split}
    \end{equation}
    for  $v_j = a_{jh} h + a_{je} e + a_{jf} f \in \mathfrak{sl}_2$ $(j \in \{1,2,3, 4\})$.
    \end{enumerate}
\end{example}

Then, by definition, we have the identity we call {\it anti-symmetric identity} for the uni-trivalent tree $Y$ if we reverse the cyclic order of half-edges for the only trivalent vertex, as in Fig. \ref{fig:AS},
\begin{equation}\label{eq:AS-indentity}
    W_{Y} + W_{\overline{Y}} = 0.
\end{equation}
More generally, for uni-trivalent trees $T_1, T_2$ which are the same except for some small regions where $Y$, $\overline{Y}$ diagrams are inserted respectively. Then, we have
\begin{equation}
    W_{T_1} + W_{T_2} = 0.
\end{equation}
\begin{figure}[h]
\captionsetup{margin=2cm}
 \centering
 \begin{tikzpicture}
\begin{scope}
\node (1) at (-1/4-0.25,1.73/4+0.25) {1};
\node (2) at (1/4+0.25,1.73/4+0.25) {2};
\node (3) at (0,-0.75) {3};
\draw[densely dotted] (0,0) -- (1/4,1.73/4);
\draw[densely dotted] (0.0,0) -- (-1/4,1.73/4);
\draw[densely dotted] (0,0) -- (0, -0.5);
\node at (0,0) {$\cdot$};
\draw[-Stealth, color=blue]  (0.0+0.15,0.0) arc(0:-320:0.15);
\node at (1.5, 0) {$+$};
\node at (4.5, 0) {$=$};
\end{scope}
\begin{scope}[xshift=3cm]
\node (1) at (-1/4-0.25,1.73/4+0.25) {1};
\node (2) at (1/4+0.25,1.73/4+0.25) {2};
\node (3) at (0,-0.75) {3};
\node at (0,0) {$\cdot$};
\draw[densely dotted] (0,0) -- (1/4,1.73/4);
\draw[densely dotted] (0.0,0) -- (-1/4,1.73/4);
\draw[densely dotted] (0,0) -- (0, -0.5);
\draw[-Stealth, color=blue]  (0.0-0.15,0.0) arc(180:180+320:0.15);
\end{scope}
\begin{scope}[xshift=6cm]
\node at (0,0) {$0.$};
\end{scope}
\end{tikzpicture}
 \caption[Anti-symmetric identity for weight systems at one vertex]{Anti-symmetric identity for weight systems. The graphs are denoted by $Y$ and $\overline{Y}$ respectively. In this figure, the blue curved arrows indicate the cyclic ordering for each trivalent vertex (i.e., vertex-wise orientation), and the number labels $1,2,3$ exhibit an ordering for external edges or possibly for the half-edges.}
 \label{fig:AS}
\end{figure}

\begin{figure}[h]
\captionsetup{margin=2cm}
 \centering
 \begin{tikzpicture}
\begin{scope}
\node (1) at (-0.25,0.5) {1};
\node (2) at (1.25,0.5) {2};
\node (3) at (1.25,-0.5) {3};
\node (4) at (-0.25,-0.5) {4};
\draw[densely dotted] (1.0,0.5) -- (-0,0.5);
\draw[densely dotted] (-0.0,-0.5) -- (1.0,-0.5);
\draw[densely dotted] (0.5,0.5) -- (0.5,-0.5);
\node at (0.5,0.5) {$\cdot$};
\node at (0.5,-0.5) {$\cdot$};
\draw[-Stealth, color=blue]  (0.15+0.5,-0.5) arc(0:-320:0.15);
\draw[-Stealth, color=blue]  (0.15+0.5,0.5) arc(0:-320:0.15);
\node at (2, 0) {$+$};
\node at (5, 0) {$+$};
\node at (8, 0) {$=$};
\end{scope}
\begin{scope}[xshift=3cm]
\node (1) at (-0.25,0.5) {1};
\node (2) at (1.25,0.5) {2};
\node (3) at (1.25,-0.5) {3};
\node (4) at (-0.25,-0.5) {4};
\draw[densely dotted] (-0.0,-0.5) -- (-0,0.5);
\draw[densely dotted] (1.0,0.5) -- (1.0,-0.5);
\draw[densely dotted] (-0.0,0) -- (1.0,0);
\node at (0,0) {$\cdot$};
\node at (1,0) {$\cdot$};
\draw[-Stealth, color=blue]  (0.15,0) arc(0:320:0.15);
\draw[-Stealth, color=blue]  (0.15+0.7,0) arc(180:180-320:0.15);
\end{scope}
\begin{scope}[xshift=6cm]
\node (1) at (-0.25,0.5) {1};
\node (2) at (1.25,0.5) {2};
\node (3) at (1.25,-0.5) {3};
\node (4) at (-0.25,-0.5) {4};
\draw[densely dotted] (1.0,-0.5) --(-0,0.5);
\draw[densely dotted] (1.0,0.5) -- (-0.0,-0.5);
\draw[densely dotted] (0.2,-0.3) -- (0.8,-0.3);
\node at (0.2,-0.3) {$\cdot$};
\node at (0.8,-0.3) {$\cdot$};
\draw[-Stealth, color=blue]  (0.35,-0.3) arc(0:-320:0.15);
\draw[-Stealth, color=blue]  (0.95,-0.3) arc(0:-320:0.15);
\end{scope}
\begin{scope}[xshift=9cm]
\node at (0,0) {$0\,.$};
\end{scope}
\end{tikzpicture}
 \caption[Jacobi identity for weight systems -- IHX relation]{Jacobi identity for weight systems. The graphs are denoted by $I, H, X$ respectively.}
 \label{fig:Jacobi}
\end{figure}

\begin{lemma}{(Jacobi identity)}
\begin{enumerate}[(1)]
\item 
For uni-trivalent trees $I,H,X$ given as in Fig. \ref{fig:Jacobi}, we have the following identity on the associated weight systems 
\begin{equation} \label{eq:Jacobi-indentity}
   W_{I} + W_H + W_X = 0.
\end{equation}
\item More generally, for uni-trivalent trees $T_1, T_2, T_3$ which are the same except for some region where $I$, $H$ and $X$ diagrams are inserted respectively. Then, their associated weight systems satisfy the identity
\begin{equation}
     W_{T_1} + W_{T_2} + W_{T_3} = 0.
\end{equation}
\end{enumerate}
\end{lemma}
\begin{proof}
    (1) The proof is given by direct computation as follows. For $v_1, v_2, v_3, v_4 \in \mathfrak{g}$, we have
\begin{equation}
    W_{I}[v_1,v_2, v_3, v_4] = B([v_1,v_2], [v_3, v_4]) = B(v_1, [v_2,[v_3,v_4]]).
\end{equation}
Similarly, 
\begin{equation}
\begin{split}
     & W_{H}[v_1,v_2, v_3, v_4] = B([v_1,v_4], [v_2, v_3]) = B(v_1, [v_4,[v_2,v_3]]),\\
     & W_{X}[v_1,v_2, v_3, v_4] = B([v_1,v_3], [v_4, v_2]) = B(v_1, [v_3,[v_4,v_2]]).
\end{split}
\end{equation}
Therefore, 
\begin{equation}
\begin{split}
     &(W_{I} +  W_{H} +  W_{X})[v_1,v_2, v_3, v_4]\\
     =& B(v_1, [v_2,[v_3,v_4]]) + B(v_1, [v_4,[v_2,v_3]]) + B(v_1, [v_3,[v_4,v_2]])\\
     =& B(v_1, [v_2,[v_3,v_4]] + [v_4,[v_2,v_3]] + [v_3,[v_4,v_2]])\\
     =& 0.
\end{split}
\end{equation}
(2) follows immediately from (1), since the computation of a weight system can be decomposed into that of several pieces of weight systems by construction.
\end{proof}


\subsection{Decorated graphs}\label{ss6.3decorated}
In this subsection, we define a decorated graph which is a variant of one defined in \cite{BC2}. 

\begin{definition}{(Decorated graph)}\label{def:6.3.1August}
 Let $\Gamma$ be a connected graph whose vertices have valency $\geq 3$ (hence no external vertices). A {\it decorated graph} is a graph $\Gamma$ endowed with the following data:
\begin{itemize}
    \item enumerations on the set of edges $e(\Gamma)$ and the set of vertices $v(\Gamma)$, that is, $\Gamma$ is endowed with fixed bijections 
    $$e(\Gamma) \simeq \{1,2,\ldots, |e(\Gamma)|\} ,\; v(\Gamma) \simeq \{1,2,\ldots, |v(\Gamma)|\};$$
    \item directions on edges;
    \item induced order on the set of $h_{\Gamma}(i)$ of half-edges at each vertex $i \in v(\Gamma)$ from the order of $e(\Gamma)$. Here, for a self-loop $e$ connecting the vertex $i$, the order of two half-edges $(i, e, +1)$ and $(i, e, -1)$ is defined so that $(i, e, +1)$ is putted just before $(i, e, -1)$; Note that this order on $h_{\Gamma}(i)$ defines the vertex-wise orientation at vertex $i$;
    \item for each vertex $i$ with valency $n\geq 4$, information of an insertion of oriented uni-trivalent tree $T_i$ with exactly $(2n-3)$ edges and $n$ of them are ordered external vertices corresponding to the $n$ incident half-edges at this vertex $i$. More precisely, we consider a small ball centered at $v$ which intersects on the boundary with half-edges at $n$ distinct points. Then, we endow it with the information of embedding of $T_i$ into the ball so that the $n$-external vertices are put on the intersection points disjointly on the boundary. Here, we also require the embedding of $T_i$ be done so that 
    \begin{enumerate}[(i)]
        \item the order of the external vertices of $T_i$ are given by the order of their corresponding half-edges of $\Gamma$ attached to $i$;
        \item the cyclic order of $T_i$ at a trivalent vertex which connects more than one external vertices are compatible with the order of the half-edges of $\Gamma$ attached to $i$;
    \end{enumerate}
    \item for each vertex $i$, we equip it with the weight system defined as in \eqref{eq:6.2.1August} associated to $T_i$, which is a $\pi_1(M)$-equivariant homomorphism, 
    $$W_i:=W_{T_i}:\otimes_{h\in h_{\Gamma}(i) } \mathfrak{g}_h \rightarrow \bR$$
    which, sometimes, is also regarded as the map
    \begin{equation}
    W_{i} :\mathbb{R} \rightarrow \otimes_{h\in h_{\Gamma}(i) }\mathfrak{g}^\ast_{h}, \; 1\mapsto W_{T_i}
    \label{eq:6.2.1}
    \end{equation}
    where $\mathfrak{g}_h$ (resp. $\mathfrak{g}^\ast_{h}$) is a copy of $\mathfrak{g}$ (resp. $\mathfrak{g}^\ast$). To unify the notation, when $i$ is a trivalent vertex, then we set $T_i$ to be the $Y$-shape uni-trivalent tree and $W_{i}:=\mathrm{Tr}_i\,$.
\end{itemize}
\end{definition}

\begin{remark}
Note that for a vertex $i$ of valency $n\geq 4$, the inserted uni-trivalent tree $T_i$ is required to have exactly $(2n-3)$ edges with $n$ of them being external, and this condition forces the choices of such tree to lie in a finite list of uni-trivalent trees. In this way, if we fix the number of edges and vertices, we only have finitely many different decorated graphs that meet our definition.

\end{remark}

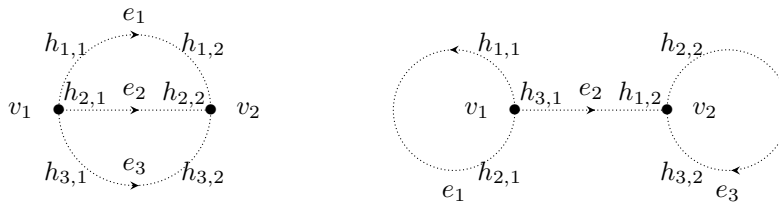
\begin{figure}[h]
\captionsetup{margin=2cm}
 \centering
 \begin{tikzpicture}[font=\footnotesize]
 \begin{scope}
    \draw[densely dotted] (0,0) circle (1.0cm);
    \draw[densely dotted] (-1,0) -- (1,0);
    \draw[->,>=stealth] (0.05,1)--(0.051,1);
    \draw[->,>=stealth] (0.05,0)--(0.051,0);
    \draw[->,>=stealth] (0.05,-1)--(0.051,-1);
    \node at (-1,0) {$\bullet$};
    \node at (1,0) {$\bullet$};
    \node at (-1.5, 0) {$v_1$};
    \node at (1.5,0) {$v_2$};
    \node at (0,1.25) {$e_1$};
    \node at (0,0.25) {$e_2$};
    \node at (0,-0.75) {$e_3$};
    \node at (-0.90,0.85) {${h_{1,1}}$};
    \node at (-0.65,0.20) {${h_{2,1}}$};
    \node at (-0.90,-0.85) {${h_{3,1}}$};
     \node at (0.90,0.85) {${h_{1,2}}$};
    \node at (0.65,0.20) {${h_{2,2}}$};
    \node at (0.90,-0.85) {${h_{3,2}}$};
    \end{scope}
    \begin{scope}[xshift=6cm]
        \draw[densely dotted] (-1.8,0) circle (0.8cm);
        \draw[densely dotted] (-1, 0) -- (1,0);
        \draw[densely dotted] (1.8, 0) circle (0.8cm);
        \node at (0,0.25) {$e_2$};
        \node at (-0.65,0.20) {$h_{3,1}$};
         \node at (0.65,0.20) {$h_{1,2}$};
         \draw[->,>=stealth] (0.05,0)--(0.051,0);
         \node at (-1.2,0.85) {$h_{1,1}$};
         \node at (-1.2,-0.85) {$h_{2,1}$};
         \node at (-1.8,-1.1) {$e_1$};
         \node at (1.2,0.85) {$h_{2,2}$};
         \node at (1.2,-0.85) {$h_{3,2}$};
         \node at (1.8,-1.1) {$e_3$};
         \node at (-1,0) {$\bullet$};
    \node at (1,0) {$\bullet$};
    \node at (-1.5, 0) {$v_1$};
    \node at (1.5,0) {$v_2$};
    \draw[->,>=stealth] (-1.85,0.8)--(-1.851,0.8);
    \draw[->,>=stealth] (1.850,-0.8)--(1.849,-0.8);
    \end{scope}
\end{tikzpicture}

 \caption[Decorated theta graph and dumbbell graph]{Examples of decorated trivalent graphs, whose underlying topological graphs are called theta graph and dumbbell graph respectively. Here, $h_{i,a}$ denotes the $i$-th half-edge at the vertex $v_a$.}
 \label{fig:6.1.theta}
\end{figure}

In this paper, we depict decorated graphs with dashed curves as in Fig. \ref{fig:6.1.theta}. As long as we have the ordering on the half-edges, the decoration $\mathrm{Tr}$ (or the corresponding $\pi_1(M)$-equivariant homomorphism) for each vertex is determined uniquely by the above conventions. Note that in the sequel, we sometimes omit the numberings of vertices, edges, and half-edges or the equivariant homomorphisms for simplicity when depicting decorated graphs (see Lemma \ref{lm:6.6.2}).

\begin{example}\label{ex:trdeco}
For a vertex with valency 4, there are only three ways (up to a sign) to insert uni-trivalent trees as in Fig. \ref{fig:insertion}. Each tree in the small balls carries a vertex-wise orientation uniquely induced from the order of these four half-edges and corresponds to one of the cases in Fig. \ref{fig:Jacobi} up to a sign (see also \eqref{eq:AS-self-loop} and \eqref{eq:self-loop_jacobi}).
\begin{figure}[h]
\captionsetup{margin=2cm}
\centering
 \adjustbox{scale=1.4,center}{
\begin{tikzpicture}
    \draw[fill=black!40, black!40] (4.5, 0) circle (0.4cm);
   \node at (4.5, 0) {$\bullet$};
 	\draw[densely dotted] (4, -0.5) -- (5, 0.5);
 	\draw[densely dotted] (4, 0.5) -- (5, -0.5);
  \draw[densely dashed, ->] (6,0) -- (5.5,0);
  \begin{scope}[xshift=7cm, rotate=90] 
    \draw (0,0) circle (0.4cm);
    \draw[densely dotted] (-0.2, 0) -- (0.2, 0);
     \node at (0.2,0) {$\cdot$};
      \node at (-0.2,0) {$\cdot$};
    \draw[densely dotted] (-0.2, 0) -- (-0.3,0.25);
    \draw[densely dotted] (-0.2, 0) -- (-0.3,-0.25);
    \draw[densely dotted] (0.2, 0) -- (0.3,0.25);
    \draw[densely dotted] (0.2, 0) -- (0.3,-0.25);
    \end{scope}
     \begin{scope}[xshift=8.5cm] 
    \draw (0,0) circle (0.4cm);
     \node at (0.2,0) {$\cdot$};
      \node at (-0.2,0) {$\cdot$};
    \draw[densely dotted] (-0.2, 0) -- (0.2, 0);
    \draw[densely dotted] (-0.2, 0) -- (-0.3,0.25);
    \draw[densely dotted] (-0.2, 0) -- (-0.3,-0.25);
    \draw[densely dotted] (0.2, 0) -- (0.3,0.25);
    \draw[densely dotted] (0.2, 0) -- (0.3,-0.25);
    \end{scope}
     \begin{scope}[xshift=10cm] 
    \draw (0,0) circle (0.4cm);
    \draw[densely dotted] (-0.2, -0.17) -- (0.2, -0.17);
    \draw[densely dotted] (0.3,-0.25) -- (-0.3,0.25);
    \draw[densely dotted] (0.3,0.25) -- (-0.3,-0.25);
     \node at (0.2,-0.17) {$\cdot$};
      \node at (-0.2,-0.17) {$\cdot$};
    \end{scope}
    \end{tikzpicture}
    }
     \caption[Weight systems for a vertex of valency $4$]{Three ways to insert a uni-trivalent tree in place of a 4-valent vertex. These graphs are denoted by $I, H$, and $X$ respectively.}
 \label{fig:insertion}
\end{figure}
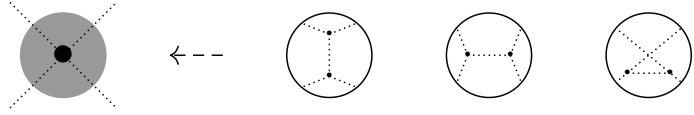

Generally, one can decorate a vertex with valency $4$ with a more complicated weight system, for example, as in Fig. \ref{fig:insertionMore}, but we exclude this case from our definition.
\begin{figure}[h]
\captionsetup{margin=2cm}
\centering
 \adjustbox{scale=1.4,center}{
\begin{tikzpicture}
    \draw[fill=black!40, black!40] (4.5, 0) circle (0.4cm);
   \node at (4.5, 0) {$\bullet$};
 	\draw[densely dotted] (4, -0.5) -- (5, 0.5);
 	\draw[densely dotted] (4, 0.5) -- (5, -0.5);
  \draw[densely dashed, ->] (6,0) -- (5.5,0);
  \begin{scope}[xshift=7.3cm] 
    \draw (0,0) circle (0.8cm);
    \draw[densely dotted]  (0.4, 0.68) -- (0, 0.4);
    \draw[densely dotted]  (-0.4, 0.68) -- (0,0.4);
    \draw[densely dotted]  (-0.4, -0.68) -- (0,-0.4);
    \draw[densely dotted]  (0.4, -0.68) -- (0,-0.4);
      \draw[densely dotted]  (0, 0.2) -- (0,0.4);
    \draw[densely dotted]  (0, -0.2) -- (0,-0.4);
          \draw[densely dotted]  (-0.4, 0) -- (0.4,0);
       \draw[densely dotted]  (0, 0.2) -- (0.4,0);
              \draw[densely dotted]  (-0.4, 0) -- (0,0.2);
       \draw[densely dotted]  (0, -0.2) -- (0.4,0);
              \draw[densely dotted]  (-0.4, 0) -- (0,-0.2);
              \node at (0,0.4) {$\cdot$};
                   \node at (0,-0.4) {$\cdot$};
      \node at (-0.4,0) {$\cdot$};
           \node at (0.4,0) {$\cdot$};
      \node at (0,-0.2) {$\cdot$};
           \node at (0,0.2) {$\cdot$};
    \end{scope}
    \end{tikzpicture}
    }
     \caption[A complicated weight system that is not allowed]{A more complicated weight system, which is not allowed.}
 \label{fig:insertionMore}
\end{figure}
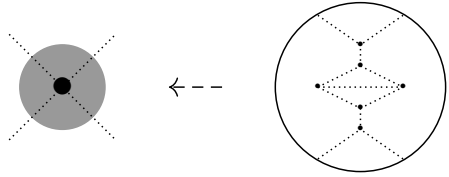
\end{example}

\subsection{Graph complex for an acyclic local system}\label{ss6.4complex}
Now let $G$ be a connected semi-simple Lie group with Lie algebra $\mathfrak{g}$ and the Killing form $B$ as considered in Section \ref{section3}. In this subsection, we introduce a graph complex for the acyclic local systems associated to $\mathfrak{g}$ and $B$, more precisely, the trace form $\mathrm{Tr}$. 

Let $\widetilde{\mathcal{GC}}_{\mathfrak{g}}$ be the vector space spanned by all the decorated graphs as in Definition \ref{def:6.3.1August} over $\mathbb{Q}$. These spaces of decorated graphs are bigraded by the following order and degree:
\begin{equation}
\begin{split}
& \ord(\Gamma) = |e(\Gamma)| - |v(\Gamma)|,  \\
& \deg(\Gamma) = 2  |e(\Gamma)| - 3 |v(\Gamma)| .
\end{split}
\label{eq:orddegJuly} 
\end{equation}
Note that the decorated graphs that we consider here always have $\deg \geq 0$.

In some context, we also like to talk about the loops for a connected graph $\Gamma$. Viewing the (topological) graph as a $CW$-complex, then the Euler characteristic number is
\begin{equation}
\chi(\Gamma)=|v(\Gamma)|-|e(\Gamma)|=-\ord(\Gamma)=1-\ell,
\end{equation}
where $\ell$ corresponds to the first Betti number of $\Gamma$ hence the number of loops in $\Gamma$. Note that the number of loops as above only makes sense for a connected graph, if the graph is not connected, we also need to consider the number of connected components to conclude the number of loops in a topological sense.

\begin{remark}\begin{itemize}
 \item Any trivalent graph $\Gamma$, that is, a graph whose all the vertices have valency $3$, is of $\deg(\Gamma) = 0$. Since we assume that any graph in $\widetilde{\mathcal{GC}}_{\mathfrak{g}}$ have internal vertices with valencies at least $3$, a finite graph without any external edges and of degree $0$ has to be a trivalent graph. A (nonempty) trivalent graph has at least $2$ vertices and $3$ edges, so the least order is $1$. The trivalent graphs of order $1$ have only two possibilities: the theta graph and the dumbbell graph, both are connected.
    \item For a trivalent graph $\Gamma$, its order $\ord(\Gamma)=\frac{1}{2}|v(\Gamma)|$ defined in \eqref{eq:orddegJuly} agrees with its {\it degree} customarily used in the theory of finite-type (Vassiliev) invariants. 
\end{itemize}
\end{remark}

We define an equivalence relation on $\widetilde{\mathcal{GC}}_{\mathfrak{g}}$ as follows: if two decorated graphs $\Gamma$ and $\Gamma'$ differ by only  
\begin{enumerate}[(1)]
    \item permutation of numberings for all edges which induces $k$ times change in the total of cyclic orders of associated trees at vertices, where the changes on the cyclic orders are forced by the compatibility of cyclic orders on the associated trees with the new numberings on their external edge, 
    \item edge (including self-loop edges and non-self-loop edges) direction reversals of times $m$; let $(-1)^{m'}$ denote the total sign change of the cyclic orders of the associated trees induced by the direction changes on the self-loop edges (since the cyclic orders are not affected by direction reversals on non-self-loop edges),
   \item  permutation of numberings of vertices, let $(-1)^d$ denote the sign, 
\end{enumerate}
then we set  $\sign(\Gamma, \Gamma') = (-1)^{k+m+m'+d}$, and
\begin{equation}
	\Gamma = \sign(\Gamma, \Gamma') \cdot \Gamma'.
 \label{eq:6.2.4}
\end{equation}

We also introduce a relation, called  {\it internal vertex-wise AS relation}, of connected decorated graphs
as follows.  
Let $\Gamma_Y$, $\Gamma_{\overline{Y}}$ be two connected decorated graphs with the same underlying topological graph and the same decoration except that, at one fixed internal vertex,  embedded uni-trivalent trees are different but related by anti-symmetric identity as \eqref{eq:AS-indentity}. Then, we set
\begin{equation} \label{eq:IV-AS}
    \Gamma_Y + \Gamma_{\overline{Y}} =0.
\end{equation}

Similarly, we define  {\it internal vertex-wise IHX relation}. Let $\Gamma_I$, $\Gamma_H$, $\Gamma_X$ be three connected decorated graphs with the same underlying topological graph and the same decoration except that, at one fixed internal vertex,  embedded uni-trivalent trees are different but related by Jacobi identity as \eqref{eq:Jacobi-indentity}. Then, we set
\begin{equation}\label{eq:IV-IHX}
    \Gamma_I + \Gamma_H + \Gamma_X =0.
\end{equation}

For a self-loop edge, we now give more details to clarify the equivalence relation under the change of direction. Let $\Gamma$ be a decorated graph, and let $v$ be a vertex in $\Gamma$ of valency $3$ and attached by a self-loop edge, let $\Gamma'$ be the decorated graph obtained by reverse the direction of this self-loop edge attached to $v$, then our equivalence relation shows 
\begin{equation}
\Gamma=\Gamma'.
\end{equation}

However, this situation might be different for the vertex with higher valency.  For a vertex with valency $4$, this relation is presented by the following figure:
\begin{equation} \label{eq:AS-self-loop}
 \adjustbox{scale=1.2,center}{
    \begin{tikzpicture}[font=\scriptsize]
    \draw[fill=black!40, black!40] (4.5, 0) circle (0.4cm);
    \node at (3,0) {$\Big ($};
    \node at (8.6,0) {$\Big )$};
   \node at (4.5, 0) {$\bullet$};
   \node at (4.7, 0) {$v$};
 	\draw[densely dotted] (4, -0.5) -- (4.5, 0);
   \node at (3.8, -0.4) {${h_{1,v}}$};
 	\draw[densely dotted] (5, -0.5) -- (4.5, 0);
   \node at (5.2, -0.4) {${h_{2,v}}$};
  \draw[densely dotted] (4.5, 0.4) circle (0.4cm);
  \node at (3.8, 0.4) {${h_{3,v}}$};
  \node at (5.2, 0.4) {${h_{4,v}}$};
  \draw[->,>=stealth] (4.49,0.8)--(4.51,0.8);
  \draw[densely dashed, ->] (6,0) -- (5.5,0);
  \node at (9, 0) {$=$};
  \node at (9.5, 0) {$-$};
    \node at (7.9,0) {{ $:= -X$}};
     \begin{scope}[xshift=7cm] 
    \draw (0,0) circle (0.4cm);
         \node at (0.2,0) {$\cdot$};
      \node at (-0.2,0) {$\cdot$};
    \draw[densely dotted] (-0.2, 0) -- (0.2, 0);
    \draw[densely dotted] (-0.2, 0) -- (-0.3,0.25);
        \node at (-0.45,0.35) {$3$};
    \draw[densely dotted] (-0.2, 0) -- (-0.3,-0.25);
    \node at (-0.45,-0.35) {$1$};
    \draw[densely dotted] (0.2, 0) -- (0.3,0.25);
     \node at (0.45,0.35) {$4$};
    \draw[densely dotted] (0.2, 0) -- (0.3,-0.25);
     \node at (0.45,-0.35) {$2$};
    \node at (0, 0.6) {$T_v$};
    \end{scope}
    \begin{scope}[xshift=7cm]
    \draw[fill=black!40, black!40] (4.5, 0) circle (0.4cm);
    \node at (3,0) {$\Big ($};
    \node at (8.6,0) {$\Big )$};
   \node at (4.5, 0) {$\bullet$};
   \node at (4.7, 0) {$v$};
 	\draw[densely dotted] (4, -0.5) -- (4.5, 0);
   \node at (3.8, -0.4) {${h_{1,v}}$};
 	\draw[densely dotted] (5, -0.5) -- (4.5, 0);
   \node at (5.2, -0.4) {${h_{2,v}}$};
  \draw[densely dotted] (4.5, 0.4) circle (0.4cm);
  \node at (3.8, 0.4) {${h_{4,v}}$};
  \node at (5.2, 0.4) {${h_{3,v}}$};
  \draw[<-,>=stealth] (4.49,0.8)--(4.51,0.8);
  \draw[densely dashed, ->] (6,0) -- (5.5,0);
    \node at (7.9,0) {{$:= H$}};
     \begin{scope}[xshift=7cm] 
    \draw (0,0) circle (0.4cm);
         \node at (0.2,0) {$\cdot$};
      \node at (-0.2,0) {$\cdot$};
    \draw[densely dotted] (-0.2, 0) -- (0.2, 0);
    \draw[densely dotted] (-0.2, 0) -- (-0.3,0.25);
        \node at (-0.45,0.35) {$4$};
    \draw[densely dotted] (-0.2, 0) -- (-0.3,-0.25);
    \node at (-0.45,-0.35) {$1$};
    \draw[densely dotted] (0.2, 0) -- (0.3,0.25);
     \node at (0.45,0.35) {$3$};
    \draw[densely dotted] (0.2, 0) -- (0.3,-0.25);
     \node at (0.45,-0.35) {$2$};
    \node at (0, 0.6) {$T'_v$};
    \end{scope}
    \end{scope}
\end{tikzpicture}
}
\end{equation}
where black-filled circles denote the small balls centered at an internal vertex which carries the data of embedded uni-trivalent trees, and the figures for $T_v$ or $T'_v$ mean that the embedded uni-trivalent trees are the respective $-X$ and $H$ defined in Fig. \ref{fig:Jacobi}. 
Note that by our conditions for $T_v$ or $T'_v$ given in Definition \ref{def:6.3.1August}, the vertex-wise orientations of $T_v$ and $T'_v$ are uniquely determined by the ordering of the external vertices. 
Thus, we do not emphasize the cyclic order for each vertex in the figures of \eqref{eq:AS-self-loop}. 
This is also the reason why we need to put the minus sign in front of $X$. 
Indeed, noting that the cyclic orders in $T_v$ remain the same after the exchange of labels $3$ and $4$,  we have to put a minus sign on the right-hand side of \eqref{eq:AS-self-loop}.

In the two sides of \eqref{eq:AS-self-loop}, the different directions of the self-loop yield a minus sign. As a consequence, we conclude an equivalence between the weight systems $X$ and $H$ at this vertex with valency $4$ and attached by a self-loop. Then combining it with the internal vertex-wise IHX relation \eqref{eq:IV-IHX}, we can conclude the following nontrivial identity.

\begin{equation}\label{eq:self-loop_jacobi}
 \adjustbox{scale=1.2,center}{
     \begin{tikzpicture}[font=\scriptsize]
    \draw[fill=black!40, black!40] (4.5, 0) circle (0.4cm);
    \node at (3.1,0) {$-2\Big ($};
    \node at (8.4,0) {$\Big )$};
   \node at (4.5, 0) {$\bullet$};
   \node at (4.7, 0) {$v$};
 	\draw[densely dotted] (4, -0.5) -- (4.5, 0);
   \node at (3.8, -0.4) {${h_{1,v}}$};
 	\draw[densely dotted] (5, -0.5) -- (4.5, 0);
   \node at (5.2, -0.4) {${h_{2,v}}$};
  \draw[densely dotted] (4.5, 0.4) circle (0.4cm);
  \node at (3.8, 0.4) {${h_{3,v}}$};
  \node at (5.2, 0.4) {${h_{4,v}}$};
  \draw[->,>=stealth] (4.49,0.8)--(4.51,0.8);
  \draw[densely dashed, ->] (6,0) -- (5.5,0);
  \node at (8.9, 0) {$+$};
    \node at (7.9,0) {{$:= -X$}};
     \begin{scope}[xshift=7cm] 
    \draw (0,0) circle (0.4cm);
         \node at (0.2,0) {$\cdot$};
      \node at (-0.2,0) {$\cdot$};
    \draw[densely dotted] (-0.2, 0) -- (0.2, 0);
    \draw[densely dotted] (-0.2, 0) -- (-0.3,0.25);
        \node at (-0.45,0.35) {$3$};
    \draw[densely dotted] (-0.2, 0) -- (-0.3,-0.25);
    \node at (-0.45,-0.35) {$1$};
    \draw[densely dotted] (0.2, 0) -- (0.3,0.25);
     \node at (0.45,0.35) {$4$};
    \draw[densely dotted] (0.2, 0) -- (0.3,-0.25);
     \node at (0.45,-0.35) {$2$};
    \end{scope}
    \begin{scope}[xshift=6.2cm]
    \draw[fill=black!40, black!40] (4.5, 0) circle (0.4cm);
    \node at (3.2,0) {$\Big ($};
    \node at (8.6,0) {$\Big )\;=0$};
   \node at (4.5, 0) {$\bullet$};
   \node at (4.7, 0) {$v$};
 	\draw[densely dotted] (4, -0.5) -- (4.5, 0);
   \node at (3.8, -0.4) {${h_{1,v}}$};
 	\draw[densely dotted] (5, -0.5) -- (4.5, 0);
   \node at (5.2, -0.4) {${h_{2,v}}$};
  \draw[densely dotted] (4.5, 0.4) circle (0.4cm);
  \node at (3.8, 0.4) {${h_{3,v}}$};
  \node at (5.2, 0.4) {${h_{4,v}}$};
  \draw[->,>=stealth] (4.49,0.8)--(4.51,0.8);
  \draw[densely dashed, ->] (6,0) -- (5.5,0);
    \node at (7.9,0) {{ $:= I$}};
     \begin{scope}[xshift=7cm, rotate=90] 
    \draw (0,0) circle (0.4cm);
         \node at (0.2,0) {$\cdot$};
      \node at (-0.2,0) {$\cdot$};
    \draw[densely dotted] (-0.2, 0) -- (0.2, 0);
    \draw[densely dotted] (-0.2, 0) -- (-0.3,0.25);
        \node at (-0.45,0.35) {$1$};
    \draw[densely dotted] (-0.2, 0) -- (-0.3,-0.25);
    \node at (-0.45,-0.35) {$2$};
    \draw[densely dotted] (0.2, 0) -- (0.3,0.25);
     \node at (0.45,0.35) {$3$};
    \draw[densely dotted] (0.2, 0) -- (0.3,-0.25);
     \node at (0.45,-0.35) {$4$};
    \end{scope}
    \end{scope}
\end{tikzpicture}
}
\end{equation}

\begin{remark}
When considering trivial local systems, as is well known, graphs with self-loop edges are zero by AS relation (that is, antisymmetry of internal vertices given in \cite[Theorem 6 (1)]{BN}) by using arguments in \cite[Section 2.4]{BN}. However, this is not the case for non-trivial local systems. This can be explained as follows. For a trivial local system, (if we only consider the Lie algebra factor of associated integrations) every internal edge corresponds to the element $\bf{1}$ which lies in the symmetric part of $\mathfrak{g} \otimes \mathfrak{g}$, whereas for a non-trivial local system concerned here, we associate self-loop edges anti-symmetric element of $E_{\rho} \otimes E_{\rho}$ so that the AS relation does not imply vanishings of graphs with self-loops. See also \cite[Page 180]{AS2}.
\end{remark}

\begin{figure}[h]
\captionsetup{margin=2cm}
 \centering
 \begin{tikzpicture}[scale=1.2,font=\footnotesize]
  \begin{scope}[xshift=-5cm] 
    \draw (0,0) circle (0.4cm);
    \node at (-0.8,0) {{$H= $}};
         \node at (0.2,0) {$\cdot$};
      \node at (-0.2,0) {$\cdot$};
    \draw[densely dotted] (-0.2, 0) -- (0.2, 0);
    \draw[densely dotted] (-0.2, 0) -- (-0.3,0.25);
        \node at (-0.45,0.35) {$4$};
    \draw[densely dotted] (-0.2, 0) -- (-0.3,-0.25);
    \node at (-0.45,-0.35) {$1$};
    \draw[densely dotted] (0.2, 0) -- (0.3,0.25);
     \node at (0.45,0.35) {$3$};
    \draw[densely dotted] (0.2, 0) -- (0.3,-0.25);
     \node at (0.45,-0.35) {$2$};
    \node at (0, 0.6) {$T_v$};
    \end{scope}
    \begin{scope}
    \draw[densely dashed, ->] (-3.9,0) -- (-3.4,0);
     \draw[fill=black!40, black!40] (-2.6, 0) circle (0.4cm);
        \draw[densely dotted] (-1.8,0) circle (0.8cm);
        \draw[densely dotted] (-1, 0) -- (0,0);
        \draw[densely dotted] (-3.6, 0.8) -- (-2.6,0);
        \draw[densely dotted] (-3.6, -0.8) -- (-2.6,0);
        \node at (-2.9,0.55) {${h_{1,v}}$};
         \node at (-2.9,-0.55) {${h_{2,v}}$};
              \node at (-1.95,0.55) {${h_{3,v}}$};
         \node at (-1.95,-0.55) {${h_{4,v}}$};
         \node at (-1.8,-1.1) {$e'$};
          \node at (-1.8,1.1) {$e$};
         \node at (-1,0) {$\bullet$};
         \node at (-2.6,0) {$\bullet$};
    \node at (-1.2, 0) {$v'$};
     \node at (-2.4, 0) {$v$};
    \draw[->,>=stealth] (-1.80,0.8)--(-1.79,0.8);
    \draw[->,>=stealth] (-1.80,-0.8)--(-1.79,-0.8);
    \end{scope}
\end{tikzpicture}
 \caption[Valency-$4$ vertex with two incident edges attached to the same vertex]{An example of two non-self-loop edges with the same ending vertices, we assume $e < e'$ in the given numberings on the edges of $\Gamma$.}
 \label{fig:Sept7}
\end{figure}
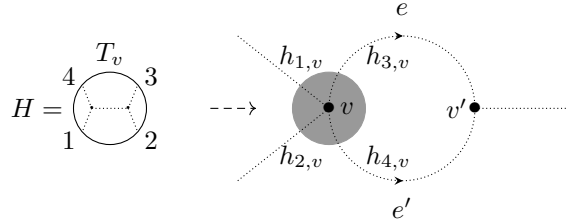

Let us consider another example with a vertex $v$ of valency $4$ as in Fig. \ref{fig:Sept7}, where
the vertex $v$ of valency $4$ has two attached edges connecting to the same trivalent vertex $v'$.

Now we exchange the numberings $e$ and $e'$ for these two edges described in Fig. \ref{fig:Sept7} to a new decorated graph $\Gamma'$. As a consequence, the number labels $3$ and $4$ in $T_v$ are exchanged, but the induced cyclic orders at each vertex in $T_v$ remain the same, so there is no sign produced for $v$. However, since $v'$ is trivalent, exchanging $e$ and $e'$ produces a factor $(-1)$ for the equivalence relation between $\Gamma$ and $\Gamma'$.  Then, combining the internal IHX relation, we conclude an identity in Fig. \ref{fig:Sept8}.


\begin{figure}[h!]
\captionsetup{margin=2cm}
 \centering
 \begin{tikzpicture}[scale=1.2, font=\footnotesize]
  \begin{scope}[xshift=-5cm] 
    \draw (0,0) circle (0.4cm);
    \node at (-1.3,0) {$2\;\Big ($};
    \node at (-0.8,0) {{$H= $}};
         \node at (0.2,0) {$\cdot$};
      \node at (-0.2,0) {$\cdot$};
    \draw[densely dotted] (-0.2, 0) -- (0.2, 0);
    \draw[densely dotted] (-0.2, 0) -- (-0.3,0.25);
        \node at (-0.45,0.35) {$4$};
    \draw[densely dotted] (-0.2, 0) -- (-0.3,-0.25);
    \node at (-0.45,-0.35) {$1$};
    \draw[densely dotted] (0.2, 0) -- (0.3,0.25);
     \node at (0.45,0.35) {$3$};
    \draw[densely dotted] (0.2, 0) -- (0.3,-0.25);
     \node at (0.45,-0.35) {$2$};
    \node at (0, 0.6) {$T_v$};
    \end{scope}
    \begin{scope}
    \draw[densely dashed, ->] (-3.9,0) -- (-3.4,0);
     \draw[fill=black!40, black!40] (-2.6, 0) circle (0.4cm);
        \draw[densely dotted] (-1.8,0) circle (0.8cm);
        \draw[densely dotted] (-1, 0) -- (-0.5,0);
        \node at (-0.15,0) {$\Big )\;+\;$};
        \draw[densely dotted] (-3.6, 0.8) -- (-2.6,0);
        \draw[densely dotted] (-3.6, -0.8) -- (-2.6,0);
        \node at (-2.9,0.55) {${h_{1,v}}$};
         \node at (-2.9,-0.55) {${h_{2,v}}$};
              \node at (-2,0.55) {${h_{3,v}}$};
         \node at (-2,-0.55) {${h_{4,v}}$};
         \node at (-1.8,-1.1) {$e'$};
          \node at (-1.8,1.1) {$e$};
         \node at (-1,0) {$\bullet$};
         \node at (-2.6,0) {$\bullet$};
    \node at (-1.2, 0) {$v'$};
     \node at (-2.4, 0) {$v$};
    \draw[->,>=stealth] (-1.80,0.8)--(-1.79,0.8);
    \draw[->,>=stealth] (-1.80,-0.8)--(-1.79,-0.8);
    \end{scope}

    \begin{scope}[xshift=1.5cm] 
    \draw (0,0) circle (0.4cm);
    \node at (-1.3,0) {$\;\Big ($};
    \node at (-0.8,0) {{ $I= $}};
         \node at (0.2,0) {$\cdot$};
      \node at (-0.2,0) {$\cdot$};
    \draw[densely dotted] (-0.2, 0) -- (0.2, 0);
    \draw[densely dotted] (-0.2, 0) -- (-0.3,0.25);
        \node at (-0.45,0.35) {$1$};
    \draw[densely dotted] (-0.2, 0) -- (-0.3,-0.25);
    \node at (-0.45,-0.35) {$2$};
    \draw[densely dotted] (0.2, 0) -- (0.3,0.25);
     \node at (0.45,0.35) {$3$};
    \draw[densely dotted] (0.2, 0) -- (0.3,-0.25);
     \node at (0.45,-0.35) {$4$};
    \node at (0, 0.6) {$T_v$};
    \end{scope}
    \begin{scope}[xshift=6.5cm]
    \draw[densely dashed, ->] (-3.9,0) -- (-3.4,0);
     \draw[fill=black!40, black!40] (-2.6, 0) circle (0.4cm);
        \draw[densely dotted] (-1.8,0) circle (0.8cm);
        \draw[densely dotted] (-1, 0) -- (-0.5,0);
        \node at (0,0) {$\Big )\;=0.$};
        \draw[densely dotted] (-3.6, 0.8) -- (-2.6,0);
        \draw[densely dotted] (-3.6, -0.8) -- (-2.6,0);
        \node at (-2.9,0.55) {${h_{1,v}}$};
         \node at (-2.9,-0.55) {${h_{2,v}}$};
              \node at (-2,0.55) {${h_{3,v}}$};
         \node at (-2,-0.55) {${h_{4,v}}$};
         \node at (-1.8,-1.1) {$e'$};
          \node at (-1.8,1.1) {$e$};
         \node at (-1,0) {$\bullet$};
         \node at (-2.6,0) {$\bullet$};
    \node at (-1.2, 0) {$v'$};
     \node at (-2.4, 0) {$v$};
    \draw[->,>=stealth] (-1.80,0.8)--(-1.79,0.8);
    \draw[->,>=stealth] (-1.80,-0.8)--(-1.79,-0.8);
    \end{scope}
\end{tikzpicture}
 \caption[A special case of internal IHX relation]{A special case of internal IHX relation.}
 \label{fig:Sept8}
\end{figure}

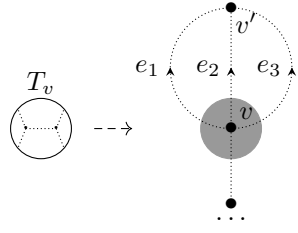
\begin{figure}[h!]
\captionsetup{margin=2cm}
 \centering
 \begin{tikzpicture}[scale=1., font=\footnotesize]
  \begin{scope}[xshift=-2.5cm] 
    \draw (0,0) circle (0.4cm);
       \node at (0.2,0) {$\cdot$};
      \node at (-0.2,0) {$\cdot$};
    \draw[densely dotted] (-0.2, 0) -- (0.2, 0);
    \draw[densely dotted] (-0.2, 0) -- (-0.3,0.25);

    \draw[densely dotted] (-0.2, 0) -- (-0.3,-0.25);

    \draw[densely dotted] (0.2, 0) -- (0.3,0.25);

    \draw[densely dotted] (0.2, 0) -- (0.3,-0.25);

    \node at (0, 0.6) {$T_v$};
\draw[densely dashed, ->] (0.7,0) -- (1.2,0);

    \end{scope}
    \begin{scope}[rotate=90, xshift=2.6cm]
    
     \draw[fill=black!40, black!40] (-2.6, 0) circle (0.4cm);
        \draw[densely dotted] (-1.8,0) circle (0.8cm);
        
\draw[densely dotted] (-2.6, 0) -- (-1,0);
       \node at (-3.6, 0) {$\bullet$};
\node at (-1.8, 0.3) {$e_2$};

        \draw[densely dotted] (-3.6, 0) -- (-2.6,0);
\node at (-3.8, 0) {$\ldots$};

     \node at (-1.8,-0.5) {$e_3$};
          \node at (-1.8,1.1) {$e_1$};
         \node at (-1,0) {$\bullet$};
         \node at (-2.6,0) {$\bullet$};
    \node at (-1.2, -0.2) {$v'$};
     \node at (-2.4, -0.2) {$v$};
    \draw[->,>=stealth] (-1.80,0.8)--(-1.79,0.8);
    \draw[->,>=stealth] (-1.80,-0.8)--(-1.79,-0.8);
 \draw[->,>=stealth] (-1.80,0)--(-1.79,0);
  \end{scope}
\end{tikzpicture}
 \caption[Valency-$4$ vertex with three incident edges attached to the same vertex]{Vertex $v$ with valency $4$, vertex $v'$ with valency $3$, $T_v$ is one of the cases $\{I,H,X\}$, the above decorated graph is identified to be $0$ by the internal IHX relation and sign relation.}
 \label{fig:Sept9}
\end{figure}

Another situation for a vertex $v$ with valency $4$ is given in Fig. \ref{fig:Sept9}, where the decorated graph is always identified to be zero by the internal IHX relations.

\begin{definition}\label{def:CGA_ac_gc}
We define $\mathcal{GC}_{\mathfrak{g}}=\widetilde{\mathcal{GC}}_{\mathfrak{g}}/{\sim}$ as the graded commutative algebra over $\mathbb{Q}$ generated by equivalent classes of decorated connected graphs in $\widetilde{\mathcal{GC}}_{\mathfrak{g}}$, subject to 
\begin{itemize}
\item the sign relation \eqref{eq:6.2.4},
\item internal vertex-wise AS relation \eqref{eq:IV-AS},
    \item internal vertex-wise IHX relation \eqref{eq:IV-IHX},
\end{itemize}
the (graded commutative) algebra structure on $\mathcal{GC}_{\mathfrak{g}}$ is given by disjoint union (denoted by $\cup$), which is defined as follows, the numberings on the edges and internal vertices of $\Gamma \cup \Gamma'$ are given as keeping the same for $\Gamma$ and shifting the numberings for $\Gamma'$ by adding $|e(\Gamma)|$, $|v(\Gamma)|$ respectively. Note that the disjoint union of two decorated connected graphs, viewed as a newly decorated graph, leads to the summation of their respective orders and degrees so that a multiple of a connected decorated graph with order $n$ and degree $t$ is still considered as having the same order and degree, which is a different object from the disjoint union of multiple copies of this graph (see Fig. \ref{fig:6.4.4}). We have the following commutative relation
\begin{equation}
    \Gamma \cup \Gamma' =(-1)^{\mathrm{deg}(\Gamma)\mathrm{deg}(\Gamma')}\Gamma' \cup \Gamma.
\end{equation}
In particular, if $\Gamma$ has an odd degree (equivalently, has odd number of vertices), then we have (in $\mathcal{GC}_{\mathfrak{g}}$)
\begin{equation}
    \Gamma \cup \Gamma = 0.
\end{equation}

\begin{figure}[h]
\captionsetup{margin=2cm}
 \centering
 \begin{tikzpicture}[font=\footnotesize]
 \begin{scope}
    \draw[densely dotted] (0,0) circle (1.0cm);
    \draw[densely dotted] (-1,0) -- (1,0);
    \draw[->,>=stealth] (0.05,1)--(0.051,1);
    \draw[->,>=stealth] (0.05,0)--(0.051,0);
    \draw[->,>=stealth] (0.05,-1)--(0.051,-1);
    \node at (-1,0) {$\bullet$};
    \node at (1,0) {$\bullet$};
    \node at (-1.3, 0) {$v_1$};
    \node at (-2,0) {${\large 2\; \Big(}$};
    \node at (1.3,0) {$v_2$};
    \node at (1.8,0) {${\large\Big)}$};
    \node at (2.4,0) {${\neq}$};
    \node at (0,1.25) {$e_1$};
    \node at (0,0.25) {$e_2$};
    \node at (0,-0.75) {$e_3$};
    \node at (-0.9,0.85) {$h_{1,1}$};
    \node at (-0.65,0.20) {$h_{2,1}$};
    \node at (-0.9,-0.85) {$h_{3,1}$};
     \node at (0.9,0.85) {$h_{1,2}$};
    \node at (0.65,0.20) {$h_{2,2}$};
    \node at (0.9,-0.85) {$h_{3,2}$};
    \end{scope}
    \begin{scope}[xshift=4.5cm]
\draw[densely dotted] (0,0) circle (1.0cm);
    \draw[densely dotted] (-1,0) -- (1,0);
    \draw[->,>=stealth] (0.05,1)--(0.051,1);
    \draw[->,>=stealth] (0.05,0)--(0.051,0);
    \draw[->,>=stealth] (0.05,-1)--(0.051,-1);
    \node at (-1,0) {$\bullet$};
    \node at (1,0) {$\bullet$};
    \node at (-1.3, 0) {$v_1$};
    \node at (1.3,0) {$v_2$};
    \node at (0,1.25) {$e_1$};
    \node at (0,0.25) {$e_2$};
    \node at (0,-0.75) {$e_3$};
    \node at (-0.9,0.85) {$h_{1,1}$};
    \node at (-0.65,0.20) {$h_{2,1}$};
    \node at (-0.9,-0.85) {$h_{3,1}$};
     \node at (0.9,0.85) {$h_{1,2}$};
    \node at (0.65,0.20) {$h_{2,2}$};
    \node at (0.9,-0.85) {$h_{3,2}$};    
    \node at (2,0) {$\large\cup$};
    \end{scope}

     \begin{scope}[xshift=8.5cm]
\draw[densely dotted] (0,0) circle (1.0cm);
    \draw[densely dotted] (-1,0) -- (1,0);
    \draw[->,>=stealth] (0.05,1)--(0.051,1);
    \draw[->,>=stealth] (0.05,0)--(0.051,0);
    \draw[->,>=stealth] (0.05,-1)--(0.051,-1);
    \node at (-1,0) {$\bullet$};
    \node at (1,0) {$\bullet$};
    \node at (-1.3, 0) {$v_3$};
    \node at (1.3,0) {$v_4$};
    \node at (0,1.25) {$e_4$};
    \node at (0,0.25) {$e_5$};
    \node at (0,-0.75) {$e_6$};
    \node at (-0.9,0.85) {$h_{1,3}$};
    \node at (-0.65,0.20) {$h_{2,3}$};
    \node at (-0.9,-0.85) {$h_{3,3}$};
     \node at (0.9,0.85) {$h_{1,4}$};
    \node at (0.65,0.20) {$h_{2,4}$};
    \node at (0.9,-0.85) {$h_{3,4}$};    
    \end{scope}
\end{tikzpicture}
 \caption[Union operator of decorated graphs]{Multiple of a connected graph considered different from the disjoint union of multiple copies of the graph.}
 \label{fig:6.4.4}
\end{figure}
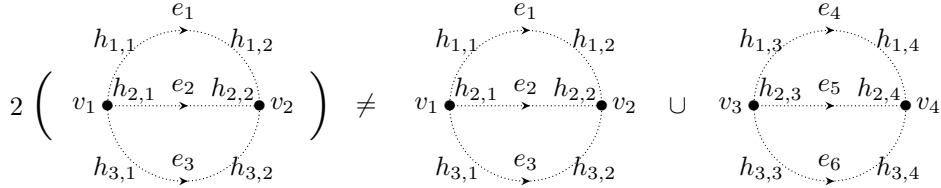
\end{definition}

Now we introduce an operator $\delta$ on $\mathcal{GC}_{\mathfrak{g}}$ as follows. Set
\begin{equation}
    \delta \Gamma = 
			\sum_{e=(ij): \text{non-self-loop edge}} \sigma(i,j) \cdot \Gamma / e,
\label{eq:6.2.5}
\end{equation}
where $e=(ij)$ denotes a non-self-loop edge connecting the vertex $i$ to the vertex $j$, $\Gamma / e$ means the decorated graph obtained from $\Gamma$ by contracting $e=(ij)$ to the original vertex $i$ then equipped with the consistent renumbering of edges and vertices and with the obvious information of insertion of one edge in place of the resulting  vertex, and the sign $\sigma(i,j)$ is defined as follows:

\begin{equation}
    \sigma(i,j) = \begin{cases}
    (-1)^{j} & \text{if $j > i$} \\
    (-1)^{i+1}  & \text{if $j < i$}.
    \end{cases}
    \label{eq:6.2.6}
\end{equation}

More concretely, the renumbering of $\Gamma/e$ is defined as follows. If $e=(ij)$ is the $k$-th edge of $\Gamma$, we renumber edges $e_l$ with $ k<l$ by letting them decrease by one. We renumber the vertices $v_l$ with $\max\{i,j\} \leq l$ by letting decrease by one and label the resulting vertex where the contraction has happened by $\min\{i,j\}$.

For a non-self-loop edge $e=(ij)$ connecting the vertex $i$ to the vertex $j$, the resulting vertex $i':=\min\{i,j\}$ by contracting $e=(ij)$ is attached with 
the equivariant homomorphism
\begin{equation}
    W_{i'}:  \otimes_{h\in h_{\Gamma/e}(i')}\mathfrak{g}_{h} \rightarrow  \mathbb{R}
\end{equation}
defined as follows. For defining $W_{i'}$, it is enough to define the corresponding oriented uni-trivalent tree $T_{i'}$ inserted at vertex $i'$ in $\Gamma/e$. Assume $T_i$, $T_j$ to be the inserted oriented uni-trivalent trees attached to vertex $i$ and $j$ respectively, then the inserted tree $T_{i'}$ is defined as the tree given by connecting $T_i$ and $T_j$ via the external edges corresponding to the edge $e=(ij)$. Note that the external edges of $T_{i'}$ are ordered according to the numberings on the edges and the directions of self-loops (if attached to $i'$), whose ordering is compatible with the ones of $T_i$ and $T_j$. The vertex-wise orientation on $T_{i'}$ is the one inherited from $T_i$ and $T_j$. Let $n_i$, $n_j$ denote the valencies of $i$, $j$ respectively, then the valency for this vertex $i'$ in $\Gamma/e$ is $n_{i'}:=n_i+n_j-2$, and the total edge number of $T_{i'}$ is $(2n_i-3)+(2n_j-3)-1=2n_{i'}-3$, this way, we confirm that $\Gamma/e$ with the above weight system $T_{i'}$ at $i'$ satisfies Definition \ref{def:6.3.1August}, i.e., $\Gamma/e\in \mathcal{GC}_{\mathfrak{g}}$.

\begin{remark}\label{rmk:6.4.4August}
    Assuming that $e=(ij)$ connecting $k$-th half-edge of $|h_{\Gamma}(i)|$ half-edges at $i$ and $l$-th half-edge of $|h_{\Gamma}(j)|$ half-edges at $j$, then the weight system $W_{i'}$ can be computed by (after re-order the tensor factors according to the ordering of half-edges) 
\begin{equation}
    W_{i'} = B_{k, |h_{\Gamma}(i)| +l} (W_i \otimes W_j) 
\end{equation}
where $B_{r,s}$ denotes the bilinear form $B$ acting on $r$-th and $s$-th components of tensor products $(\otimes_{h\in h_{\Gamma}(i)}\mathfrak{g}_{h})\otimes (\otimes_{h\in h_{\Gamma}(j)}\mathfrak{g}_{h}).$
\end{remark}
 
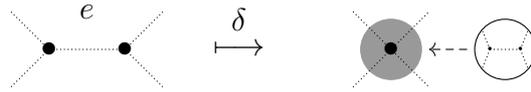
\begin{figure}[h]
 \centering
 \begin{tikzpicture}
 	\draw[densely dotted] (0,0) -- (1,0);
 	\draw[densely dotted] (0,0) -- (-0.5, 0.5);
 	\draw[densely dotted] (0,0) -- (-0.5, -0.5);
 	\draw[densely dotted] (1,0) -- (1.5, 0.5);
 	\draw[densely dotted] (1,0) -- (1.5, -0.5);
 	\node at (0.5, 0.5) {$e$};
  \node at (0,0) {$\bullet$};
  \node at (1,0) {$\bullet$};
  \draw[fill=black!40, black!40] (4.5, 0) circle (0.4cm);
   \node at (4.5, 0) {$\bullet$};
   \node at (2.5,0) {$\longmapsto$};
 	\node at (2.5,0.4) {$\delta$};
 	\draw[densely dotted] (4, -0.5) -- (5, 0.5);
 	\draw[densely dotted] (4, 0.5) -- (5, -0.5);
  \draw[densely dashed, ->] (5.5,0) -- (5,0);
  \begin{scope}[xshift=6cm]
    \draw (0,0) circle (0.4cm);
         \node at (0.2,0) {$\cdot$};
      \node at (-0.2,0) {$\cdot$};
    \draw[densely dotted] (-0.2, 0) -- (0.2, 0);
    \draw[densely dotted] (-0.2, 0) -- (-0.3,0.25);
    \draw[densely dotted] (-0.2, 0) -- (-0.3,-0.25);
    \draw[densely dotted] (0.2, 0) -- (0.3,0.25);
    \draw[densely dotted] (0.2, 0) -- (0.3,-0.25);
    \end{scope}
 \end{tikzpicture}
 \caption[The map $\delta$ for an internal edge]{The map $\delta$ for an internal edge $e$}
 \label{fig:part_delta}
\end{figure}

One simple example of the above contraction of one edge is illustrated as Fig. \ref{fig:part_delta}

Then, we have the following proposition analogous to \cite[Proposition 3.4]{BC2}. 
\begin{proposition}\label{prop:general_graph_complex}
The operator $\delta$ is a well-defined linear operator on $\mathcal{GC}_{\mathfrak{g}}$ and satisfies $\delta^2 = 0$. Moreover, for each $t\in\mathbb{Z}$, denoting by $\mathcal{GC}^t_{\mathfrak{g}}$  the subspace of $\mathcal{GC}_{\mathfrak{g}}$  spanned by the decorated graphs of degree $t$, we have
\begin{equation}
	\delta: \mathcal{GC}^t_{\mathfrak{g}} \rightarrow \mathcal{GC}^{t+1}_{\mathfrak{g}}.
\end{equation}
That is, the pair $(\bigoplus_{t} \mathcal{GC}^t_{\mathfrak{g}}, \delta)$ forms a complex.

If $\Gamma$ and $\Gamma'$ are two connected decorated graphs, then we have
\begin{equation}
\delta (\Gamma\cup \Gamma')=(\delta \Gamma)\cup \Gamma' +(-1)^{\mathrm{deg}(\Gamma)}\Gamma\cup (\delta \Gamma').
\label{eq:6.4.15Sept}
\end{equation}
\end{proposition}
\begin{proof}
By \eqref{eq:6.2.6} and Remark \ref{rmk:6.4.4August}, the well-definedness of $\delta$ and $\delta^2 =0$ follows from the same arguments as in the proof of \cite[Proposition 3.4]{BC2}, the existence of self-loops in the graphs does not produce any new obstacles. The identity \eqref{eq:6.4.15Sept} follows from the explicit computations for the operations $\delta$ and $\cup$.
\end{proof}

\begin{remark}\label{rk:new24}
    Now we can compare our graph complex for acyclic local systems with the Lie graph complex defined in \cite[\S 3]{MR2026331}. As mentioned at the beginning of \cite[\S 3]{MR2026331}, the objects in the Lie graph complex are finite graphs such that each vertex is decorated with a vertex-oriented trivalent tree, modulo AS and IHX relations. So their graphs are almost the same as ours in Definition \ref{def:6.3.1August}, apart from directions on edges and the equivariant weight $W_T$. In the main part of \cite[\S 3]{MR2026331}, they used a simplification for these graphs, the trivalent graphs equipped with a forest that contains all the vertices (considering one single vertex as a tree in the forest). From our definition, this can be obtained by inserting all the trees attached to each vertex into the graph. The grading used in \cite[\S 3]{MR2026331} is the total number of the trees in the given forest, their boundary map, equivalent to \eqref{eq:6.2.5}, will perform the connection of two different trees in the forest hence decreases the grading by $1$. However, to serve our purpose on the configuration space integral, we grade our decorated graphs according to the degree defined in \eqref{eq:orddegJuly}, and the boundary map increases the degree by $1$. For instance, in \cite[\S 3]{MR2026331}, a trivalent graph with the forest given as the set of all vertices belongs to the graph space graded by the total number of vertices, but in our case, it is of degree $0$.
\end{remark}

\subsection{Graph complexes of decorated graphs without self-loops}\label{sec:gc_acyclic}
The self-loops of the graph concerned here play a different role from the regular edges so we will investigate separately the subspaces of $\mathcal{GC}_{\mathfrak{g}}$ consisting of graphs without any self-loop and with at least one self-loop.

\begin{definition}
 Let $\mathcal{GC}_{\mathfrak{g}} = (\mathcal{GC}_{\mathfrak{g}}, \delta)$ be the graph complex defined in Subsection \ref{ss6.4complex}. Then, we similarly let $\mathcal{G}_{\mathfrak{g}}$ and $\mathcal{GC}'_{\mathfrak{g}}$ be the $\mathbb{Q}$-vector subspaces of $\mathcal{GC}_{\mathfrak{g}}$ spanned by the equivalent classes of decorated graphs, respectively, without self-loops and with at least one self-loop.

 These spaces of decorated graphs are bigraded by their order and degree.
For $n, t\in \mathbb{Z}$, let $\mathcal{G}^t_{\mathfrak{g}:n}$, $\mathcal{GC}^t_{\mathfrak{g}:n}$, $\mathcal{GC}^{\prime,t}_{\mathfrak{g}:n}$ denote the subspaces of $\mathcal{G}_{\mathfrak{g}}$, $\mathcal{GC}_{\mathfrak{g}}$, $\mathcal{GC}^{\prime}_{\mathfrak{g}}$, respectively, spanned by all the equivalent classes of decorated graphs with order $n$ and degree $t$.  
\end{definition}

\begin{figure}[h]
 \centering
 \begin{tikzpicture}
 \draw[fill=black!40, black!40] (3,0) circle (0.4cm);
    \draw[densely dotted] (0,0) circle (0.5cm);
    \draw[densely dotted] (-0.5,0) -- (0.5,0);
    \node at (-0.5,0) {$\bullet$};
    \node at (0.5,0) {$\bullet$};
    \node at (1.25,0) {$\longmapsto$};
    \node at (1.25,0.4) {$\delta$};
    \draw[densely dotted] (3,0.5) circle (0.5cm);
    \draw[densely dotted] (3,-0.5) circle (0.5cm);
    \node (p1) at (3,0) {$\bullet$};
    \node at (p1.90) {$v$};
    \node at (2,0) {$3$};
    \draw[densely dashed, ->] (4,0) -- (3.5,0);
    \begin{scope}[xshift=4.5cm]
    \draw (0,0) circle (0.4cm);
         \node at (0.2,0) {$\cdot$};
      \node at (-0.2,0) {$\cdot$};
    \draw[densely dotted] (-0.2, 0) -- (0.2, 0);
    \draw[densely dotted] (-0.2, 0) -- (-0.3,0.25);
    \draw[densely dotted] (-0.2, 0) -- (-0.3,-0.25);
    \draw[densely dotted] (0.2, 0) -- (0.3,0.25);
    \draw[densely dotted] (0.2, 0) -- (0.3,-0.25);
    \end{scope}
\end{tikzpicture}
 \caption[The action of $\delta$ on the theta graph]{The action of the operator $\delta$ on theta graph}
 \label{fig:6.2.2}
\end{figure}
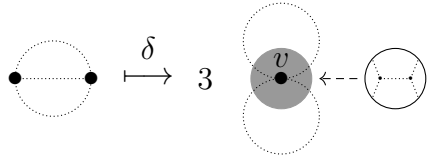

Note that $\delta$-action yields a self-loop when $\delta$ acts on a non-self-loop edge which is non-regular. Therefore, $\delta$ does not preserve the subspace $\mathcal{G}_{\mathfrak{g}}$. For example, Fig. \ref{fig:6.2.2} shows that $\delta$-action on the theta graph produces $3$ copies of the same figure-eight graphs as depicted. 

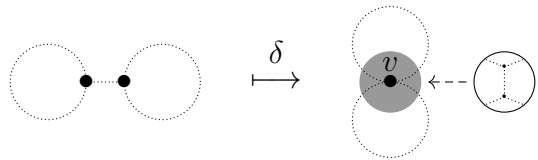
\begin{figure}[h]
 \centering
 \begin{tikzpicture}
 \begin{scope}[xshift=-0.5cm]
    \draw[densely dotted] (0,0) circle (0.5cm);
    \draw[densely dotted] (0.5,0) -- (1,0);
    \draw[densely dotted] (1.5,0) circle (0.5cm);
    \node at (0.5, 0) {$\bullet$};
    \node at (1, 0) {$\bullet$};
    \end{scope}
    \draw[fill=black!40, black!40] (4,0) circle (0.4cm);
    \node at (2.5,0) {$\longmapsto$};
    \node at (2.5,0.4) {$\delta$};
    \draw[densely dotted] (4,0.5) circle (0.5cm);
    \draw[densely dotted] (4,-0.5) circle (0.5cm);
    \node (p2) at (4,0) {$\bullet$};
    \node at (p2.90) {$v$};
    \draw[densely dashed, ->] (5,0) -- (4.5,0);
    \begin{scope}[xshift=5.5cm, rotate=90]
    \draw (0,0) circle (0.4cm);
         \node at (0.2,0) {$\cdot$};
      \node at (-0.2,0) {$\cdot$};
    \draw[densely dotted] (-0.2, 0) -- (0.2, 0);
    \draw[densely dotted] (-0.2, 0) -- (-0.3,0.25);
    \draw[densely dotted] (-0.2, 0) -- (-0.3,-0.25);
    \draw[densely dotted] (0.2, 0) -- (0.3,0.25);
    \draw[densely dotted] (0.2, 0) -- (0.3,-0.25);
    \end{scope}
\end{tikzpicture}
 \caption[The action of $\delta$ on dumbbell graph]{The action of the operator $\delta$ on dumbbell graph}
 \label{fig:6.2.3}
\end{figure}

On the other hand, $\delta$ preserves the subspace of decorated graphs with self-loops. In Fig. \ref{fig:6.2.3}, $\delta$-action on the dumbbell graph gives the figure-eight graph with the weight system at the unique vertex $v$ as given in Fig. \ref{fig:6.2.3}. 


Then, noting that $\delta$ preserves the order of a connected decorated graph, by Proprositin \ref{prop:general_graph_complex}, we conclude the following results. 
\begin{proposition}\label{prop:graph_complex}
\begin{enumerate}[(1)]
\item For each $n \in \mathbb{Z}$, the pair $(\bigoplus_{t} \mathcal{GC}^t_{\mathfrak{g}: n}, \delta)$ forms a complex.
\item  Taking the graphs always with at least one self-loop, $(\bigoplus_{t} \mathcal{GC}^{\prime, t}_{\mathfrak{g}: n}, \delta)$ form a subcomplex of $(\bigoplus_{t} \mathcal{GC}^t_{\mathfrak{g}: n}, \delta)$.
\end{enumerate}
\end{proposition}

By Proposition \ref{prop:graph_complex} (2), we can define a complex by
\begin{equation}
	(\bigoplus_{t} \mathcal{G}^{ t}_{\mathfrak{g}: n}, \delta^\sharp) = \Big(\bigoplus_{t} \big(\mathcal{GC}^t_{\mathfrak{g}: n}/\mathcal{GC}^{\prime, t}_{\mathfrak{g}: n}\big), \delta\Big).
 \label{eq:6.1.13paris}
\end{equation}
Then, the quotient complex $(\bigoplus_{t} \mathcal{G}^{ t}_{\mathfrak{g}: n}, \delta^\sharp) $ is the direct analog of the graph complex defined in \cite[Proposition 3.4]{BC2}, where the differential $\delta^\sharp$ acts only on a decorated graph $\Gamma$ without self-loops by
\begin{equation}
    \delta^\sharp \Gamma = \sum_{e=(ij): \text{regular edge}} \sigma(i,j) \cdot \Gamma / e
\end{equation}
where $e=(ij)$ denotes a regular edge connecting the vertex $i$ to the vertex $j$, and the sum is set to be zero if there is no regular edge in $\Gamma$.

Note that we only consider $t\geq 0$, $n\geq 1$. Fix a $n\in\mathbb{N}^\ast$, for $i\in\mathbb{N}$, let $H^i(\mathcal{GC}^\bullet_{\mathfrak{g}: n}, \delta)$ (resp. $H^i(\mathcal{G}^\bullet_{\mathfrak{g}: n}, \delta^\sharp)$, $H^i(\mathcal{GC}^{\prime,\bullet}_{\mathfrak{g}: n}, \delta)$) denote the $i$-th cohomology group of the complex $( \mathcal{GC}^\bullet_{\mathfrak{g}: n}, \delta)$ (resp. $(\mathcal{G}^\bullet_{\mathfrak{g}: n}, \delta^\sharp)$, $(\mathcal{GC}^{\prime, \bullet}_{\mathfrak{g}: n}, \delta)$). We have the following exact sequence
\begin{equation}
   0 \rightarrow H^0(\mathcal{GC}^{\prime,\bullet}_{\mathfrak{g}: n}, \delta) \rightarrow  H^0(\mathcal{GC}^\bullet_{\mathfrak{g}: n}, \delta) \rightarrow H^0(\mathcal{G}^\bullet_{\mathfrak{g}: n}, \delta^\sharp)\rightarrow  H^1(\mathcal{GC}^{\prime,\bullet}_{\mathfrak{g}: n}, \delta) \rightarrow \cdots
\label{eq:6.5.3sept}
\end{equation}

\begin{definition}
An element in $\Ker \delta \subset \mathcal{GC}_{\mathfrak{g}}$ or in $\Ker \delta^\sharp \subset \mathcal{G}_{\mathfrak{g}}$ is called a {\it graph cocycle} in the respective graph complexes. In particular, the graph cocycles of degree $0$ are exactly the elements in $H^0(\mathcal{GC}^\bullet_{\mathfrak{g}: n}, \delta)$, $H^0(\mathcal{G}^\bullet_{\mathfrak{g}: n}, \delta^\sharp)$.
\end{definition}

In fact, we are mainly concerned with graph cocycles (or simply cocycles) of degree $0$. Before we proceed to see some examples, we give several easy facts.
\begin{itemize}
    \item For $n=1$, we have $\mathcal{GC}^{\prime,0}_{\mathfrak{g}: 1}=\mathbb{Q}\cdot \text{\,dumbbell}$, $\mathcal{GC}^{\prime,1}_{\mathfrak{g}: 1}= \mathbb{Q}\cdot \text{\,figure-eight}$, which are $1$-dimensional. We can conclude
    \begin{equation}
        H^0(\mathcal{GC}^{\prime,\bullet}_{\mathfrak{g}: 1},\delta)=H^1(\mathcal{GC}^{\prime,\bullet}_{\mathfrak{g}: 1},\delta)=0.
    \end{equation}
    So that $H^0(\mathcal{GC}^\bullet_{\mathfrak{g}: 1}, \delta)=H^0(\mathcal{G}^\bullet_{\mathfrak{g}: 1}, \delta^\sharp)$, and they are also $1$-dimensional (over $\mathbb{Q}$).
    \item Any cocycle in $\mathcal{GC}_{\mathfrak{g}}$, by taking its quotient class or equivalently, by removing all the terms including self-loops, gives a cocycle in $\mathcal{G}_{\mathfrak{g}}$.
    \item If $\Gamma$, $\Gamma'$ are two cocycles, then so is $\Gamma\cup \Gamma'$. So that the spaces of cocycles carry the induced structure of graded commutative algebra.
\end{itemize}

\begin{example}[2-loop cocycles]\label{ex:6.2.1}
\begin{enumerate}[(1)]
\item In the graph complex $\mathcal{GC}_{\mathfrak{g}}$, the following linear combination gives a cocycle of degree $0$ with $2$-loops:
\begin{equation}
\Theta  - \frac{3}{2}\text{O--O}
\end{equation}
where $\Theta$  and $\text{O--O}$ decorated as in Fig. \ref{fig:6.1.theta}. In fact, $H^0(\mathcal{GC}^\bullet_{\mathfrak{g}: 1}, \delta)$ is exactly spanned by the above cocycle over $\mathbb{Q}$. Here we need to put coefficient $\frac{3}{2}$ instead of $3$ because the factor $2$ appeared in \eqref{eq:self-loop_jacobi}, following from the internal IHX relation.
\item  In the graph complex $\mathcal{G}_{\mathfrak{g}}$, the $\Theta$-graph itself gives a cocycle since $\delta$-action on it yields graph with self-loops which is defined to be zero in the quotient space $\mathcal{G}_{\mathfrak{g}}\simeq \mathcal{GC}_{\mathfrak{g}}/\mathcal{GC}^{\prime}_{\mathfrak{g}}$. So that $H^0(\mathcal{G}^\bullet_{\mathfrak{g}: 1}, \delta^\sharp)$ is spanned by the $\Theta$-graph.

\item The decorated graph $\Theta \cup \Theta$ is a nontrivial cocycle in  $\mathcal{G}^0_{\mathfrak{g}: 2}$, but it has $4$ loops in topological sense.

\end{enumerate}
\end{example}

\begin{example}[3-loop cocycles]\label{ex:3loop}
In \cite[Example 4.6]{BC}, for the trivial local system on a framed homology $3$-sphere, Bott--Cattaneo gave an example of cocycle with degree $0$ and order $2$ (hence with $3$ loops) by the following linear combination:
    \begin{equation}
    \Gamma'  = \frac{1}{12} \Gamma_1 + \frac{1}{4} \Gamma_2,
    \end{equation}
    where $\Gamma_1$ and $\Gamma_2$ are given in Fig. \ref{fig:6.2.4} (without the weight systems induced from Lie algebra $\mathfrak{g}$ or the numbering on the edges).

As an element in $\mathcal{G}^0_{\mathfrak{g}:2}$, the following linear combination is a cocycle (the coefficient of $\Gamma_2$ has been changed to $-\frac{1}{8}$):
 \begin{equation}
    \Gamma  = \frac{1}{12} \Gamma_1 - \frac{1}{8} \Gamma_2,
    \label{eq:6.5.7new}
    \end{equation}
    where $\Gamma_1$ and $\Gamma_2$ are given in Fig. \ref{fig:6.2.4} with all the decorations.
\begin{figure}[h]
 \centering
    \begin{tikzpicture}[scale=1.2, font=\footnotesize]
    \begin{scope}
    \draw[densely dotted] (0,0) circle (1.0cm);
    \draw[densely dotted, ->, >=stealth] (1.414/2, 1.414/2) -- (-1.414/4, -1.414/4);
    \draw[densely dotted] (1.414/4, 1.414/4) -- (-1.414/2, -1.414/2);
    \draw[densely dotted, ->, >=stealth] (-1.414/2, 1.414/2) -- (1.414/4, -1.414/4);
    \draw[densely dotted] (-1.414/4, 1.414/4) -- (1.414/2, -1.414/2);
    \node at (1.414/2, 1.414/2) {$\bullet$};
    \node at (-1.414/2, 1.414/2) {$\bullet$};
    \node at (1.414/2, -1.414/2) {$\bullet$};
    \node at (-1.414/2, -1.414/2) {$\bullet$};
    \node at (0,-1.7) {$\Gamma_1$};
    \draw[->,>=stealth] (1,-0.05)--(1,-0.051);
    \draw[->,>=stealth] (-1,0.05)--(-1,0.051);
    \draw[->,>=stealth] (0.05,1)--(0.051,1);
    \draw[->,>=stealth] (-0.05,-1)--(-0.051,-1);
    \node at (-0.9, 0.9) {$v_1$};
    \node at (0.9, 0.9) {$v_2$};
    \node at (0.9, -0.9) {$v_3$};
    \node at (-0.9, -0.9) {$v_4$};
    \node at (0, 1.2) {$e_1$};
    \node at (1.25, 0) {$e_2$};
    \node at (0, -1.2) {$e_3$};
    \node at (-1.25, 0) {$e_4$};
    \node at (-0.4, 0.2) {$e_5$};
    \node at (0.4, 0.2) {$e_6$};
 \end{scope}
    \begin{scope}[xshift=5cm]
         \draw[densely dotted] (0,0) circle (1.0cm);
    \draw[densely dotted] (0,0) circle (0.5cm);
    \draw[densely dotted, ->, >=stealth] (0.5,0) -- (0.8, 0);
    \draw[densely dotted] (0.8,0) -- (1, 0);
    \draw[densely dotted, ->, >=stealth] (-1,0) -- (-0.7, 0);
    \draw[densely dotted] (-0.7,0) -- (-0.5, 0);
    \node at (0.5, 0) {$\bullet$};
    \node at (-0.5, 0) {$\bullet$};
     \node at (1, 0) {$\bullet$};
    \node at (-1, 0) {$\bullet$};
    \draw[->,>=stealth] (0.05,1)--(0.051,1);
     \draw[->,>=stealth] (0.05,0.5)--(0.051,0.5);
     \draw[->,>=stealth] (0.05,-0.5)--(0.051,-0.5);
    \draw[->,>=stealth] (0.05,-1)--(0.051,-1);
    \node at (0,-1.7) {$\Gamma_2$};
    \node at (-1.3, 0) {$v_1$};
     \node at (-0.2, 0) {$v_2$};
     \node at (0.2, 0) {$v_3$};
     \node at (1.3, 0) {$v_4$};
 \node at (-0.75, 0.2) {$e_1$};
    \node at (0.75, 0.2) {$e_2$};
    \node at (0, 0.7) {$e_3$};
    \node at (0, -0.7) {$e_4$};
    \node at (0, 1.2) {$e_5$};
    \node at (0, -1.2) {$e_6$};
    \end{scope}
\end{tikzpicture}	
 \caption[Decorated trivalent 3-loop graphs without self-loop]{Two examples of decorated trivalent 3-loop graphs without self-loop}
 \label{fig:6.2.4}
\end{figure}
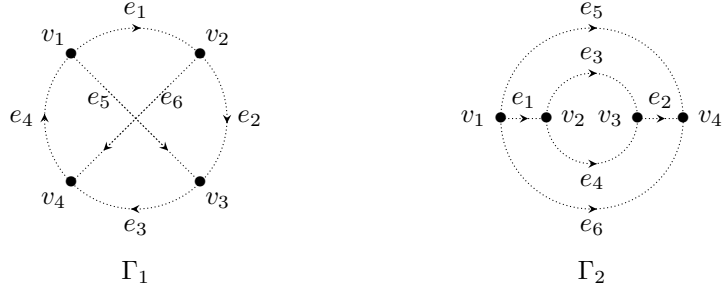

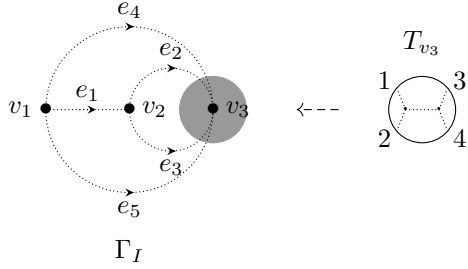
\begin{figure}[h]
 \centering
    \begin{tikzpicture}[scale=1.1, font=\footnotesize]
  \begin{scope}
 \draw[fill=black!40, black!40] (1, 0) circle (0.4cm);
       \draw[densely dotted] (0,0) circle (1.0cm);
    \draw[densely dotted] (0.5,0) circle (0.5cm);
    \draw[densely dotted, ->, >=stealth] (-0.6,0) -- (-0.4, 0);
    \draw[densely dotted] (-1,0) -- (0, 0);
    \node at (0, 0) {$\bullet$};
     \node at (1, 0) {$\bullet$};
    \node at (-1, 0) {$\bullet$};
    \draw[->,>=stealth] (0.05,1)--(0.051,1);
     \draw[->,>=stealth] (0.55,0.5)--(0.551,0.5);
     \draw[->,>=stealth] (0.55,-0.5)--(0.551,-0.5);
    \draw[->,>=stealth] (0.05,-1)--(0.051,-1);
    \node at (0,-1.7) {$\Gamma_I$};
    \node at (-1.3, 0) {$v_1$};
     \node at (0.3, -0) {$v_2$};
   \node at (1.3, 0) {$v_3$};
 \node at (-0.5, 0.2) {$e_1$};
  \node at (0.5, 0.7) {$e_2$};
    \node at (0.5, -0.7) {$e_3$};
    \node at (0, 1.2) {$e_4$};
    \node at (0, -1.2) {$e_5$};
 \draw[densely dashed, ->] (2.5,0) -- (2,0);

    \end{scope}

    \begin{scope}[xshift=3.5cm] 
    \draw (0,0) circle (0.4cm);
      \node at (0.2,0) {$\cdot$};
      \node at (-0.2,0) {$\cdot$};
    \draw[densely dotted] (-0.2, 0) -- (0.2, 0);
    \draw[densely dotted] (-0.2, 0) -- (-0.3,0.25);
        \node at (-0.45,0.35) {$1$};
    \draw[densely dotted] (-0.2, 0) -- (-0.3,-0.25);
    \node at (-0.45,-0.35) {$2$};
    \draw[densely dotted] (0.2, 0) -- (0.3,0.25);
     \node at (0.45,0.35) {$3$};
    \draw[densely dotted] (0.2, 0) -- (0.3,-0.25);
     \node at (0.45,-0.35) {$4$};
    \node at (0, 0.8) {$T_{v_3}$};
    \end{scope}

\end{tikzpicture}	
 \caption[Graph $\Gamma_I$ in the computation of $\delta^\sharp \Gamma$]{Graph $\Gamma_I$ in the computation of $\delta^\sharp \Gamma$}
 \label{fig:6.2.5sept}
\end{figure}

\begin{figure}[h]
\captionsetup{margin=2cm}
 \centering
    \begin{tikzpicture}
    \begin{scope}
    \foreach \x in {0, 2, 4}{
     \draw[densely dotted] (\x,0) circle (0.5cm);
    }
   \draw[densely dotted] (0.5, 0) -- (1.5,0);
   \draw[densely dotted] (2.5, 0) -- (3.5,0);
   \foreach \x in {0.5, 1.5, 2.5, 3.5}{
    \node at (\x, 0) {$\bullet$};
   }
   \Darr{1}{0};
   \Darr{0}{0.5};
   \Darr{2}{0.5};
   \Darr{4}{0.5};
    \Darr{3}{0};
   \node at (2,-1.5) {$\Gamma_3$};
 \end{scope}
 \begin{scope}[xshift=6cm]
 \foreach \x in {0, 2}{
     \draw[densely dotted] (\x,0) circle (0.5cm);
    }
   \draw[densely dotted] (0.5, 0) -- (1.5,0);
   \draw[densely dotted] (2, 0.5) -- (2,-0.5);
    \node at (1,-1.5) {$\Gamma_4$};
    \foreach \x in {0.5, 1.5}{
     \node at (\x, 0) {$\bullet$};
    }
    \foreach \x in {0.5, -0.5}{
     \node at (2, \x) {$\bullet$};
    }
   \Darr{1}{0};
   \Darr{0}{0.5};
   \draw[->, >=stealth] (1.75, 0.45) -- (1.75+0.001,0.45+0.001);
    \draw[->, >=stealth] (1.75, -0.45) -- (1.75+0.001,-0.45-0.001);
    \draw[->, >=stealth] (2, 0) -- (2, 0+0.001);
    \draw[->, >=stealth] (2.5, 0) -- (2.5, 0-0.001);
   \end{scope}
   \begin{scope}[xshift=11cm]
   \foreach \i in {90, 210, 330}{%
\draw[densely dotted] (0,0) --  (\i:0.5cm);
\draw[densely dotted]  (\i:1cm) circle (0.5cm);
\draw[->, >=stealth] (\i:0.35cm) -- (\i:0.35cm+0.001cm);
\node at (\i:0.5cm) {$\bullet$};
}
\Darr{-0.8}{-0.99};
\Darr{0.8}{-0.99};
\Darr{0.0}{1.5};
\node at (0,0) {$\bullet$};
 \node at (0,-1.5) {$\Gamma_5$};
   \end{scope}
\end{tikzpicture}	
 \caption[Decorated trivalent 3-loop graphs with self-loops]{Three examples of decorated trivalent 3-loop graphs with self-loops.}
 \label{fig:6.2.4new}
\end{figure}
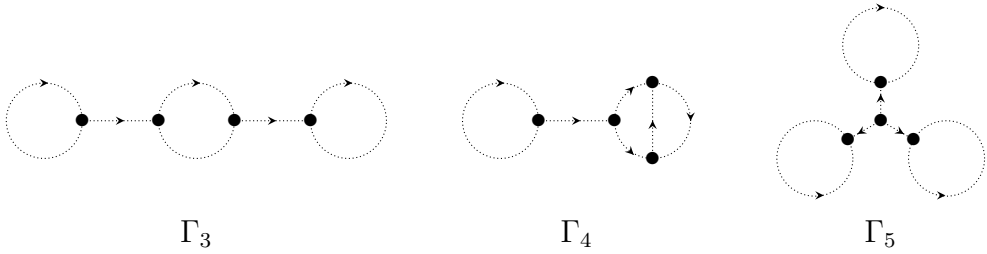

To show $\Gamma$ being a cocycle in $\mathcal{G}^0_{\mathfrak{g}:2}$, we need to use the identity in Fig. \ref{fig:Sept8} and the following equations in $\mathcal{G}^0_{\mathfrak{g}:2}$, 
\begin{equation}
\delta^\sharp \Gamma_1=3 \Gamma_I,\; \delta^\sharp \Gamma_2=2 \Gamma_I.
\label{eq:6.5.8sept}
\end{equation}

However, if we view $\Gamma$ as an element in $\mathcal{GC}^0_{\mathfrak{g}:2}$, then it is no longer a cocycle since the $\delta$-action on $\Gamma_2$ produces terms with self-loops which do not cancel each other out.

In fact, we have three more connected trivalent graphs of order $2$, which are displayed in Fig. \ref{fig:6.2.4new} (one needs to assign numberings on their vertices or edges to make a decorated graph). We can add these graphs with self-loops into $\Gamma$ to obtain a cocycle $\widetilde{\Gamma}$ in $\mathcal{GC}^0_{\mathfrak{g}:2}$:
\begin{equation}
\widetilde{\Gamma}=\Gamma + \mathbb{Q}\text{-linear combination of } \Gamma_3, \Gamma_4 \text{ and }\Gamma_5.
\label{eq:6.5.9new}
\end{equation}
A precise formula for $\widetilde{\Gamma}$ can be deduced from Example \ref{ex:6.6.7kks}. 
\end{example}


\subsection{Topological trivalent graph and Chern--Simons perturbation theory}\label{ss6.6kk}

Now we focus on the connected trivalent graphs. 
For integer $n\geq 1$, let $\mathfrak{G}$ be a topological trivalent graph with $2n$ vertices and $3n$ edges, or equivalently, consider a set $h(\mathfrak{G})$ of $6n$ (abstract) half-edges, then a trivalent graph means the couple of partitions of  $h(\mathfrak{G})$:
\begin{itemize}
    \item A partition into pairs of half-edges which we call edges.
   \item A partition into sets of cardinality (=valency) 3 which we call vertices.
\end{itemize}
If any two vertices can be connected by a consecutive path of edges (any neighboring edges have only one common vertex), then we call the graph to be connected. If $\mathfrak{G}_1$, $\mathfrak{G}_2$ are two topological trivalent graphs, they are called equivalent to each other if there is a bijection between $h(\mathfrak{G}_1)$ and $h(\mathfrak{G}_2)$ which maps the couple of partitions of $h(\mathfrak{G}_1)$ to the ones of $h(\mathfrak{G}_2)$. We will always identify the equivalent graphs as the same one. Note that the order of $\mathfrak{G}$ is defined by the same formula as in \eqref{eq:orddegJuly} (since $\mathfrak{G}$ is trivalent, it always has degree $0$).

\begin{definition}[Automorphism group of topological trivalent graph]\label{def:6.6.3}
    let $\mathfrak{G}$ be a topological trivalent graph with $2n$ vertices and $3n$ edges, then an automorphism of $\mathfrak{G}$ is an element of the permutation group of $h(\mathfrak{G})$ which preserves both partitions of $h(\mathfrak{G})$ for the edges and vertices of $\mathfrak{G}$, we denote the group of all automorphism of $\mathfrak{G}$ by $\mathrm{Aut}(\mathfrak{G})$.
\end{definition}

\begin{remark}
    If $\mathfrak{G}$ is a topological trivalent graph that is not connected, then by our definition of $\mathrm{Aut}(\mathfrak{G})$, its action always preserves the non-equivalent connected components of $\mathfrak{G}$. For example, suppose that $\mathfrak{G}_1$, $\mathfrak{G}_2$ are two connected topological trivalent graphs, then
   \begin{equation}
        \mathrm{Aut}(\mathfrak{G}_1\cup \mathfrak{G}_2)=\begin{cases}
           \mathrm{Aut}(\mathfrak{G}_1)\times \mathrm{Aut}( \mathfrak{G}_2)\rtimes \mathbb{Z}_2,  &\text{ if }\mathfrak{G}_1 =\mathfrak{G}_2\neq \emptyset;\\
           \mathrm{Aut}(\mathfrak{G}_1)\times \mathrm{Aut}( \mathfrak{G}_2),   &\text{ if else}.
           \end{cases}
    \end{equation}
\end{remark}

\begin{definition}[Relative orientation]\label{lm:6.6.2kk}
    Suppose that $\Gamma$ is a connected decorated trivalent graph, and let $h(\Gamma)$ denote the set of all half-edges of $\Gamma$. Then we have the induced vertex orientation and the induced edge orientation on the same set $h(\Gamma)$, the relative orientation of $\Gamma$, denoted by $\mathrm{or}_\Gamma\in \{\pm 1\}$, is defined as the sign of the permutation that maps the edge orientation to the vertex orientation.
\end{definition}

The following Lemma \ref{lm:6.6.2} is an analog of \cite[Corollary 1]{MR2026331} and \cite[\S 3]{AS}.
\begin{lemma}\label{lm:6.6.2}
    Let $\Gamma$ be a connected decorated trivalent graph, let $\Gamma'$ be the decorated trivalent graph given by the same underlying topological graph as of $\Gamma$ but with different numberings on the edges and vertices, then we have the identity in $\mathcal{GC}^0_{\mathfrak{g}}\;$,
    \begin{equation}
        \mathrm{or}_\Gamma \cdot \Gamma = \mathrm{or}_{\Gamma'} \cdot \Gamma'.
    \end{equation}
     Let $\Gamma_1$, $\Gamma_2$ be two connected decorated trivalent graphs, then we have
        \begin{equation}
        \mathrm{or}_{\Gamma_1\cup \Gamma_2}  = \mathrm{or}_{\Gamma_1} \cdot \mathrm{or}_{\Gamma_2}. 
    \end{equation}
\end{lemma}
\begin{proof}
    It is enough to check that the term $\mathrm{or}_{\Gamma}$ cancels the sign change coming from the numbering change of vertices, edges, and direction reversals of edges. If we permute the numbering on $v(\Gamma)$ by a permutation of order $p$, then edge orientation of $h(\Gamma)$ is fixed but vertex orientation of $h(\Gamma)$ gives sign change of $(-1)^p$ so that $\mathrm{or}_{\Gamma}$ also change sign by $(-1)^p$. If we change the direction of a non-self-loop edge $e=(ij)$, then the edge orientation of $h(\Gamma)$ differs by $(-1)$ under this direction reversal whereas the vertex orientation of $h(\Gamma)$ is fixed. If we change the direction of a self-loop edge $e=(ii)$ incident to vertex $i$, then both edge orientation and vertex orientation of $h(\Gamma)$ differ by $(-1)$ under this change, hence the total change is $1$. Similarly, if numbering change of $e(\Gamma)$ gives rise to $k$-times changes of cyclic ordering at trivalent vertices, then edge orientation of $h(\Gamma)$ is fixed but vertex orientation of $h(\Gamma)$ changes by the same manner. 
    
    Therefore, $\mathrm{or}_{\Gamma}$ cancels any sign change from permutations of numberings of edges and vertices, and direction reversals of edges.
\end{proof}

The above lemma indicates that for decorated trivalent graphs, the underlying topological graph uniquely determines its equivalence class up to a sign. Let $\mathcal{TG}_n$ denote the vector space spanned by all the topological trivalent graphs with $2n$ vertices over $\mathbb{Q}$. We consider the linear map
\begin{equation}\label{eq:topg_to_decg}
    \Psi_n: \mathcal{TG}_n\rightarrow \mathcal{GC}^0_{\mathfrak{g}: n}, \;\;\mathfrak{G} \mapsto \mathrm{or}_{\Gamma(\mathfrak{G})} \cdot\Gamma(\mathfrak{G})
\end{equation}
where $\Gamma(\mathfrak{G})$ is any decorated trivalent graph whose underlying topological graph is $\mathfrak{G}$.

As a consequence of Lemma \ref{lm:6.6.2}, we have the following result.
\begin{corollary}\label{Cor:4.6.5-july}
For each $n\geq 1$, the linear map $\Psi_n$ is an isomorphism of finite-dimensional vector spaces.
\end{corollary}

Now Let us consider the generating series of perturbative invariants for a framed closed $3$-manifold defined from the perturbative Chern--Simons theory (see \cite[Section 2]{Ko}, \cite{AS, AS2}), which, in terms of the trivalent graph, is formally given by
\begin{equation}
    \sum_{n\geq 0}\hbar^n \sum_{\substack{\text{trivalent graph }\mathfrak{G}\\ \text{ of order } n}} \frac{1}{|\mathrm{Aut}(\mathfrak{G})|}\mathfrak{G}=\exp\Big(\sum_{n\geq 1}\hbar^n \sum_{\substack{\text{connected trivalent graph }\mathfrak{G}\\ \text{ of order } n}} \frac{1}{|\mathrm{Aut}(\mathfrak{G})|}\mathfrak{G}\Big),
\end{equation}
where for $n=0$ we take $\mathfrak{G}=\emptyset$ viewed as a unit element, $|\mathrm{Aut}(\mathfrak{G})|=1$, and the multiplication of the topological graphs is given by the disjoint union $\cup$ (it is commutative).

\begin{proposition}\label{prop:CSseries}
    \begin{enumerate}
        \item[(1)] For each integer $n$ with $n\geq 1$, there is a cocycle of order $n$ in $\mathcal{GC}^0_{\mathfrak{g}:n}$ given as the form
\begin{equation}
    \sum_{\text{connected } \mathfrak{G}} \frac{1}{|\Aut(\mathfrak{G})|}  \Psi_n(\mathfrak{G}) \in H^0(\mathcal{GC}^\bullet_{\mathfrak{g}:n},\delta),
\label{eq:6.6.6new}
\end{equation}
where the sum runs over all the connected topological trivalent graph $\mathfrak{G}$ of order $n$.
\item[(2)] For each integer $n$ with $n\geq 1$, there is a cocycle of order $n$ in $\mathcal{G}^0_{\mathfrak{g}:n}$ given as the form
\begin{equation}
\sum_{\substack{\text{connected } \mathfrak{G} \\ \text{without self-loops}}} \frac{1}{|\Aut(\mathfrak{G})|}  \Psi_n(\mathfrak{G}) \in H^0(\mathcal{G}^\bullet_{\mathfrak{g}:n},\delta^\sharp),
    \label{eq:6.6.7new}
\end{equation}
where the sum runs over all the connected topological trivalent graph $\mathfrak{G}$ without self-loops and of order $n$.
    \end{enumerate}
\end{proposition}
\begin{proof}
Fix a partition $V$ of $(6n)$ half-edges into the sets of cardinality $3$ (viewed as vertices), let $P_{n}$ denote the set of partitions of this set of $(6n)$ half-edges into pairs. Then, consider a surjective map
   \begin{equation}
       \pi_n: P_{n} \rightarrow \mathcal{TG}_n,
   \end{equation}
   which sends a partition $E$ in $P_n$ to the topological graph $\mathfrak{G}(E,V)$ given by the equivalent class of the couple of partitions $(E,V)$.

   For a topological graph $\mathfrak{G}\in \mathcal{TG}_n$, let $G_{h(\mathfrak{G})}$ be the permutation group of the set $h(\mathfrak{G})$ of half-edges of $\mathfrak{G}$, and let $G_{v(\mathfrak{G})}, G_{e(\mathfrak{G})} \subset G_{h(\mathfrak{G})}$ denote the subgroups preserving the partitions of $h(\mathfrak{G})$ for the vertices and edges of $\mathfrak{G}$ respectively. 
   With these notations, we get
   \begin{equation}
       |\pi_n^{-1}(\mathfrak{G})| = \frac{|G_{v(\mathfrak{G})}|}{|\Aut(\mathfrak{G})|}
   \end{equation}

Note that if $\mathfrak{G}$ is connected, then all partitions of edges $E\in \pi_n^{-1}(\mathfrak{G})$ together with $V$ defines a connected trivalent graph. Then the sum in \eqref{eq:6.6.6new} can be written as

   \begin{equation}
   \begin{split}
      &   \sum_{\text{ connected } \mathfrak{G}} \frac{1}{|\Aut(\mathfrak{G})|}  \Psi_n(\mathfrak{G})\\
       &= \sum_{\text{ connected } \mathfrak{G}} \frac{|\pi_n^{-1}(\mathfrak{G})|}{|G_{v(\mathfrak{G})}|}\Psi_n(\mathfrak{G})\\
       & =\sum_{\text{ connected } \mathfrak{G}} \frac{|\pi_n^{-1}(\mathfrak{G})|}{ (3!)^{|v(\mathfrak{G})|} |v(\mathfrak{G})|!} \Psi_n(\mathfrak{G})\\
       & =   \frac{1}{(3!)^{2n} (2n)!} \sum_{E \in P_n \,;\, \text{ connected}} \Psi_n(\mathfrak{G}(E,V))
   \end{split}
   \label{eq:4.6.10-july-24}
   \end{equation}
   where in the last summation $E$ runs over all the partitions in $P_n$ such that the resulting graphs $\mathfrak{G}(E,V)$  are connected.

   
   We focus on a vertex with valency 4 after the $\delta$-action on the last line of \eqref{eq:4.6.10-july-24}, that is, we focus on a pair of half-edges to be contracted in the $\delta$-action, so that the other $4$ half-edges in these two associated vertices determine a vertex with valency $4$. Then, there are essentially only 3 possible ways to insert weight systems at this vertex, that is, IHX type graphs as in Fig. \ref{fig:insertion}, in this way we find back the connected trivalent graphs before the $\delta$-action. The insertions at each vertex with valency $4$ produce $4$ triplets of partitions $E$ in $P_n$, each triplet exactly consists of compatible $I,H,X$-insertions at this vertex (hence in total $12$ different partitions in $P_n$). 
   
   Note that the $X$ graph in Fig. \ref{fig:Jacobi} has opposite cyclic ordering from our convention, but one sees that this sign emerges from $\mathrm{or}_{\Gamma}$ in $\Psi_n$ by direct computation for these local graphs. Therefore, the sign of these three weight systems is compatible with Jacobi identity \eqref{eq:Jacobi-indentity}. Finally, the $\delta$-action on the last line of \eqref{eq:4.6.10-july-24} has to be $0$ by the vertex-wise IHX relation, so the conclusion (1) holds. The conclusion (2) is immediate from (1) by removing graphs with self-loops. This completes the proof.
   \end{proof}


\begin{example}\label{exm:6.6.7sss}
    In particular, considering the order-1 part in Proposition \ref{prop:CSseries}, we can recover Cattaneo--Shimizu's 2-loop term \eqref{eq:3.0.4}. Indeed, let $\Theta$ and $\text{O--O}$ be the theta graph and dumbbell graph decorated as Fig. \ref{fig:6.1.theta} (they are the only connected trivalent graphs of order $1$), then $\mathrm{or}_{\Theta} = -1$ and $\mathrm{or}_{\text{O--O}}=1$. For their underlying topological graphs, we have $|\mathrm{Aut}(\Theta)|=12$, $|\mathrm{Aut}(\text{O--O})|=8$. Thus, we get a cocycle in $\mathcal{GC}^0_{\mathfrak{g}:1}$,
    \begin{equation}
        \Gamma=\frac{1}{12} \mathrm{or}_{\Theta} \cdot \Theta + \frac{1}{8} \mathrm{or}_{\text{O--O}} \cdot\text{O--O} = - \frac{1}{12}\Theta  + \frac{1}{8} \text{O--O}.
    \end{equation}
    Applying the linear map $(-12) Z_{\Gamma}(\omega)$, we obtain $Z_1(\rho; \omega, \xi)$ in \eqref{eq:3.0.4}.
    
\end{example}

\begin{example}\label{ex:6.6.7kks}
    For the connected topological trivalent graph of order $2$, there are two cases without self-loops as given in Fig. \ref{fig:6.2.4}, and there are $3$ other cases with at least one self-loop edge given in Fig. \ref{fig:6.2.4new}. Let $\mathfrak{G}_j$, $j=1,\ldots, 5$, denote the respective underlying topological graphs of $\Gamma_j$, $j=1,\ldots, 5$ in Fig. \ref{fig:6.2.4} \& Fig. \ref{fig:6.2.4new}, then we have 
    \begin{equation}
    \begin{split}
    & |\mathrm{Aut}(\mathfrak{G}_1)|=24, \; |\mathrm{Aut}(\mathfrak{G}_2)|=16, \\
        &|\mathrm{Aut}(\mathfrak{G}_3)|=16, \;|\mathrm{Aut}(\mathfrak{G}_4)|=8, \;|\mathrm{Aut}(\mathfrak{G}_5)|=48.
    \end{split}
    \end{equation}
    Meanwhile, we have $\mathrm{or}_{\Gamma_1} = -1$ and $\mathrm{or}_{\Gamma_2}=1$, this way, from \eqref{eq:6.6.7new} we get a cocycle (without self-loops) in $\mathcal{G}^0_{\mathfrak{g}:2}$, which is proportional to the one defined in \eqref{eq:6.5.7new}. If we include the other $3$ cases with self-loops ($\Gamma_3$, $\Gamma_4$, $\Gamma_5$), we can work out explicitly a cocycle $\Gamma'$ mentioned in \eqref{eq:6.5.9new}.
\end{example}


\section{Higher-loop integral invariants for $3$-manifolds}\label{section6}
This section aims to define the higher-loop integral invariants associated with graph cocycles in the complex of decorated graphs. In general, the higher-loop integral invariants that extend Cattaneo--Shimizu's $Z_1$-invariant (see \cite{CS} and Definition \ref{defn:3.4.1-24}) require the graph complex with self-loops. Using the different graph complexes introduced in Section \ref{section:graph} and the adapted propagators as in Definition \ref{def:6.1.1}, we show that the map $Z(M,\rho,[f])$ in Theorem \ref{thm:CC} factors through a quotient graph complex that coincides with the one in Bott--Cattaneo \cite{BC2} without self-loops. This means that when an acyclic local system is given by $\rho:\pi_1(M) \rightarrow G \overset{\Ad}{\longrightarrow} \Aut(\mathfrak{g})$, the graph complex without self-loops is enough to define the integral invariants. 

In this section, a framing $f$ of $M$ and an orientation $o(M)$ are always fixed, and we always assume the local system $E_\rho$ to be acyclic.

\subsection{Integral invariants of higher order associated to acyclic local systems}

In this subsection, we study graph cocycle invariants of a framed closed 3-manifold with acyclic local system associated with a representation $\pi_1(M) \rightarrow G \overset{\Ad}{\longrightarrow} \Aut(\mathfrak{g})$.


Assume the local system $E_\rho$ to be acyclic.  Given a propagator $\omega$ as defined in Subsection \ref{ss4.3sept}, let us define a $\mathbb{Q}$-linear map $Z_{-}(\omega)$ on $\mathcal{GC}^0_{\mathfrak{g}}$. 

At first, associated to each edge $e\in e(\Gamma)$ of a decorated trivalent graph $\Gamma$ (hence of degree $0$), we define the 2-form $\omega_{e}$ on $C_{2n}(M)$ as follows:
\begin{equation}\label{eq:6.2.2june}
	\omega_{e} = \begin{cases}
 p_{ij}^{\ast} \omega 	& \text{if $e=(ij)$ with $i \neq j$},\\
  q^{\ast} p_i^{\ast} \xi & \text{if $e=(ii)$ is a directed self-loop},
 \end{cases}
\end{equation}
where $p_{ij} : C_{2n}(M) \rightarrow C_2(M)$ is the natural projection map induced by $M^{2n} \rightarrow M\times M$ which sends $(x_1,\ldots, x_{2n}) \mapsto (x_i,x_j)$, $q :C_{2n}(M) \rightarrow M^{2n}$ is the blow-down map (by abuse of notation), and $p_i: M^{2n} \rightarrow M$ is the natural $i$-th projection map. Note that when $e=(ij)$, $i\neq j$, the coefficient of form $\omega_{e}$ is in $p_i^\ast E_\rho \otimes p_j^\ast E_\rho$; when $e=(ii)$ is a self-loop with the orientation given by the ordered half-edges $h_+=(i,e,+1)< h_-=(i,e,-1)$, then the form $\omega_{(ii)}$ is valued in $p^\ast_i (E_{\rho, h_+}\otimes E_{\rho, h_-})$.
In particular, since $T^\ast\omega=-\omega$, we conclude for $i\neq j$,
\begin{equation}
    \omega_{(ij)}=-\omega_{(ji)}.
\label{eq:6.2.2sss}
\end{equation}

Taking a decorated trivalent graph $\Gamma$ of order $n$, then $2n =2 \ord(\Gamma) = |v(\Gamma)|$ is the number of vertices of $\Gamma$. 
At each vertex, the weight system $\mathrm{Tr}$ at each vertex of $\Gamma$ can be viewed as a flat skew-symmetric section of $(E^\vee_\rho)^{\otimes 3}\rightarrow M$ (also cf. Subsection \ref{ss3.3Sept}).


The order of the product manifold $M^{2n}$ (which has an induced orientation from $o(M)$) coincides with the given numbering on $v(\Gamma)$, or equivalently, we may write $(x_1,\ldots,x_{2n})=(x_i)_{i\in v(\Gamma)}\in M^{2n}$. At each point $(x_1,\ldots,x_{2n})\in M^{2n}$, we have the tensor product of vector bundles:
\begin{equation}
    (E^\vee_\rho)^{\otimes 3}_{x_1}\otimes (E^\vee_\rho)^{\otimes 3}_{x_2}\otimes \cdots \otimes (E^\vee_\rho)^{\otimes 3}_{x_{2n}},
    \label{eq:product}
\end{equation}
then by considering the set of half-edges $h(\Gamma)$ of $\Gamma$, each factor $E^\vee_{\rho,x_j}$ in \eqref{eq:product} can be regard as a copy of $E^\vee_\rho$ indexed by a half-edge $h$ attached to vertex $j$. Then we consider the pull-back of $E^\vee_\rho\boxtimes E^\vee_\rho \rightarrow M\times M$ by $p_{ij}$ for a non-self-loop edge $e=(ij)$ and the pull-back $E^\vee_\rho\otimes E^\vee_\rho \rightarrow M$ by $p_i$ for a self-loop edge $e=(ii)$, then each copy of $E^\vee_\rho$ is clearly index by the half-edges of $e$, therefore we get again the tensor product of vector bundles as in \eqref{eq:product}. We always identify these two perspectives for the vector bundle $p_1^\ast (E^\vee_\rho)^{\otimes 3} \otimes p_2^\ast (E^\vee_\rho)^{\otimes 3} \otimes \cdots\otimes p_{2n}^\ast (E^\vee_\rho)^{\otimes 3}$ over $M^{2n}$ or $C_{2n}(M)$.

Set
\begin{equation}
	Z_{\Gamma}(\omega) = \int_{C_{2n}(M)} \Big(\bigotimes_{i \in v(\Gamma)} \Tr_i \Big) \Big[\bigwedge_{e \in e(\Gamma)} \omega_{e}\Big],
\label{eq:6.2.4paris}
\end{equation}
where $e=(ij)$ means the edge connecting the vertex $i$ to the vertex $j$ (which always carry an orientation when $i=j$). Note that in \eqref{eq:6.2.4paris}, the factor $\bigotimes_{i \in v(\Gamma)} \Tr_i $ corresponds to the tensor product of the decoration $\Tr$ at each vertex given in Definition \ref{def:6.3.1August}. Note that by our convention, to apply $\bigotimes_{i \in v(\Gamma)} \Tr_i $ on $\bigwedge_{e \in e(\Gamma)} \omega_{e}$, we need to pair the factor of $E^\vee_\rho$ in $\bigotimes_{i \in v(\Gamma)} \Tr_i$ corresponding to a half-edge $h\in h(\Gamma)$ with the factor $E_\rho$ in $\bigwedge_{e \in e(\Gamma)} \omega_{e}$ that corresponds to the same half-edge $h$.

\begin{proposition}[Definition of $Z_{-}(\omega)$] \label{prop:6.2.1paris}
Let $\Gamma$ be a decorated trivalent graph of order $n$. If $\Gamma'$ is another decorated trivalent graph of order $n$ which is equivalent to $\Gamma$ via the equivalence relation of \eqref{eq:6.2.4}, then $\sign(\Gamma, \Gamma') Z_{\Gamma'}(\omega)=Z_\Gamma(\omega)$. 

Therefore, the following linear map is well-defined: 
\begin{equation}
	Z_{-}(\omega): \mathcal{GC}^0_{\mathfrak{g}:n} \rightarrow \mathbb{R}, \quad \Gamma \mapsto Z_{\Gamma}(\omega).
 \label{eq:7.1.4sept}
\end{equation} 
\end{proposition}
\begin{proof}
Note that if we permute the numbering on $v(\Gamma)$ by a permutation of order $p$, this gives an auto-identification of $C_{2n}(M)$ with the orientation change by $(-1)^p$. If we change an orientation of a non-self-loop edge $e=(ij)$, it is equivalent to change $\omega_{(ij)}$ to $\omega_{(ji)}$, we obtain a factor $(-1)$ by \eqref{eq:6.2.2sss}, let $(-1)^m$ denoete the total change by this kind of operation. If we change the orientation of a self-loop edge at vertex $i$, we obtain the same term by our convention and the property $T^\ast\xi=-\xi$. If we permute the numbering on $e(\Gamma)$ which implies $k$-times change cyclic orders at trivalent vertices, we obtain a factor $(-1)^k$ by corresponding sign change on associated cubic traces. This way, we obtain the sign $(-1)^{p+m+k}$ when we compare $Z_{\Gamma'}(\omega)$ with $Z_\Gamma(\omega)$, it completes the proof of our proposition.
\end{proof}

\begin{lemma}\label{lm:7.1.2sk}
    Fix a propagator $\omega$. If $\Gamma_1$, $\Gamma_2$ are two decorated trivalent graphs, then
\begin{equation}
    Z_{\Gamma_1\cup \Gamma_2}(\omega)=Z_{\Gamma_1}(\omega)Z_{\Gamma_2}(\omega).
\end{equation}
\end{lemma}
\begin{proof}
Set $n_1=\mathrm{ord}(\Gamma_1)$, $n_2=\mathrm{ord}(\Gamma_2)$, then $\mathrm{ord}(\Gamma_1\cup\Gamma_2)=n_1+n_2$. Consider the smooth map
\begin{equation}
    \Psi: \mathrm{Conf}_{2n_1+2n_2}(M)\rightarrow  \mathrm{Conf}_{2n_1}(M)\times  \mathrm{Conf}_{2n_2}(M).
\end{equation}
It induces a diffeomorphism between $\mathrm{Conf}_{2n_1+2n_2}(M)$ and $\mathrm{Image}(\Psi)$, and $\mathrm{Image}(\Psi)$ has full measure in $\mathrm{Conf}_{2n_1}(M)\times  \mathrm{Conf}_{2n_2}(M)$, that is, $\big(\mathrm{Conf}_{2n_1}(M)\times  \mathrm{Conf}_{2n_2}(M)\big)\setminus \mathrm{Image}(\Psi)$ has the Lebesgue measure zero. Moreover, the tangent map of $\Psi$ acts as identity map on each copy of $TM$.

In the same time, on $\mathrm{Conf}_{2n_1+2n_2}(M)\simeq \mathrm{Image}(\Psi)$, we have the identity
\begin{equation}
    \Big(\bigotimes_{i \in v(\Gamma_1\cup\Gamma_2)} \Tr_i \Big)\Big[\bigwedge_{e \in e(\Gamma_1\cup\Gamma_2)} \omega_{e} \Big]=\Big(\bigotimes_{i \in v(\Gamma_1)} \Tr_i \Big)\Big[\bigwedge_{e \in e(\Gamma_1)} \omega_{e}\Big] \wedge \Big(\bigotimes_{i \in v(\Gamma_2)} \Tr_i \Big)\Big[\bigwedge_{e \in e(\Gamma_2)} \omega_{e}\Big].
\label{eq:7.1.7sksk}
\end{equation}
In the definition \eqref{eq:6.2.4paris}, we can replace the integrals on $C_{2n}(M)$ by the integrals on $\mathrm{Conf}_{2n}(M)$ or on an open dense subset with full measure. Therefore, our lemma follows from the relation \eqref{eq:7.1.7sksk}.
\end{proof}

    In the above definition, $Z_\Gamma(\omega)$ depends on the decorations of connected graphs. Following Lemma \ref{lm:6.6.2}, for a fixed propagator, we can get configuration integrals depending only on the underlying topological graph as follows. 
    
\begin{lemma}
    Fix a propagator $\omega$. Let $\Gamma$ be a decorated trivalent graph, and let $\mathrm{or}_{\Gamma}$ be the relative orientation as in Definition \ref{lm:6.6.2kk}. Then, the quantity
    \begin{equation}\label{eq:super_prop}
       \mathrm{or}_{\Gamma} \cdot Z_{\Gamma}(\omega)
    \end{equation}
    is independent of the choice of numbering of $v(\Gamma)$,  $e(\Gamma)$ and orientations of edges. In other words, it only depends on the underlying topological graph of $\Gamma$.
\end{lemma}

\begin{remark}
    Definition in \eqref{eq:super_prop} is essentially the same as \cite{AS2} where they use super propagator to define their integral invariants. The key difference is that we are allowed to permute freely each factor $E^\vee_\rho$ in $p_1^\ast (E^\vee_\rho)^{\otimes 3} \otimes p_2^\ast (E^\vee_\rho)^{\otimes 3} \otimes \cdots\otimes p_{2n}^\ast (E^\vee_\rho)^{\otimes 3}$ to make $\left(\bigotimes_{i \in v(\Gamma)} \Tr_i \right)$ pair with $\bigwedge_{e \in e(\Gamma_2)} \omega_{e}$, while permutations produce nontrivial signs in the formalism of super propagators of \cite{AS2}. Our definition here is inspired by that of \cite{Les, MR4521898} for the integral invariants of rational homology $3$-spheres.  
\end{remark}

We first give the following theorem which can be viewed as a direct higher-order extension of the $Z_1$-invariant of Cattaneo--Shimizu \cite{CS}, where we use the general propagators to define the integral invariants for the cocycles in $\mathcal{GC}^0_{\mathfrak{g}:n}$. 


\begin{theorem}\label{thm:6.2.5ss}
Fix a homotopy class $[f]$ of framing of $M$ and an orientation $o(M)$. Let $E_{\rho}$ be an acyclic local system on $M$ associated with a representation $\rho : \pi_1(M) \rightarrow G \overset{\Ad}{\longrightarrow} \Aut(\mathfrak{g})$. Let $\Gamma \in \mathcal{GC}^0_{\mathfrak{g}:n}$ be a cocycle (that is, $\delta \Gamma=0$). Then the number $Z_{\Gamma}(\omega)\in\bR$ is independent of the choice of the propagator $\omega$ or the framing $f\in [f]$, which is called the integral invariant associated with the cocycle $\Gamma$.

Therefore, the linear functional 
$$Z(M,\rho,[f]): \ker(\delta|_{\mathcal{GC}^0_{\mathfrak{g}:n}})=H^0(\mathcal{GC}^\bullet_{\mathfrak{g}:n},\delta) \rightarrow \bR,$$ given by $Z(M,\rho,[f])(\Gamma):=Z_\Gamma(\omega)$ with any propagator $\omega$ constructed from a given framing $f\in [f]$, is an invariant of $(M,o(M),[f])$ and of the acyclic local system $E_\rho$.
\end{theorem}

The proofs of  Theorem \ref{thm:6.2.5ss}  will be given in Subsection \ref{ss6.4paris}.

Next, we connect the integral invariants associated with graph complex $\mathcal{GC}_{\mathfrak{g}}$ possibly with self-loops and those associated with $\mathcal{G}_{\mathfrak{g}}$ without self-loops. This extends the idea of Theorem \ref{prop:3.2}, where the introduction of an adapted propagator is the key step. Recall that $\mathcal{GC}^{\prime, 0}_{\mathfrak{g}}\subset \mathcal{GC}^{0}_{\mathfrak{g}}$ is the subspace consisting of all decorated trivalent graphs always with self-loops.

\begin{proposition} \label{prop:6.3.1}
Let $\omega^\sharp$ be an adapted propagator with $	\mathfrak{i}_\partial^{*}(\omega^\sharp) = I(\eta) + q_\partial^{\ast}(\xi^\sharp)$ as in Definition \ref{def:6.1.1}, the map $Z_{-}(\omega^\sharp)$ restricts to zero on $\mathcal{GC}^{\prime, 0}_{\mathfrak{g}}$.
\end{proposition}

\begin{proof}
The vanishing argument is almost the same as the case of the dumbbell graph. We focus only on integrand $q^{\ast} p_i^{\ast} \xi^\sharp$ associated with a self-loop $(ii)$. Note that by our definition of an adapted propagator $\omega^\sharp$, we have $\mathfrak{L}(\xi^\sharp)=0$.

Suppose that the vertex $i$ is connected by an edge $(ij)$ with $j\neq i$. By \eqref{eq:1.5.2bis} and \eqref{eq:3.0.1}, the integrand associated with the vertex $i$ and the edges $(ii)$ and $(ij)$ becomes
\begin{equation}
\begin{split}
	 & \Tr_i [\omega^\sharp_{(ij)} q^{\ast} p_i^{\ast} \xi^\sharp] \\
	 = &B_{i}(\mathfrak{L}(\xi^\sharp), \omega^\sharp_{(ij)})\\
	 =& 0
\end{split}
\end{equation}
where we suppress other forms associated with edges connecting the vertex $j$ and the associated cubic trace for simplicity. Therefore, it means that, if a decorated trivalent graph $\Gamma$ has at least one self-loop, $Z_{\Gamma}(\omega^\sharp) = 0$. This completes the proof.
\end{proof}

For each order $n$, we have the quotient graph complex $(\mathcal{G}^\bullet_{\mathfrak{g}: n},\delta^\sharp)$ defined in \eqref{eq:6.1.13paris} for the decorated graphs without self-loops.

\begin{theorem}\label{thm:6.3.2}
Assume $E_\rho$ to be acyclic. For any adapted propagator $\omega^\sharp$ in Definition \ref{def:6.1.1}, the map $Z_{-}(\omega^\sharp)$ factors through the quotient $\mathcal{GC}^0_{\mathfrak{g}: n}/\mathcal{GC}^{\prime, 0}_{\mathfrak{g}: n}\simeq \mathcal{G}^0_{\mathfrak{g}: n}$:
	
	\begin{equation}
		\begin{tikzcd}
 \mathcal{GC}^0_{\mathfrak{g}: n} \arrow[r,"Z_{-}(\omega^\sharp)"]\arrow[d] & \mathbb{R}\\
\mathcal{GC}^0_{\mathfrak{g}: n}/\mathcal{GC}^{\prime, 0}_{\mathfrak{g}: n}\simeq \mathcal{G}^0_{\mathfrak{g}: n} \arrow[ru] & 
\end{tikzcd}
	\end{equation}
	This way, we have a linear map
	\begin{equation}
		Z_{-}(\omega^\sharp) : \mathcal{G}^0_{\mathfrak{g}: n} \rightarrow \mathbb{R}.
	\end{equation}
\end{theorem}

\begin{proof}
	This follows from \eqref{eq:6.1.13paris} and Propositions \ref{prop:6.2.1paris} \& \ref{prop:6.3.1}.
\end{proof}

Our main results for this section are as follows. The first part refines \cite[Theorem 1.1]{BC2} for the acyclic local systems defined via the adjoint representation of $\pi_1(M)$ on a semi-simple Lie algebra. In the second part, we connect the integral invariants defined by graph cocycles with self-loops to those without self-loops. 



 \begin{theorem}\label{thm:6.3.1}
Fix a homotopy class $[f]$ of framing of $M$ and an orientation $o(M)$. Let $E_{\rho}$ be an acyclic local system over $M$ associated with a representation $\rho : \pi_1(M) \rightarrow G \overset{\Ad}{\longrightarrow} \Aut(\mathfrak{g})$.  

\begin{enumerate}[(1)]
    \item Let $\Gamma \in \mathcal{G}^0_{\mathfrak{g}: n}$ be a cocycle of order $n$ (that is, $\delta^\sharp \Gamma=0$). Then the number $Z_{\Gamma}(\omega^\sharp)\in\bR$ is independent of the choice of the adapted propagator $\omega^\sharp$ or the framing $f\in [f]$, which is called the integral invariant associated to the cocycle $\Gamma$.

Therefore, the linear functional 
$$Z^\sharp(M,\rho,[f]): \ker(\delta^\sharp|_{\mathcal{G}^0_{\mathfrak{g}: n}})=H^0(\mathcal{G}^\bullet_{\mathfrak{g}: n}, \delta^\sharp)\rightarrow \bR,$$ given by $Z^\sharp(M,\rho,[f])(\Gamma)=Z_\Gamma(\omega^\sharp)$ with any adapted propagator $\omega^\sharp$ defined with a framing $f\in [f]$, is an invariant of $(M,o(M),[f])$ and local system $E_\rho$.

\item  With the notation in (1), we have the following commutative diagram:

\begin{equation}\label{eq:7.2.1}
\begin{tikzcd}
H^0(\mathcal{GC}^\bullet_{\mathfrak{g}: n}, \delta)\arrow[r, "{Z(M,\rho, [f])}"]\arrow[d] &[1.5cm] \bR \arrow[d,"="]\\
H^0(\mathcal{G}^\bullet_{\mathfrak{g}: n}, \delta^\sharp)\arrow[r, "{Z^\sharp(M,\rho, [f])}"] & \bR
\end{tikzcd}
\end{equation}
where the left vertical map is induced by the quotient map $\mathcal{GC}^0_{\mathfrak{g}: n}\rightarrow \mathcal{GC}^0_{\mathfrak{g}: n}/\mathcal{GC}^{\prime, 0}_{\mathfrak{g}: n}\simeq \mathcal{G}^0_{\mathfrak{g}: n}$, which is already explained in \eqref{eq:6.5.3sept}.

    More precisely, let $\Gamma\in \mathcal{GC}^0_{\mathfrak{g}: n}$ be a cocycle, and let $\omega$ be any propagator as in Definition \ref{defn:prop1}. Let $\Gamma'\in \mathcal{G}^0_{\mathfrak{g}: n}$ be the cocycle given by removing the terms with self-loops from $\Gamma$. Then we have, for any adapted propagator $\omega^\sharp$,
    \begin{equation}
        Z_\Gamma(\omega)=Z(M,\rho,[f])(\Gamma)=Z^\sharp(M,\rho,[f])(\Gamma')=Z_{\Gamma'}(\omega^\sharp).
    \end{equation}
    \end{enumerate}
\end{theorem}

The proof of Theorem \ref{thm:6.3.1}  will be given in Subsection \ref{ss7.3s}.


Finally, combining Theorems \ref{thm:6.2.5ss} \& \ref{thm:6.3.1} with Proposition \ref{prop:CSseries}, we obtain a generating series of perturbative invariants of a closed 3-manifold associated with acyclic representation $\rho : \pi_1(M) \rightarrow G \overset{\Ad}{\longrightarrow} \Aut(\mathfrak{g})$.

\begin{corollary}\label{cor:generating_func}
Fix a homotopy class $[f]$ of smooth framing of $M$ and an orientation $o(M)$. Let $E_{\rho}$ be an acyclic local system over $M$ associated with a representation $\rho : \pi_1(M) \rightarrow G \overset{\Ad}{\longrightarrow} \Aut(\mathfrak{g})$. 
\begin{enumerate}[(1)]
\item Let $\omega$ be a propagator. Consider the formal sum
\begin{equation}
   \log \mathcal{Z}_{\mathrm{CS}}(M, \rho, [f])= \sum_{\text{connected } \mathfrak{G}} \frac{\hbar^{\mathrm{ord}(\mathfrak{G})}}{|\Aut(\mathfrak{G})|}  Z_{\Psi_n(\mathfrak{G})}(\omega) \in \bR[[\hbar]],
\label{eq:7.6.6new}
\end{equation}
where the sum runs over all the connected topological trivalent graph $\mathfrak{G}$.
  Then, it is independent of the choice of propagator $\omega$. Therefore, $\mathcal{Z}_{\mathrm{CS}}(M, \rho, [f])$ is an invariant of $(M, o(M), \rho, [f])$.
\item Let $\omega^\sharp$ be any adapted propagator, then the formal sum in \eqref{eq:7.6.6new} satisfies the following identity
\begin{equation}
   \log \mathcal{Z}_{\mathrm{CS}}(M, \rho, [f])= \sum_{\substack{\text{connected } \mathfrak{G} \\ \text{without self-loops}}} \frac{\hbar^{\mathrm{ord}(\mathfrak{G})}}{|\Aut(\mathfrak{G})|}  Z_{\Psi_n(\mathfrak{G})}(\omega^\sharp) \in \bR[[\hbar]],
\label{eq:7.1.15new}
\end{equation}
where the sum runs over all the connected topological trivalent graph $\mathfrak{G}$ without self-loops.
\end{enumerate}
\end{corollary}

\begin{remark}\label{rk:7.2.10ss}
Corollary \ref{cor:generating_func} (1) is an analogous result given in \cite{AS} and \cite{AS2} that they used propagators constructed from the de Rham--Hodge Laplacian.
\end{remark}

\subsection{A variation formula and proof of Theorem \ref{thm:6.2.5ss}}\label{ss6.4paris}
 We will use some ideas from the proofs of \cite[Theorem 4.7]{BC} and \cite[Theorem 1.1]{BC2} to achieve our proofs of Theorem \ref{thm:6.2.5ss} and Theorem \ref{thm:6.3.1}. We will provide the necessary details for completeness. One of the differences from theirs is that we compute the graphs obtained by contracting non-regular edges in detail, which involves self-loops. For the case of $2$-loop trivalent graphs, such a computation was already explained in \cite[\S 4.2]{CS}.

 The proof of Theorem \ref{thm:6.2.5ss} is given in a way similar to that in \cite[Theorem 4.7]{BC}. Here is the outline of the proof: we consider a smooth one-parameter family of propagators over the unit interval $I=[0,1]$ as a parameter space. Then, this family of propagators gives rise to a family of integrals associated with a given decorated trivalent graph. To prove the independence of the choice of propagators, it suffices to show that this family of integrals is constant on $I$, or equivalently, its differential on $I$ is identically $0$. Stokes' formula and Kontsevich's lemma (Lemma \ref{lem:Kontsevich}) tell that there are non-vanishing boundary contributions, but they can be made zero by graph cocycle relation. In this way, we finally obtain Theorem \ref{thm:6.2.5ss}. In fact, we at first will prove a result analogous to \cite[Corollary 4.12]{BC} from which Theorem \ref{thm:6.2.5ss} follows clearly. This result gives us a formula for the variations of $Z_\Gamma(\omega)$ as $\omega$ varies smoothly and for any decorated trivalent graph $\Gamma$ which is not necessary to be a cocycle.

 Note that in our construction the map $Z(M, \rho, [f])$ exactly gives rise to an invariant of framed 3-manifold with acyclic representation $\rho$ associated to a graph cocycle $\Gamma$. This is different from \cite{BC2} where, to obtain a graph cocycle invariant, a modification of $Z(M, \rho, [f])$ is required to cancel a boundary contribution by adding correction terms.

\subsubsection{A variation formula for a family of propagators}

Let $I=[0,1]$ denote the unit interval with the standard coordinate $\tau\in [0,1]$. The vector bundles on $M$, $C_2(M)$, etc, are viewed naturally as vector bundles on $I\times M$, $I\times C_2(M)$, etc, respectively, and so do the differential forms. We also extend the action of $T$ on $I\times \cdots$ by trivial action on the factor $I$. Let $d^{\mathrm{tot}}=d\tau\wedge\frac{\partial}{\partial \tau}+d^M$ denote the total differential on the product space $I\times M$. We will use the same notation for the spaces $I\times C_2(M)$, $I\times \partial C_2(M)$, $I\times M\times M$, etc.

If $\Gamma$ is a connected decorated graph with degree $1$ and without any external edges, due to our convention that the minimal valency at each vertex is at least $3$, we can conclude that $\Gamma$ has exactly one vertex with valency $4$ and all other vertices are trivalent. Set $m=|v(\Gamma)|$. Then $m$ has to be an odd integer, and we have
\begin{equation}
\mathrm{ord}(\Gamma)=\frac{1}{2}(m+1).
\end{equation}

We consider a pair of differential $2$-forms $(\widetilde{\omega},\widetilde{\xi})\in \Omega^2 (I\times C_2(M); F_\rho)\times \Omega^2 (I\times M; E_\rho\otimes E_\rho)$ such that 
\begin{equation}
    T^\ast\widetilde{\omega}=-\widetilde{\omega},\; T^\ast\widetilde{\xi}=-\widetilde{\xi}.
    \label{eq:7.2.2}
\end{equation}
As in \eqref{eq:6.2.2june}, we associate a $2$-form on $I \times C_m(M)$ to each $e=(ij)\in e(\Gamma)$ as follows
\begin{equation}\label{eq:7.2.2sept}
	\widetilde{\omega}_{e} := \begin{cases}
 p_{ij}^{\ast} \widetilde{\omega}	& \text{if $e=(ij)$ with $i \neq j$},\\
  q^{\ast} p_i^{\ast} \widetilde{\xi} & \text{if $e=(ii)$ is a directed self-loop},
 \end{cases}
\end{equation}
Let $\sigma: I \times C_{m}(M) \rightarrow I$ denote the obvious projection, and let $\sigma_{\ast}: \Omega^{\bullet+3m}(I \times C_{m}(M)) \rightarrow \Omega^\bullet(I)$ be the fiber integration (Definition \ref{def:2.4.1}).
Similar to \eqref{eq:6.2.4paris}, we define
\begin{equation}
	Z_{\Gamma}(\widetilde{\omega},\widetilde{\xi}) = \sigma_\ast\Big\{\Big(\bigotimes_{i \in v(\Gamma)} W_i \Big)\Big[\bigwedge_{e \in e(\Gamma)} \widetilde{\omega}_{e}\Big]\Big\}\in \Omega^1(I),
\label{eq:7.2.4paris}
\end{equation}
where $W_i$ denotes the weight system at vertex $i$, it is $\mathrm{Tr}_i$ when vertex $i$ is trivalent and is $\pm W_I$, $\pm W_H$, $\pm W_X$ at the only vertex of valency $4$.

If $\Gamma$ is a decorated trivalent graph, then $Z_{\Gamma}(\widetilde{\omega},\widetilde{\xi})\in \Omega^0(I)$ can also be defined by considering a smooth family of the integrations as in \eqref{eq:6.2.4paris}. In summary, we have the following result.
\begin{lemma}
    The following linear map is well-defined for $j=0,1$,
\begin{equation}
	Z_{-}(\widetilde{\omega},\widetilde{\xi}): \mathcal{GC}^j_{\mathfrak{g}} \rightarrow \Omega^j(I), \quad \Gamma \mapsto Z_{\Gamma}(\widetilde{\omega},\widetilde{\xi}).
 \label{eq:7.2.5sept}
\end{equation} 
\end{lemma}
 \begin{proof}
 If $j=0$, this is exactly a family version of Proposition \ref{prop:6.2.1paris}, which follows from the same proof since the boundary condition \eqref{eq:2.57} for $\{\widetilde{\omega}|_{\{\tau\}\times C_2(M)}\}_{\tau\in I}$ is not needed.

 If $j=1$, by \eqref{eq:7.2.2} and \eqref{eq:7.2.2sept}, the same arguments in the proof of Proposition \ref{prop:6.2.1paris} shows that the definition \eqref{eq:7.2.4paris} is compatible with the sign convention on the decorated graphs of degree $1$. For internal IHX relation, we can consider a decorated graph to be the sum of three decorated graphs $\Gamma_j$, $j=1,2,3$, which have exactly the same underlying topological graph and the decorations on the edges and vertices except for the different weight systems (I, H, X, respectively) at the vertices of valency $4$, then by \eqref{eq:Jacobi-indentity}, the sum of $\left(\bigotimes_{i \in v(\Gamma_j)} W_i \right) \left[\bigwedge_{e \in e(\Gamma_j)} \widetilde{\omega}_{e}\right]$ vanishes identically. Then the lemma holds.
 \end{proof}

Now we can state our result for the variations of the integrals $Z_\Gamma(\omega)$ defined in \eqref{eq:6.2.4paris}, when $\omega$ varies smoothly, in the spirit of the second part of \cite[Corollary 4.12]{BC}.
\begin{proposition}\label{prop:7.2.2s}
    Let the pair $(\widetilde{\omega},\widetilde{\xi})\in \Omega^2 (I\times C_2(M); F_\rho)\times \Omega^2 (I\times M; E_\rho\otimes E_\rho)$ be such that
    \begin{itemize}
        \item $d^\mathrm{tot}\widetilde{\omega}=0,\; d^\mathrm{tot}\widetilde{\xi}=0$;
        \item  $T^\ast\widetilde{\omega}=-\widetilde{\omega},\; T^\ast\widetilde{\xi}=-\widetilde{\xi};$
        \item there exists a closed smooth $2$-form $\widetilde{\mu}\in \Omega^2(I\times \partial C_2(M);\mathbb{R})$ such that
        \begin{itemize}
            \item $\widetilde{\mu}$ is a vertical $2$-form with respect to the submersion 
            $$I\times \partial C_2(M)\simeq I\times M\times \mathbb{S}^2\rightarrow M;$$
            \item $\widetilde{q_\partial}_\ast\widetilde{\mu}=1$ on $I\times M$, where $\widetilde{q_\partial}:=(\mathrm{Id}_I, q_\partial): I\times \partial C_2(M)\simeq I\times M\times \mathbb{S}^2\rightarrow I\times M$;
\item let $\widetilde{\mathfrak{i}_\partial}$ be the inclusion $I\times \partial C_2(M)\rightarrow I\times C_2(M)$, analogous to \eqref{eq:2.57}, we have
    \begin{equation}
\widetilde{\mathfrak{i}_\partial}^\ast(\widetilde{\omega})=\widetilde{\mu}\otimes \mathbf{1}+\widetilde{q_\partial}^\ast (\widetilde{\xi}),
\label{eq:7.2.6s}
    \end{equation}
where $\mathbf{1}$ is the flat section in Lemma \ref{lm:1.1}.
 \end{itemize}
    \end{itemize}
Then we have the following identity for any $\Gamma\in \mathcal{GC}^0_{\mathfrak{g}}$,
\begin{equation}
    d Z_\Gamma(\widetilde{\omega},\widetilde{\xi})=Z_{\delta \Gamma}(\widetilde{\omega},\widetilde{\xi})\in\Omega^1(I).
    \label{eq:7.2.7ssss}
\end{equation}
\end{proposition}

Let the pair $(\widetilde{\omega},\widetilde{\xi})$ be as given in the above proposition.
Let $\Gamma$ be a trivalent connected decorated graph with order $n$, then $\delta \Gamma$ is a linear combination of connected decorated graphs in $\mathcal{GC}^1_{\mathfrak{g}: n}$. For simplicity, set
\begin{equation}
    \Tr^{\Gamma}(\widetilde{\omega},\widetilde{\xi}) =  \Big(\bigotimes_{i \in v(\Gamma)} \Tr_i \Big)\Big[\bigwedge_{e \in e(\Gamma)} \widetilde{\omega}_{e}\Big] \in \Omega^{6n}(I \times C_{2n}(M)).
\label{eq:7.2.8kk}
\end{equation}
It is clear that $\Tr^{\Gamma}(\widetilde{\omega},\widetilde{\xi})$ is $d^{\mathrm{tot}}$-closed. By Stokes' formula (Proposition \ref{prop:2.4.1}), we get
\begin{equation}
\begin{split}
d Z_{\Gamma}(\widetilde{\omega},\widetilde{\xi}) &=d\,\sigma_\ast \Tr^{\Gamma}(\widetilde{\omega},\widetilde{\xi})\\
& = \sigma_{\ast}\big(	d^\mathrm{tot} \Tr^{\Gamma}(\widetilde{\omega},\widetilde{\xi}) \big) + \sigma_{\ast}^{\partial} \big(\widetilde{\mathfrak{i}_\partial}^\ast\Tr^{\Gamma}(\widetilde{\omega},\widetilde{\xi})\big)\\
&=  \sigma_{\ast}^{\partial} \big(\widetilde{\mathfrak{i}_\partial}^\ast\Tr^{\Gamma}(\widetilde{\omega},\widetilde{\xi})\big),
\end{split}
\label{eq:7.2.9kkk}
\end{equation}
where $\widetilde{\mathfrak{i}_\partial}$ denotes the inclusion $I\times \partial C_{2n}(M)\rightarrow I\times C_{2n}(M)$, and $\sigma^\partial:I\times \partial C_{2n}(M)\rightarrow I$ denotes the obvious projection. 

Therefore, our calculation reduces to that of $\sigma_{\ast}^{\partial} \big(\widetilde{\mathfrak{i}_\partial}^\ast\Tr^{\Gamma}(\widetilde{\omega},\widetilde{\xi})\big)$, where we need to investigate the geometry of the codimension-$1$ boundary $\partial^\ast C_{2n}(M)$, that is, $S^1(C_{2n}(M))$. Let $S$ be a subset of $\{1,2,\ldots, 2n\}$ or $v(\Gamma)$ with $\ell=|S|\geq 2$. Let $\partial_S C_{2n}(M)$ denote the component of $\partial^\ast C_2(M)$ corresponding to $M(\{S\})$ in the notation of  \ref{app:FMAS}, they are defined by collapsing points $\{\mathbf{x}_i\}_{i\in S}\in M^\ell$ into the same point. 

Note that a point in the open strata $M(\{S\})^0$ can be represented by
\begin{equation}
\begin{split}
    \Big(\mathbf{x}_S=(z,\ldots,z)\in\Delta_S\simeq M \, ;\, u_S\in \bR^*_+\backslash \big((T_z M)^S/T_zM - \{0\}\big)\, ;\hspace{2cm}&
    \\ \, \{\mathbf{x}_j\}_{j\not\in S}\in \mathrm{Conf}_{2n-\ell}(M) \text{ with } \mathbf{x}_j\neq z\Big),&
\end{split}
\end{equation}
where the normal vector $u_S$ does not have any two components which are equal. Let $\mathrm{Conf}_S(T_zM)$ denote the configuration space of vectors in $T_zM$ indexed by $S$, let $T_zM$ act on $\mathrm{Conf}_S(T_zM)$ by on-diagonal translations and let $\bR^\ast_+$ act on $\mathrm{Conf}_S(T_zM)$ by rescalings. Then we can rewrite the above requirements on $u_S$ as $u_S\in \mathrm{Conf}_S(T_zM)/T_zM\rtimes \bR^*_+$.
Consider the smooth projection $\mathrm{Pr}_S: M(\{S\})^0 \rightarrow \mathrm{Conf}_{2n-\ell+1}(M)$ which sends the above point in $M(\{S\})^0$ to the point $(z,\mathbf{x}_j, j\not\in S)\in \mathrm{Conf}_{2n-\ell+1}(M)$, then the fibre of this projection is given by $\mathrm{Conf}_S(T_zM)/T_zM\rtimes \bR^*_+\simeq \mathrm{Conf}_\ell(\bR^3)/\bR^3\rtimes \bR^*_+$.

Then we extend it smoothly to the projection, denoted by the same notation,
$$\mathrm{Pr}_S: \partial_S C_{2n}(M) \rightarrow C_{2n-\ell+1}(M).$$
The generic fibre of $\mathrm{Pr}_S$ is given by $F_S\simeq C_\ell(\mathbb{R}^3)/\mathbb{R}^3 \rtimes \mathbb{R}^\ast_{+}$. In particular, $\dim_\bR F_S=3\ell -4$.

To compute the contribution of $\partial_S C_{2n}(M)$ in $\sigma_{\ast}^{\partial} \big(\widetilde{\mathfrak{i}_\partial}^\ast\Tr^{\Gamma}(\widetilde{\omega},\widetilde{\xi})\big)$, we need the following Kontsevich's vanishing lemma.

\begin{lemma}[{\cite[Lemma 2.1]{Ko} \cite[Lemma 4.9]{BC}}]\label{lem:Kontsevich}

	Let $F_{S}$ denote the fiber of the face $\partial_{S} C_{2n}(M)$ corresponding to the collapse of $\ell$ points with coordinate $\mathbf{x}_j, j\in S$, that is $F_{S}$ is the generic fibre of $\mathrm{Pr}_S: \partial_S C_{2n}(M) \rightarrow C_{2n-\ell+1}(M)$. Fix a smooth framing $f: TM\rightarrow M\times \bR^3$, then it induces an identification $F_S=f^\ast( C_\ell(\mathbb{R}^3)/\mathbb{R}^3 \rtimes \mathbb{R}^\ast_{+})$.  Let $\eta \in \Omega^2(\mathbb{S}^2; \bR)$ be any volume form of $\mathbb{S}^2$ with $T^{\ast} \eta = - \eta$. For $i,j\in S$, $i\neq j$, let $\pi_{ij}: F_{S} \rightarrow \mathbb{S}^2$ be the projection defined as
 \begin{equation}
     \pi_{ij}: F_{S} \rightarrow \mathbb{S}^2; \quad (\mathbf{x}_j)_{j\in S} \mapsto \frac{\mathbf{x}_j - \mathbf{x}_i}{|\mathbf{x}_j - \mathbf{x}_i|}\quad (i \neq j)
 \end{equation}
 and $\pi_{ij}^{\ast} \eta$ be the pullback $\eta$ via $\pi_{ij}$. Then, any triple of indices $i,j,k$ in $S$ with $i\neq j$ and $i\neq k$, the integral vanishes:
 \begin{equation}
     \int_{\mathbf{x}_i} \pi_{ij}^{\ast}\eta \wedge  \pi_{ik}^{\ast}\eta = 0,
 \end{equation}
 where $\int_{\mathbf{x}_i}\cdots$ means the fibre-wise integration with respect to the projection of forgetting $\mathbf{x}_i$-coordinate: $C_\ell(\mathbb{R}^3)/\mathbb{R}^3 \rtimes \mathbb{R}^\ast_{+}\rightarrow C_{\ell-1}(\mathbb{R}^3)/\mathbb{R}^3 \rtimes \mathbb{R}^\ast_{+}$ (provided $\ell\geq 3$).
\end{lemma}

\begin{remark}
    With the same notation as above, note that Lemma \ref{lem:Kontsevich} immediately implies the following (original) statement. For any two sequences $s_i, t_i$ $(i=1,\ldots, L)$ of integers with $s_i \neq t_i$ $(1 \leq s_i,t_i \leq \ell)$, the integral vanishes:
    \begin{equation}
         \int_{(\mathbf{x}_1,\ldots, \mathbf{x}_\ell) \in F_{S}} \bigwedge_{i=1}^L \pi_{s_i t_i}^{\ast} \eta  =0.
    \end{equation}
\end{remark}

\begin{proof}[Proof of Proposition \ref{prop:7.2.2s}]
We use the above notation and consider the face $\partial_{S} C_{2n}(M)$. Note that sthe coordinate $(\mathbf{x}_1,\ldots, \mathbf{x}_{2n})$ corresponds to the vertices of the decorated trivalent graph $\Gamma$. We will regard the vertices in $S$ as the collapsing vertices.

Using the projection $\mathrm{Pr}_S: \partial_S C_{2n}(M) \rightarrow C_{2n-\ell+1}(M)$, and for a differential form $\alpha$ on $\partial_S C_{2n}(M)$, we can decompose $\alpha$ into two parts: the vertical direction and the basic direction. For the integral $\int_{\partial_S C_{2n}(M)}\alpha$ being nonzero, the degree of $\alpha$ shall be $6n-1$ with the vertical degree of $\alpha$ being $3\ell-4$. Moreover, we have
\begin{equation}
    \int_{\partial_S C_{2n}(M)}\alpha=\int_{C_{2n-\ell+1}(M)} (\mathrm{Pr}_S)_\ast\alpha.
    \label{eq:7.2.14sss}
\end{equation}

Let's consider the differential form $\Tr^{\Gamma}(\widetilde{\omega},\widetilde{\xi})|_{\partial_S C_{2n}(M)}$, which can be obtained by taking the product of $\widetilde{\omega}_e|_{\partial_S C_{2n}(M)}$. Note that by \eqref{eq:7.2.2sept}, if $e=(ii)$ is a self-loop edge, then $\widetilde{\omega}_e|_{\partial_S C_{2n}(M)}$ is always basic differential form (with respect to the projection $\mathrm{Pr}_S$); if $e=(ij)$ is non-self-loop edge such that $i$ or $j$ does not lie in $S$, then $\widetilde{\omega}_e|_{\partial_S C_{2n}(M)}$ is also basic. To have the vertical directions in $\widetilde{\omega}_e|_{\partial_S C_{2n}(M)}$, we need $e=(ij)$ with $i,j\in S, i\neq j$, and in this case the vertical form contributed by $\widetilde{\omega}_e|_{\partial_S C_{2n}(M)}$ is $\widetilde{\mu}$ in \eqref{eq:7.2.6s}. Our assumptions on $\widetilde{\mu}$ implies that its contribution in $\widetilde{\omega}_e|_{\partial_S C_{2n}(M)}$ can be written as follows, for $\tau\in I$,
\begin{equation}
    \widetilde{\mu} = \pi^*_{ij}\eta_\tau +d\tau\wedge \pi^*_{ij}\beta_\tau,
\label{eq:7.2.15sss}
\end{equation}
where $\eta_\tau$ is a volume form on $\mathbb{S}^2$ (depending smoothly on $\tau$), and $\beta_\tau$ is a $1$-form on $\mathbb{S}^2$.

Let $\bb{e}_v$ be the total number of edges connecting two distinct collapsing vertices (in $S$) and let $\bb{e}_h$ be the total number of self-loop edges incident to the collapsing vertices. Let $\bb{e}_0$ be the number of edges connecting a collapsing vertex in $S$ with a non-collapsing one. Since we consider trivalent graphs, we have the relation $2(\bb{e}_v+\bb{e}_h) + \bb{e}_0 = 3 \ell$. Then, the maximal degree of vertical form in $\Tr^{\Gamma}(\widetilde{\omega},\widetilde{\xi})|_{\partial_S C_{2n}(M)}$ is $2\bb{e}_v$. Considering $(\mathrm{Pr}_S)_\ast \Tr^{\Gamma}(\widetilde{\omega},\widetilde{\xi})|_{\partial_S C_{2n}(M)}$, it is nonzero only if $2\bb{e}_v - (3\ell - 4) = 4 - \bb{e}_0-2\bb{e}_h\geq 0$.

Let us first consider the case $\ell=|S| \geq 3$. By \eqref{eq:7.2.15sss}, the integrand form along the vertical direction of $\mathrm{Pr}_S$ is given by a product of $\pi_{ij}^{\ast}\eta_\tau$ and $d\tau\wedge \pi_{ij}^{\ast}\beta_\tau$. Since $\ell\geq 3$, $\dim_\bR F_S=3\ell- 4\geq 5$, so that we shall have at least two non-self-loop edges attached to the collapsing vertices to reach this vertical degree, then by Kontsevich's vanishing lemma (Lemma \ref{lem:Kontsevich}), we get $(\mathrm{Pr}_S)_\ast \Tr^{\Gamma}(\widetilde{\omega},\widetilde{\xi})|_{\partial_S C_{2n}(M)}=0$.

The remaining case is that $\ell=2$ and $S=\{i,j\}$ ($i\neq j$) with $e=(ij)$ or $(ji)$ is an edge of $\Gamma$. If $e_1$, $e_2$ are two different non-self-loop edges in $\Gamma$ connecting the same vertices $S=\{i,j\}$, then $\widetilde{\omega}_{e_1}\wedge \widetilde{\omega}_{e_2}|_{\partial_S C_{2n}(M)}$ has two nontrivial terms
$$\widetilde{u}\otimes \mathbf{1} \wedge \widetilde{\xi} + \widetilde{\xi}\wedge \widetilde{u}\otimes \mathbf{1}.$$
The first term, in the computation of $\sigma_{\ast}^{\partial} \big(\widetilde{\mathfrak{i}_\partial}^\ast\Tr^{\Gamma}(\widetilde{\omega},\widetilde{\xi})\big)$, corresponds to the contraction operation on $e_1$, i.e., the term $\Gamma/e_1$ in $\delta \Gamma$, and the second term corresponds to $\Gamma/e_2$. It is similar to the case where we have three different non-self-loop edges with the same ending vertices. Note that in the definition of the weight system $W_i$ at a vertex of valency $4$, it is the same as decorating the contracted edge by the Casimir element $\mathbf{1}$ then applying the cubic trace $\mathrm{Tr}$. This way, we conclude from \eqref{eq:7.2.9kkk} and the assumption $\widetilde{q_\partial}_\ast \widetilde{\mu}=1$ that
\begin{equation}
d Z_{\Gamma}(\widetilde{\omega},\widetilde{\xi})= \sum_{\substack{e\in e(\Gamma) \\ \text{non-self-loop edge}}} \pm Z_{\Gamma/e}(\widetilde{\omega},\widetilde{\xi}).
\label{eq:7.2.16kkk}
\end{equation}
The last step is to calculate precisely the sign $\pm$ in front of each term and then check the compatibility with the sign convention \eqref{eq:6.2.6} in the definition of $\delta \Gamma$ in \eqref{eq:6.2.5}.

By Proposition \ref{prop:general_graph_complex}, the map $\delta$ is well-defined under the sign relation \eqref{eq:6.2.4}, so that we can assume that $S=\{1,2\}$, $e=(12)$ is the edge numbered as $1$. Then \eqref{eq:7.2.16kkk} can be written as
\begin{equation}
d Z_{\Gamma}(\widetilde{\omega},\widetilde{\xi})= Z_{\Gamma/e}(\widetilde{\omega},\widetilde{\xi})+\sum_{\substack{\text{other }e'\in e(\Gamma) \\ \text{non-self-loop edge}}} \pm Z_{\Gamma/e'}(\widetilde{\omega},\widetilde{\xi}),
\end{equation}
while we have $\delta\Gamma= \Gamma/e +\sum_{\substack{\text{other }e'\in e(\Gamma) \\ \text{non-self-loop edge}}} \pm \Gamma/e'$. This way, we get exactly \eqref{eq:7.2.7ssss} for a connected decorated trivalent graph $\Gamma$. Then combing this result with \eqref{eq:6.4.15Sept} and Lemma \ref{lm:7.1.2sk}, we complete the proof for general $\Gamma$ in $\mathcal{GC}^0_{\mathfrak{g}}$.
\end{proof}

\subsubsection{Proof of Theorem \ref{thm:6.2.5ss}}\label{sss7.2.2}
Note that we always fix an orientation $o(M)$ of $M$.
Let $f$ and $f'$ be two smooth framings of $M$ which are homotopic, and let $\eta$, $\eta'$ be two normalized volume forms on $\mathbb{S}^2$. Let $(\omega, f,\eta,\xi)$, $(\omega',f',\eta',\xi')$ be two propagators defined for the acyclic local sytem $E_\rho$ as in Definition \ref{defn:prop1}. By Proposition \ref{prop:uniqueness}, the cohomology class $[\omega]=[\omega']$ is unique. But we need a more explicit relation between $\omega, \omega'$ with which we can apply Proposition \ref{prop:7.2.2s}.

Recall that $H^{2}(I \times \mathbb{S}^2; \mathbb{R}) \simeq H^2(\mathbb{S}^2; \mathbb{R})$.
For two $T$-asymmetric normalized volume form $\eta$, $\eta'$ on $\mathbb{S}^2$, there exists a closed $2$-form $\widetilde{\eta}\in\Omega^2(I\times \mathbb{S}^2 ;\bR)$ such that
\begin{equation}
\widetilde{\eta}_{\tau=0}=\eta,\;  \widetilde{\eta}_{\tau=1}=\eta'.
\end{equation}
The closedness of $\widetilde{\eta}$ implies that for each $\tau\in I$, 
\begin{equation}
\int_{\mathbb{S}^2}\widetilde{\eta}_\tau=1.
\end{equation}
We also require $T_{\mathbb{S}^2}^\ast \widetilde{\eta}=-\widetilde{\eta}$.

Since $f$ and $f'$ are homotopic, let $\widetilde{f}: I\times S(TM)\rightarrow I\times M\times\mathbb{S}^2$ denote the smooth path of framings which connects $f$ ($\tau=0$) and $f'$($\tau=1$). Set
\begin{equation}
I(\widetilde{\eta})=\widetilde{f}^\ast(\widetilde{\eta})\otimes \mathbf{1}\in\Omega^2_-(I\times \partial C_2(M); E_\rho\otimes E_\rho).
\end{equation}
At the same time, we have
\begin{equation}
H^2(I\times M\times M; E_\rho\boxtimes E_\rho)\simeq H^2(M\times M; E_\rho\boxtimes E_\rho)\oplus \mathbb{R}d\tau\wedge H^1(M\times M; E_\rho\boxtimes E_\rho)=0.
\end{equation}
Note that by Proposition \ref{cor:2.2.1}, the cohomology class $[\xi]$ is uniquely determined for any propagator $\omega$ associated with $E_\rho$.

Now we can follow the arguments as the proofs of Proposition \ref{cor:2.2.1} and Theorem \ref{thm:enhanced} (see also the proof of \cite[Proposition 2.1]{CS}) to construct a propagator on $I\times C_2(M)$ with analogous properties. More precisely, there exist a closed $2$-form $\widetilde{\omega}\in \Omega^2_-(I\times C_2(M); F_\rho)$ (that is, $d^{\mathrm{tot}}\widetilde{\omega}=0$) and closed $2$-form $\widetilde{\xi}\in \Omega^{2}_{-}(I\times \Delta; E_\rho\otimes E_\rho)$ such that
\begin{itemize}
    \item $\widetilde{\omega}_{\{0\}\times C_2(M)}=\omega$, $\widetilde{\omega}_{\{1\}\times C_2(M)}=\omega'$, or equivalently, $\widetilde{\xi}_{\{0\}\times C_2(M)}=\xi$, $\widetilde{\xi}_{\{1\}\times C_2(M)}=\xi'$;
    \item analogous to \eqref{eq:2.57}, we have
    \begin{equation}
\widetilde{\mathfrak{i}_\partial}^\ast(\widetilde{\omega})=I(\widetilde{\eta})+\widetilde{q_\partial}^\ast (\widetilde{\xi}).
    \end{equation}
\end{itemize}

Now we can take $\widetilde{\mu}$ to be $\widetilde{f}^\ast\widetilde{\eta}$, and then $\widetilde{\omega}$, $\widetilde{\xi}$ constructed above satisfy the conditions in Proposition \ref{prop:7.2.2s}. Let $\Gamma\in \mathcal{GC}^0_{\mathfrak{g}}$ be a cocycle, that is, $\delta \Gamma=0$, then by \eqref{eq:7.2.7ssss}, we get
\begin{equation}
 dZ_\Gamma(\widetilde{\omega},\widetilde{\xi})=Z_{\delta \Gamma}(\widetilde{\omega},\widetilde{\xi})=0.
\end{equation}
Therefore,
\begin{equation}
Z_\Gamma(\widetilde{\omega},\widetilde{\xi})|_{\tau=0}=Z_\Gamma(\widetilde{\omega},\widetilde{\xi})|_{\tau=1}.
\end{equation}
Then by \eqref{eq:6.2.4paris} and \eqref{eq:7.2.8kk}, we get $Z_\Gamma(\omega)=Z_\Gamma(\omega')$. The proof is completed.\qed

\subsection{Proof of Theorem \ref{thm:6.3.1}}\label{ss7.3s}
The proof of Theorem \ref{thm:6.3.1} (1) goes along a similar line to that of Theorem \ref{thm:6.2.5ss}, in particular, we need to check the following claims. 
\begin{enumerate}[(i)]
\item For an adapted propagator $\omega^\sharp$, we need to show that the map $Z_{-}(\omega^\sharp)$ factors through the space of decorated trivalent graphs without self-loops.
\item (Analogous to Subsection \ref{sss7.2.2}) For two different adapted propagators $\omega^\sharp$, $\omega^{\sharp,\prime}$, a closed form $\widetilde{\omega}^\sharp$ on $I\times C_2(M)$ can be constructed to connect smoothly $\omega^\sharp$ and $\omega^{\sharp,\prime}$, and that for $\tau\in I$, $\widetilde{\omega}^\sharp|_{\{\tau\}\times C_2(M)}$ is an adapted propagator.

\item Applying Stokes' theorem to one parameter family $Z_{-}(\widetilde{\omega}^\sharp,\widetilde{\xi}^\sharp)$ as defined in \eqref{eq:7.2.5sept}, we need to show that the vanishing of the contributions of the boundary terms corresponding to the collapse of two distinct vertices connected by a {\it non-regular} edge, so that $Z_{-}(\omega^\sharp)$ is invariant associated with the cocylces in $\mathcal{G}^0_{\mathfrak{g}: n}$.
\end{enumerate}

 Note that (i) follows immediately from Theorem \ref{thm:6.3.2}. For (ii), combining the proof of Theorem \ref{thm:enhanced} with the proof in Subsection \ref{sss7.2.2}, we can construct a closed $2$-form $\widetilde{\omega}^\sharp\in \Omega^2_-(I\times C_2(M); F_\rho)$ and closed $2$-form $\widetilde{\xi}^\sharp\in \Omega^{2}_{-}(I\times \Delta; E_\rho\otimes E_\rho)$ such that
\begin{itemize}
    \item $\widetilde{\omega}^\sharp_{\{0\}\times C_2(M)}=\omega^\sharp$, $\widetilde{\omega}^{\sharp}_{\{1\}\times C_2(M)}=\omega^{\sharp,\prime}$, then, $\widetilde{\xi}^{\sharp}_{\{0\}\times C_2(M)}=\xi^{\sharp}$, $\widetilde{\xi}^{\sharp,\prime}_{\{1\}\times C_2(M)}=\xi^{\sharp,\prime}$;
    \item Analogous to \eqref{eq:6.1.2}, we have
    \begin{equation}
\widetilde{\mathfrak{i}_\partial}^\ast(\widetilde{\omega}^{\sharp})=I(\widetilde{\eta})+\widetilde{q_\partial}^\ast (\widetilde{\xi}^{\sharp})
    \end{equation}
    with $\mathfrak{L}(\widetilde{\xi}^\sharp)=0$ on $I\times \Delta$.
\end{itemize}

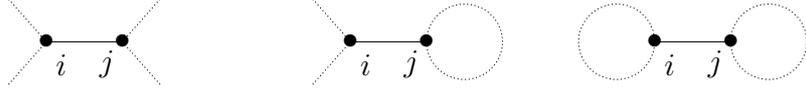
\begin{figure}[h]
\captionsetup{margin=2cm}
\centering
\begin{tikzpicture}
\begin{scope}
    \draw (-0.5,0.0)--(-1.5,0.0);
    \draw[densely dotted] (-1.5,0.0)--(-2.0,0.5*1.14);
    \draw[densely dotted] (-1.5,0.0)--(-2.0,-0.5*1.14);
    \draw[densely dotted] (-0.5,0.0)--(0,0.5*1.14);
    \draw[densely dotted] (-0.5,0.0)--(0,-0.5*1.14);
    \node at (-0.5, 0) {$\bullet$};
     \node at (-1.5, 0) {$\bullet$};
      \node at (-0.7, -0.3) {$j$};
       \node at (-1.3, -0.3) {$i$};
\end{scope}
\begin{scope}[xshift=4cm]
    \draw[densely dotted] (0,0) circle (0.5cm);
    \draw (-0.5,0.0)--(-1.5,0.0);
    \draw[densely dotted] (-1.5,0.0)--(-2.0,0.5*1.14);
    \draw[densely dotted] (-1.5,0.0)--(-2.0,-0.5*1.14);
     \node at (-0.5, 0) {$\bullet$};
     \node at (-1.5, 0) {$\bullet$};
      \node at (-0.7, -0.3) {$j$};
       \node at (-1.3, -0.3) {$i$};
\end{scope}
\begin{scope}[xshift=8cm]
    \draw[densely dotted] (0,0) circle (0.5cm);
    \draw (-0.5,0.0)--(-1.5,0.0);
    \draw[densely dotted] (-2,0) circle (0.5cm);
     \node at (-0.5, 0) {$\bullet$};
     \node at (-1.5, 0) {$\bullet$};
      \node at (-0.7, -0.3) {$j$};
       \node at (-1.3, -0.3) {$i$};
\end{scope}
\end{tikzpicture}
\caption[Three cases for an internal edge to be contracted]{Parts of trivalent graphs with regular edge $(ij)$ connecting two collapsing vertices corresponding to $(\bb{e}_0, \bb{e}_v, \bb{e}_h)=(4,1,0), (2, 1,1), (0, 1, 2)$ respectively. Here, such regular edges are depicted as solid lines.}
\label{fig:6.4_three}
\end{figure}

Therefore, it suffices to show (iii). For this, under the same arguments as in the proof of Proposition \ref{prop:7.2.2s}, we investigate in a more detailed manner the case that $\ell=|S|=2$ and $S$ correspond exactly to an edge in the graph. Note that a variation formula like \eqref{eq:7.2.7ssss} can be deduced for $Z_{\Gamma}(\widetilde{\omega}^\sharp,\widetilde{\xi}^\sharp)$ with arbitrary $\Gamma\in \mathcal{G}^0_{\mathfrak{g}: n}$, but we now focus on the proof of Theorem \ref{thm:6.3.1}, so that we will assume in the sequel that $\Gamma\in \mathcal{G}^0_{\mathfrak{g}: n}$ is a cocycle, that is, $\delta^\sharp \Gamma=0$. There are the following cases (a) and (b), where we use the notation introduced after \eqref{eq:7.2.15sss}:

\textbf{(a).} One of the edges connecting two collapsing vertices $\{i,j\}$ is regular; in this case, there are three types of local graphs corresponding to  $(\bb{e}_0, \bb{e}_v, \bb{e}_h)=(4,1,0), (2, 1,1), (0, 1, 2)$ as in Fig. \ref{fig:6.4_three}. 

For $(\bb{e}_0, \bb{e}_v, \bb{e}_h)=(2, 1,1), (0, 1, 2)$, the corresponding graphs must have at least one self-loop edge, so there is nothing to show. For the remaining case $(\bb{e}_0, \bb{e}_v, \bb{e}_h)=(4,1,0)$, the graph cocycle condition ($\delta^\sharp \Gamma=0$) gives a cancellation of these boundary contributions.

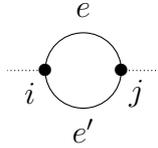
\begin{figure}[h]
\captionsetup{margin=2cm}
\centering
\begin{tikzpicture}
    \draw (0,0) circle (0.5cm);
   \draw[densely dotted] (0.5,0.0)--(1.0,0.0);
    \draw[densely dotted] (-0.5,0.0)--(-1.0,0.0);
    \node at (0.5, 0) {$\bullet$};
     \node at (-0.5, 0) {$\bullet$};
     \node at (0, 0.8) {$e$};
     \node at (0, -0.8) {$e'$};
      \node at (-0.7, -0.3) {$i$};
      \node at (0.7, -0.3) {$j$};
\end{tikzpicture}
\caption[Case of the contraction yielding a self-loop edge]{Two non-regular edge $e$ and $e'$ which connect the vertices $i$ and $j$. The contraction of the edge $e$ or $e'$  yields one self-loop edge. 
}
\label{fig:6.4.1}
\end{figure}

\textbf{(b).}  The edge connecting two collapsing vertices $i$ and $j$ is not regular; This case is further divided into two cases (b-1) and (b-2): 

\textbf{Case (b-1)}: the number of such non-regular edges is 2 as Fig. \ref{fig:6.4.1} (case that $\bb{e}_0=2, \bb{e}_v=2, \bb{e}_h=0$);
in this case integrand $\Tr_i \otimes \Tr_j (\widetilde{\omega}^\sharp_{(\bullet i)} \wedge (\widetilde{\omega}^\sharp_{(ij)})^2 \wedge \widetilde{\omega}^\sharp_{(j \bullet)})$ associated to edges connecting $i$ and $j$ (we may assume $i<j$) restricts to
\begin{eqnarray}
    &&\Tr_i \otimes \Tr_i \left[ \widetilde{\omega}^\sharp_{(\bullet i)}\wedge \big(I(\widetilde{\eta})_{(ii)} + \widetilde{q}^{\ast} \widetilde{p}^{\ast}_i \widetilde{\xi}^\sharp\big)^2 \wedge \widetilde{\omega}^\sharp_{(i \bullet)}\right]\\
    &=&\Tr_i \otimes \Tr_i \left[\widetilde{\omega}^\sharp_{(\bullet i)} \wedge\big(I(\widetilde{\eta})_{(ii)}^2 + 2 \widetilde{p}^{\ast} I(\widetilde{\eta})_{(ii)}\wedge \widetilde{q}^{\ast} \widetilde{p}^{\ast}_i \widetilde{\xi}^\sharp + \widetilde{q}^{\ast} \widetilde{p}^{\ast}_i (\widetilde{\xi}^\sharp)^2\big)\wedge \widetilde{\omega}^\sharp_{(i \bullet)}\right]
\end{eqnarray}
on the boundary component. 

Then, performing fiber-wise integration along the fiber $\mathbb{S}^2$, since $I(\widetilde{\eta})^2 = 0$ on $\mathbb{S}^2$ and $\widetilde{q}^{\ast}\widetilde{p}_i^{\ast} \widetilde{\xi}^2$ is degree 0 along the fiber, we only need to compute $2\Tr_i \otimes \Tr_i[\widetilde{\omega}^\sharp_{(\bullet i)} \wedge \mathbf{1} \widetilde{q}^{\ast}\widetilde{p}_i^{\ast} \widetilde{\xi}^\sharp \wedge \widetilde{\omega}_{(i \bullet)}]$. We need the following Lemma which is a variant of \cite[Lemma 4.4]{CS}. 
\begin{lemma}\label{lem:6.4.1}
With the same notations as above, we have the following equation.
    \begin{equation}
        \Tr_i \otimes \Tr_i\Big[\widetilde{\omega}^\sharp_{(\bullet i)} \wedge \mathbf{1} \widetilde{q}^{\ast}\widetilde{p}_i^{\ast} \widetilde{\xi}^\sharp \wedge \widetilde{\omega}^\sharp_{(i \bullet)}\Big] = \frac{1}{2} B_i\left(\mathfrak{L}_i(\widetilde{q}^{\ast}\widetilde{p}_i^{\ast} \widetilde{\xi}^\sharp), \mathfrak{L}_i(\widetilde{\omega}^\sharp_{(\bullet i)} \wedge \widetilde{\omega}^\sharp_{(i \bullet)})\right)
    \end{equation}
    where $B_i, \mathfrak{L}_i$ means that such operations occur at the vertex $i$.
\end{lemma}
\begin{proof}
    It suffices to show the claim fiber-wise, that is, for $\mathfrak{g}^{\otimes 3} \otimes \mathfrak{g}^{\otimes 3}$. Since $T^{\ast}$ acts on the diagonal subspace $\Delta$ by identity, we have $\Omega^{\bullet}_{-}(\Delta; E_{\rho}\otimes E_{\rho}) = \Omega^{\bullet}(\Delta; (E_{\rho}\otimes E_{\rho})_{-})= \Omega^{\bullet}(\Delta; \Lambda^2 E_{\rho})$. Let $e_1,\ldots, e_{\dim \mathfrak{g}}$ be basis of $\mathfrak{g}$ which is normalized with the condition $B(e_\ell,e_j)=\varepsilon_\ell\delta_{\ell j}$, $\varepsilon_\ell\in\{1,-1\}$. Then the Casimir element is given as
    $$\mathbf{1}=\sum_i \varepsilon_i e_i\otimes e_i.$$
    
    Then, for each fiber at $x \in M$, $\Lambda^2 E_{\rho, _x} = \Lambda^2 \mathfrak{g}$ has a basis $\{e_\ell \otimes e_j -e_j \otimes e_\ell \mid 1 \leq \ell < j \leq \dim \mathfrak{g}\}$. Thus, in terms of this basis, one obtains, for some $e_a$ and $e_b$ corresponding to the components of $\widetilde{\omega}^\sharp_{(\bullet i)}$ and $\widetilde{\omega}^\sharp_{(i \bullet)}$ respectively,
    \begin{equation}
    \begin{split}
        & \Tr \otimes \Tr\Big[ (e_\ell \otimes e_j -e_j \otimes e_\ell) \otimes (e_a \otimes e_b) \otimes (\sum_{n=1}^{\dim \mathfrak{g}} \varepsilon_n e_n \otimes e_n)\Big]\\
        =& \sum_{n=1}^{\dim \mathfrak{g}} B([e_\ell, e_a], \varepsilon_n e_n) B(e_n, [e_j, e_b]) - \sum_{n=1}^{\dim \mathfrak{g}} B([e_j, e_a], \varepsilon_n e_n) B(e_n, [e_\ell, e_b])\\
       =& - B([e_\ell,e_j], [e_b, e_a]) - B([e_\ell,e_b], [e_a, e_j])- B([e_j, e_a], [e_\ell, e_b])\\
        =& B([e_\ell,e_j], [e_a, e_b])\\
        =& \frac{1}{2} B(\mathfrak{L}(e_\ell \otimes e_j - e_j \otimes e_\ell), \mathfrak{L}(e_a \otimes e_b))
        \end{split}
    \end{equation}
    where the third equality follows from the Jacobi identity.
\end{proof}

By Lemma \ref{lem:6.4.1} and $\mathfrak{L}(\widetilde{\xi}^\sharp) = 0$, we conclude that the factor 
\begin{equation}
    2 \Tr_i \otimes \Tr_i\big[\widetilde{\omega}^\sharp_{(\bullet i)} \wedge \mathbf{1} \widetilde{q}^{\ast}\widetilde{p}_i^{\ast} \widetilde{\xi}^\sharp \wedge \widetilde{\omega}^\sharp_{(i \bullet)}\big]= B_i(\mathfrak{L}_i(\widetilde{q}^{\ast}\widetilde{p}_i^{\ast}\widetilde{\xi}^\sharp), \mathfrak{L}_i(\widetilde{\omega}^\sharp_{(\bullet i)} \wedge \widetilde{\omega}^\sharp_{(i \bullet)})) = 0.
\end{equation}

\textbf{Case (b-2)}: the number of such non-regular edges is 3 (case that $\bb{e}_0=0, \bb{e}_v=3, \bb{e}_h=0$); note that this case occurs only when the given connected trivalent graph is the theta graph. The integrand associated with edges connecting $i$ and $j$ becomes
\begin{equation}
\begin{split}
    &\Tr_i \otimes \Tr_i\big[\big(I( \widetilde{\eta})_{(ii)} + \widetilde{q}^{\ast} \widetilde{p}^{\ast}_i \widetilde{\xi}^\sharp\big)^3\big]\\
    =&\Tr_i \otimes \Tr_i\big[I(\widetilde{\eta})_{(ii)}^3 + 3 I(\widetilde{\eta})_{(ii)}\wedge q^{\ast} \pi^{\ast}_i(\widetilde{\xi}^\sharp)^2 + 3 I(\widetilde{\eta})_{(ii)}^2\wedge \widetilde{q}^{\ast} \widetilde{p}^{\ast}_i \widetilde{\xi}^\sharp + \widetilde{q}^{\ast} \widetilde{p}^{\ast}_i (\widetilde{\xi}^\sharp)^3\big].
    \end{split}
    \label{eq:5.3.7-july-24}
\end{equation}
After integrating along the fiber $\mathbb{S}^2$, only the term $3\Tr_i \otimes \Tr_i \big[\mathbf{1}\widetilde{q}^{\ast}\widetilde{p}_i^{\ast} (\widetilde{\xi}^\sharp)^2\big]$ remains. Then by a similar calculation as in Lemma \ref{lem:6.4.1} and condition $\mathfrak{L}(\widetilde{\xi}^\sharp) = 0$, we get
\begin{equation}
    \Tr_i \otimes \Tr_i \big[\mathbf{1}\widetilde{q}^{\ast}\widetilde{p}_i^{\ast} (\widetilde{\xi}^\sharp)^2\big]= \frac{1}{2} B_i(\mathfrak{L}(\widetilde{q}^{\ast}\widetilde{p}_i^{\ast} \widetilde{\xi}^\sharp), \mathfrak{L}(\widetilde{q}^{\ast}\widetilde{p}_i^{\ast} \widetilde{\xi}^\sharp))=0,
\end{equation}
so the configuration space integrals of the quantity in \eqref{eq:5.3.7-july-24} vanish.
 This completes the proof of Theorem \ref{thm:6.3.1} (1).

Next, we show Theorem \ref{thm:6.3.1} (2). By Theorem \ref{thm:6.2.5ss},  $Z(M, \rho, [f])(\Gamma)$ is independent of the choice of propagators. Hence, by using an adapted propagator instead of a general propagator, we obtained the same invariant as $Z(M, \rho, [f])(\Gamma)$, but associated with $\Gamma'$ obtained by removing graphs with self-loops from $\Gamma$. As we have shown in Theorem \ref{thm:6.3.1} (1), $Z^\sharp(M, \rho, [f])(\Gamma')$ is an invariant associated with cocycle $\Gamma'$ in $\mathcal{G}^0_{\mathfrak{g}}$. Thus, we obtain the commutative diagram \eqref{eq:7.2.1}. This completes the proof.\qed

\appendix
\section{Prerequisite for compactification of configuration spaces of $3$-manifolds}\label{sectionsecond}
\renewcommand{\thesubsection}{\thesection.\arabic{subsection}}
This appendix recollects the prerequisite for manifold with corners, fiberwise Stokes formula,  compactification of configuration spaces, and framing of closed orientable 3-manifolds  to fix our notation and convention. 
\subsection{Manifolds with corners, fiber integration, and fiberwise Stokes' formula}\label{ss1.1corners}
To describe the geometry of the compactified configuration spaces, we will use smooth manifolds with corners. 
For a detailed introduction to the manifold with corners, we refer to \cite{Joyce2012}, a partial note of Melrose \cite{Melrose}, and the references therein. For fiber integrations, see also \cite[Chapter VII]{GHSV}.

Set $\bR_+=[0,\infty[\;$.
The model spaces to build a manifold with corners of dimension $m$ are 
$$\bR^m_k  := \bR_+^k \times \bR^{m-k},\; 0\leq k\leq m.$$
The chart given by an open neighbourhood of the origin in $\bR^m_k  := \bR_+^k \times \bR^{m-k},\; 0\leq k\leq m$ is called $m$-dimensional chart with corners.


\begin{definition}[manifold with corners]
    A topological manifold $X$ together with a smooth structure defined by a maximal atlas consisting of compatible open $m$-dimensional charts with corners is called a (smooth) manifold with corners of dimension $m$. 
\end{definition}
One key technique to do the analysis on a manifold with corners is the Seeley Extension Theorem, which states that smooth functions on any open subset in $\bR^m_k$ always extend to its neighborhood in $\bR^m$. This result guarantees the compatibility of gluing the local charts with corners by diffeomorphisms.
For a point $x=(x_1,\ldots,x_k,x_{k+1},\ldots, x_m)\in \bR^{m}_{k}$, set $\mathrm{depth}   (x)\in\mathbb{N}$ to be $|\{ j \mid x_j = 0\}|$, $j=0,1,\ldots,k$. For an open subset $U\subset\bR^m_k$, set
\begin{equation}
    S^j(U):=\{x\in U\;|\; \mathrm{depth}_U(x)=j\}.
\end{equation}
In particular, $S^0(U)$ is the interior of $U$ in $\bR^m$. Set $\mathrm{depth}(U)=\max\{j \;|\; S^j (U)\neq \emptyset\}$. 

If $X$ is a (smooth) manifold with corners, then we can define a canonical function $\mathrm{depth}_X: X\rightarrow \mathbb{N}_0$ which associates a point $x\in X$ with its depth $\mathrm{depth}_X(x)$ by taking a compatible chart as above. The depth $j$-stratum of $X$ is the subset
\begin{equation}
    S^j(X):= \{x\in X\;|\; \mathrm{depth}_X(x)=j\}.
\end{equation}

We use $\partial X$ denote the topological boundary of $X$ (homeomorphic to $\bR_+\times \bR^{m-1}$ at origin), then as subset of $X$ we have
$$\partial X=\overline{S^1(X)}=\cup_{j\geq 1} S^j(X).$$
In general, $\partial X$ is not a smooth manifold with corners, so that we introduce a regularized boundary $\partial^* X$ of $X$ which has a structure of manifold with corners induced from $X$.

For $x\in X$, a local boundary component $\beta$ of $X$ at $x$ is a local choice of the connected component of $S^1(X)$ in a small open neighbourhood of $x$ in $X$.  The total number of local boundary components at $x$ is exactly $\mathrm{depth}_X(x)$. Set
$$\partial^*X=\{(x,\beta)\;|\; x\in \partial X, \beta \text{ a local boundary component at }x\}.$$
We have the following facts.
\begin{itemize}
    \item The charts with corners of $X$ give canonically the $(m-1)$-dimensional charts with corners on $\partial^*X$ so that $\partial^* X$ is a smooth manifold with corners of dimension $m-1$ with $S^0(\partial^*X)=S^1(X)$.
    \item
     We have a continuous projection $\partial^* X\rightarrow \partial X$ which sends $(x,\beta)$ to $x$. Composing it with inclusion $\mathfrak{i}_\partial:\partial X\hookrightarrow X$, we get a map
    $$\mathfrak{i}_\partial: \partial^* X\rightarrow X.$$
    This map $\mathfrak{i}_\partial$ is a smooth map between manifolds with corners, it is not necessary to be injective, for $x\in X$, we have
    $$|\mathfrak{i}_\partial^{-1}(x)|=\mathrm{depth}_X(x).$$
\end{itemize}

The tangent bundle $TX$ of $X$ as well as the cotangent bundle $T^*X$ is defined as in the manifold case by local charts. Then the smooth vector fields and the differential forms on $X$ are well-defined. The usual partition of unity still holds, so that if $\alpha$ is a smooth $m$-form on the oriented $X$ with compact support, the integration $\int_X \alpha$ is well-defined by considering the integration on local charts. Note that we always have $$\int_X \alpha=\int_{S^0(X)}\alpha.$$

Now let $X$ be a compact orientable (smooth) manifold with corners equipped with an orientation $o(X)$. Let $\mathbf{n}$ denote the outward normal vector 
field along $S^1(X)$ in $X$. We orient $\partial^*X$ by 
an orientation $o_{\mathrm{ind}}(\partial^* X)$ such that at all points of $S^1(X)=S^0(\partial^* X)$,
\begin{equation}
	o(X)|_{\partial^* X}=\mathbf{n}\wedge o_{\mathrm{ind}}(\partial^*X).
	\label{eq:2.1.1pf}
\end{equation}
Take 
$\alpha\in\Omega^{\bullet}(X)$, then $\mathfrak{i}_\partial^*\alpha$ is a smooth form on $\partial^\ast X$, and we have the Stokes' formula 
\begin{equation}
\int_{X}d\alpha=\int_{S^1(X)}\mathfrak{i}^{\ast}_\partial \alpha=\int_{\partial^*X}\mathfrak{i}_\partial^\ast\alpha.
	\label{eq:2.1.2pf}
\end{equation}
When there is no confusion, we can simply write the right-hand side of \eqref{eq:2.1.2pf} as $\int_{\partial X} \alpha$.


\begin{definition}
    A smooth map $f: X\rightarrow B$ is called a submersion if for all $x\in X$ with $x\in S^k(X)$, $f(x)\in S^\ell(B)$, the tangent maps $df_x: T_x X\rightarrow T_{f(x)}B$ and $df_x: T_x S^k(X)\rightarrow T_{f(x)}S^\ell(B)$ are surjective.
\end{definition}

Analogous to the usual Ehresmann's theorem (also cf. \cite[Section 5]{Joyce2012}), a proper submersion $f: X\rightarrow B$ for manifolds with corners is a locally-trivial fibration on $S^0(B)$ where the fibres are compact manifolds with corners. Note that, in general, the locally-trivial fibration could not extend to the corners of $B$, a simple counterexample is as follows: consider the submersion $f:\bR_+\times\bR_+\ni (x,y)\mapsto \frac{1}{\sqrt{2}}(x+y)\in\bR_+$, which is not locally-trivial fibration, since $f^{-1}(0)=\{(0,0)\}$, but $f^{-1}(1)$ is a nontrivial segment of line.

Fix a surjective submersion $p: (X,o(X))\rightarrow (B,o(B))$ of compact oriented manifolds with corners.  Then there is a 
unique orientation $o_{\mathrm{fibre}}(p)$ on the fibres of $p$, i.e., orientations on $p^{-1}(b)$, $b\in S^0(B)$ such that locally
\begin{equation}
	o(X)=o(B)\wedge o_{\mathrm{fibre}}(p).
	\label{eq:2.3.2pf}
\end{equation}

\begin{proposition}\label{prop:2.2.3p}
	Let $p: X\rightarrow B$ be a surjective submersion of compact orientable manifolds with corners. Then for any $\alpha\in 
	\Omega^{\bullet}(X)$, there exists a unique smooth form $\hat{\alpha}\in\Omega^\bullet(B)$ such that for $\gamma\in \Omega^{\bullet}(B)$, we have
	\begin{equation}
		\int_{(X,o(X))}p^{\ast}(\gamma)\wedge \alpha=\int_{(B,o(B))}\gamma\wedge 
		\hat{\alpha}.
		\label{eq:2.3.4pf}
	\end{equation}
\end{proposition} 

\begin{definition} \label{def:2.4.1}
The linear map $p_\ast: \Omega^\bullet(X)\ni\alpha\mapsto p_\ast(\alpha):=\hat{\alpha}\in\Omega^\bullet(B)$ is called the fibre integration of the submersion $p: X\rightarrow B$.
\end{definition}

Given a proper surjective submersion $p: X\rightarrow B$ of oriented manifolds with corners, we get a surjective submersion $p_0: X_0:=p^{-1}(S^0(B))\rightarrow S^0(B)$ for smooth manifolds. Define $X_0^\partial = \partial^*X_0$. We get a surjective submersion of manifolds with corners
$$p^\partial: X_0^\partial\rightarrow S^0(B).$$
Note that the orientation $o(X)$ induces an orientation $o(X^\partial_0) = o(B) \wedge o_\mathrm{fibre}(p^\partial)$, where the fibration orientation $o_\mathrm{fibre}(p^\partial)$ for $p^\partial$ is induced from $o_{\mathrm{fibre}}(p)$ as \eqref{eq:2.1.1pf}.


\begin{lemma}
    For $\alpha\in\Omega^\bullet(X)$, let $\mathfrak{i}_\partial^*\alpha$ denote the corresponding smooth form on $X^\partial_0$, the fibre integration $p^\partial_*( \mathfrak{i}_\partial^\ast\alpha)\in \Omega^\bullet(S^0(B))$ always extends smoothly to $B$, which we denote it by the same notation.
\end{lemma}

With the above notation, we have the following fibrewise Stokes' formula. See also \cite[Chapter VII]{GHSV}.
\begin{proposition}\label{prop:2.4.1}
	For $p\in\mathbb{N}$ and $\alpha\in\Omega^{p}(X)$, we have the following identity of smooth forms on $B$,
	\begin{equation}
		dp_{\ast}(\alpha)=p_{\ast}(d\alpha)+(-1)^{p+\dim X - \dim B +1} 
		p^{\partial}_{\ast}(\mathfrak{i}_\partial^\ast\alpha)\in 
		\Omega^{1+p+\dim B -\dim X}(B).
		\label{eq:2.3.5pf}
	\end{equation}
\end{proposition}

\subsection{Compactification of configuration spaces}\label{app:FMAS}
We recall briefly the Fulton-MacPherson compactification of configuration spaces of a closed oriented smooth $3$-manifold. For more details see \cite{FM} and \cite{AS2}. Let $M$ be a closed oriented $3$-manifold with a given orientation $o(M)$. Let $F\rightarrow M$ be a real vector bundle over $M$ with rank $r\geq 2$, and let $SF_x=\bR^*_{+}\backslash(F_x-\{0\})$, $x\in M$, be the sphere bundle of $F$. If $g^F$ is a Euclidean metric on $F$, then we can canonically identify $SF$ with the unit sphere bundle $S_{g^F}F\rightarrow M$ of $(F, g^F)$.

Let $S$ be a finite set, and put
\begin{equation}
	M^{S}=\mathrm{Maps}(S,M)=\Pi_{i\in S}M_{i},
	\label{eq:1.2.1}
\end{equation}
where $M_{i}=M$ is just a copy of $M$.
If $S=\underline{n}:=\{1,2\ldots,n\}$, $n\geq 2$, we simply denote $M^{n}=M^{S}$ to be compatible with the usual notation.

Put $\Delta_{S}\simeq M$ be the subset of $M^{S}$ consisting of constant maps, which is called the principal diagonal of $M^{S}$. Let $B\ell(M^{S},\Delta_{S})$ be the 
geometric blow-up of $M^{S}$ along $\Delta_{S}$. It can be regarded 
as replacing $\Delta_{S}$ by its sphere normal bundle 
$S\nu_{\Delta_{S}}$ of $\Delta_{S}$ in $M^{S}$. If 
$U_{S}$ is an open small tubular neighbourhood of $\Delta_{S}$ in $M^{S}$, then 
$B\ell(M^{S},\Delta_{S})$ is diffeomorphic to $M^{S}\backslash 
U_{S}$ as manifolds with boundary where $\partial B\ell(M^{S},\Delta_{S})\simeq S\nu_{\Delta_{S}}$. 
Moreover, we have a canonical smooth projection, the blow-down map, $B\ell(M^{S},\Delta_{S})\rightarrow 
M^{S}$, which, when restricting to the boundary, is given by the projection 
$S\nu_{\Delta_{S}}\rightarrow \Delta_{S}$.

For $n\in \mathbb{N}^\ast$, we denote by  $\mathrm{Conf}_n(M)$ the $n$-point configuration space of $M$, i.e.,
\begin{equation}
	\mathrm{Conf}_n(M):= \{ (x_1, \ldots, x_n) \in M^n | x_i \neq x_j (i \neq j)\}.
	\label{eq:1.2.2}
\end{equation}
We have an injective smooth map:
\begin{equation}
	\Phi_{n}:\mathrm{Conf}_n(M)\rightarrow \mathcal{B}:=M^n \times \prod_{S 
\subset \underline{n}, |S| \geq 2} B\ell(M^S, \Delta_S).
\label{eq:1.2.3}
\end{equation}
Note that the target space of $\Phi_{n}$ is a compact manifold with corners as described in \ref{ss1.1corners}.

The compactification $C_n(M)$ of $\mathrm{Conf}_n(M)$ is defined 
as the closure of the image of $\mathrm{Conf}_n(M)$ via $\Phi_{n}$ 
equipped with the induced smooth structure, that is,
\begin{equation}
	C_n(M) := \overline{\Phi_{n}(\mathrm{Conf}_n(M))} \subset \mathcal{B}=M^n \times \prod_{S 
\subset \underline{n}, |S| \geq 2} B\ell(M^S, \Delta_S).
\label{eq:1.2.4}
\end{equation}
By \cite[Section 5]{AS2}, $C_{n}(M)$, as an embedded submanifold of $\mathcal{B}$, is a compact manifold with corners of dimension $3n$. In particular, the interior of $C_{n}(M)$ is exactly $\mathrm{Conf}_n(M)$.

To understand the structure of $C_n(M)$ as a manifold with corners, we introduce the following notation. For $S\subset \underline{n}$ with $|S|>1$, for $x=(z,\ldots,z)\in\Delta_S$ with $z\in M$, the normal bundle $\nu_{\Delta_S, x}$ of $\Delta_S$ in $M^S$ at $x$ can be identified with the quotient space $(T_z M)^S/T_zM$. For $u_S=(u_i)_{i\in S}\in (T_z M)^S$, if all components $u_i$ are identical, then $[u_S]=0\in (T_z M)^S/T_zM$. For each element $[u_S]\in (T_z M)^S/T_zM$, there exists a unique representative $u_S=(u_i)_{i\in S}\in (T_z M)^S$ such that $\sum_{i\in S}u_i=0$ in $T_zM$. We can also regard such vector $u_S \in (T_z M)^{\underline{n}}$ by setting $u_j=0$ for $j\not\in S$. Note that $\bR^*_+=\;]0,+\infty]$ acts on $(T_z M)^S/T_zM$ by the diagonal re-scaling on all the components. So we have the identification $S\nu_{\Delta_S}\simeq \bR^*_+\backslash \big((T_z M)^S/T_zM - \{0\}\big)$. Then any nonzero $[u_S]\in (T_z M)^S/T_zM$ corresponds to an element in $S\nu_{\Delta_S}$, which is still denoted by $[u_S]$ if there is no confusion.


Let $q: B\ell(M^S,\Delta_S)\rightarrow M^S$ denote the obvious projection. We always use $x_{B,S}$ to denote a point in $B\ell(M^S,\Delta_S)$, such that if $x_S=q(x_{B,S})\not\in \Delta_S$, then $x_{B,S}=x_S$; otherwise, $x_{B,S}=(x_S,[u_S])$ where $x_S=(z,\ldots, z)\in \Delta_S$, $0\neq [u_S]\in (T_z M)^S/T_zM$.

As a point set, we have a characterization of $C_n(M)$ as a subset of $\mathcal{B}$: $C_n(M)$ are consisting of all points $(x,\{x_{B,S}, |S|\geq 2\})$ in $\mathcal{B}$ which satisfy the following two conditions:
\begin{itemize}
    \item $x_S=q(x_{B,S})=x|_S$, for $S\subset \underline{n}$, $|S|>1$, where $x|_S=(x_i)_{i\in S}\in M^S$ with $x=(x_1,\ldots, x_n)\in M^n$.
    \item For any subset $S$ ($|S|\geq 3$) with $x_S\in \Delta_S$, write $x_{B,S}=(x_S, [u_S])$, then for each subset $S'\subset S$ with $|S'|\geq 2$, if $S'$-components of $u_S$ are not all equal, we have $x_{B,S'}=(x_{S'}, [u_S|_{S'}])$.
\end{itemize}

Set $\mathbf{S}_n=\{S\subset \underline{n}\;|\; |S|\geq 2\}$.
\begin{definition}
    A subset $\mathcal{S}\subset \mathbf{S}_n$ is called nested if any two elements $S_1, S_2\in \mathcal{S}$ are either disjoint or else one contains the other. 
\end{definition}
For a nested subset $\mathcal{S}\subset \mathbf{S}_n$, the open strata $M(\mathcal{S})^\circ$ of $C_n(M)$ is defined as follows, it consists of the points $(x,\{x_{B,S}, |S|\geq 2\})\in C_n(M)$ such that (1) $x|_S\in \Delta_S$ if and only if $S\subset S'$ for some $S'\in\mathcal{S}$; (2) For $S'$ ($|S'|>1$) with a minimal element $S\in\mathcal{S}$ such that $S'\subset S$, then $[u_{S'}]=[u_S|_{S'}]$ in $x_{B,S'}$; (3) For $S_1, S_2\in\mathcal{S}$ such that $S_1\subsetneq S_2$, then all $S_1$-components of $u_{S_2}$ are all equal.

In \cite[Subsection 5.3 and 5.4]{AS2}, Axelrod and Singer showed the following facts:
\begin{itemize}
    \item $M(\mathcal{S})^\circ$ is a smooth manifold of dimension $3n-|\mathcal{S}|$, in particular, $M(\emptyset)^\circ=\mathrm{Conf}_n(M)$.
\item The closed strata $M(\mathcal{S})$, defined as the closure of $M(\mathcal{S})^\circ$ in $C_n(M)$, is given as 
$$M(\mathcal{S})=\cup_{\mathcal{T}\supset \mathcal{S}} M(\mathcal{T})^\circ$$
where $\mathcal{T}$ runs over all nested subsets of $\mathbf{S}_n$ that contain $\mathcal{S}$.
\item We have
$C_n(M)=\cup_{\mathcal{S}\;\mathrm{nested}} M(\mathcal{S})^\circ.$
\item For the strata of $C_n(M)$ as manifold with corners, we have for $k=0, 1, \ldots, 3n$,
\begin{equation}
    S^kC_n(M)=\cup_{\mathcal{S},\;|S|=k} M(\mathcal{S})^\circ.
\end{equation}
Then the interior of $\partial^* C_n(M)$ is given by all the single sets $S$ of $\underline{n}$ with $|S|\geq 2$.
\end{itemize}


Now we focus on the case of $n=2$, then $C_2(M)$ is a compact manifold with boundary. Let $\Delta\subset 
M\times M$ denote the diagonal. An elementary argument 
shows that
\begin{equation}
C_2(M) = (M \times M \setminus \Delta) \cup 
S\nu_{\Delta}=B\ell(M^{2},\Delta).
\label{eq:1.2.5}
\end{equation}
The blow-down map $q : C_2(M) \rightarrow M^2$ satisfies $q(S\nu_{\Delta}) = \Delta$ and $q = \Id$ otherwise. Note that the sphere normal bundle $S\nu_{\Delta}$ is given by the equivalent classes of the elements $((x,x),(-v,v))$, $x\in M, v\in T_xM$. Then it can be identified with the sphere tangent bundle $S(TM)$ by
\begin{equation}
	S\nu_{\Delta} \overset{\sim}{\rightarrow} S(TM), \quad ((x,x), 
(-v,v)) \mapsto (x, v).
\label{eq:1.2.6}
\end{equation}
We will always use the identifications $\partial C_2(M) \simeq S\nu_{\Delta}\simeq S(TM)$.

\subsection{Framings of closed orientable $3$-manifolds}\label{A.3-july}
Let $M$ be a closed orientable $3$-manifold. Then the tangent bundle $TM$ of $M$ is always parallelizable, that is, there always exists a global smooth trivialization (isomorphism of vector bundles) $f:TM\rightarrow 
M\times \mathbb{R}^{3}$ of $TM$ (see \cite{MR1509530} and \cite{MR4557281}). We call such a trivialization $f$ a framing of $M$. Fixing a smooth framing $f$ of $M$, we identify $TM$ with $M\times\bR^3$ and the sphere bundle $S(TM)$ with $M\times \mathbb{S}^2$. Set $T^V(TM)$, $T^V(S(TM))$ the vertical tangent bundles for the fibrations $TM\rightarrow M$, $S(TM)\rightarrow M$ respectively. Then the above identifications induce the splittings
\begin{equation}
\begin{split}
     &T(TM)=f^*TM\oplus  f^*T\bR^3 =:T^H_f(TM)\oplus T^V(TM),\;\\ &T(S(TM))=f^*TM\oplus f^*T\mathbb{S}^2=:T^H_f(S(TM))\oplus T^V(S(TM)).
\end{split}
\label{eq:A.3.1}
\end{equation}
A differential form $\alpha$ on $TM$ or $S(TM)$ is said to be $f$-vertical (or simply, vertical, when the framing $f$ is fixed) if for all $U\in T^H_f(TM)$ or $T^H_f(S(TM))$, we have
\begin{equation}
\iota_U\alpha=0,
\end{equation}
where $\iota_U$ denotes the contraction of $U$.

An orientation $o(M)$ gives an orientation $o(\mathbb{R}^{3})$ of $\mathbb{R}^{3}$ via $f$, and let $o_{\mathrm{ind}}(\mathbb{S}^{2})$ be the induced orientation on the unit $2$-sphere $\mathbb{S}^{2}$ with outward normal first convention, viewed as the boundary of the standard $3$-ball. Meanwhile, the complementary of open $3$-ball in $\mathbb{R}^{3}$ can be identified with the real blow-up 
$B\ell(\mathbb{R}^{3},0)$ at $0$ of 
$\mathbb{R}^{3}$ whose boundary $\partial 
B\ell(\mathbb{R}^{3},0)=\mathbb{S}^{2}$. We orient 
$B\ell(\mathbb{R}^{3},0)$ as for $\mathbb{R}^{3}$. Then
\begin{equation}
	o_{\mathrm{ind}}(\partial 
B\ell(\mathbb{R}^{3},0))=-o_{\mathrm{ind}}(\mathbb{S}^{2}).
	\label{eq:2.3.6pf}
\end{equation}
We identify $S(TM)$ with $M\times 
\mathbb{S}^{2}$ via $f$, and let $o(S(TM))$ be the orientation $o(M)\wedge o_{\mathrm{ind}}(\mathbb{S}^{2})$.



\begin{lemma}
	Under the identification $\partial C_2(M) \simeq 
	S\nu_{\Delta}\simeq S(TM)$ as explained in \eqref{eq:1.2.6}, we 
	have
	\begin{equation}
		o_{\mathrm{ind}}(\partial C_{2}(M))=o(S(TM)).
		\label{eq:2.3.8pf}
	\end{equation}
\end{lemma}
 

%


\bibliographystyle{amsalpha}
\bibliography{Construction_of_propagator}

\providecommand{\bysame}{\leavevmode\hbox to3em{\hrulefill}\thinspace}
\providecommand{\MR}{\relax\ifhmode\unskip\space\fi MR }
\providecommand{\MRhref}[2]{%
  \href{http://www.ams.org/mathscinet-getitem?mr=#1}{#2}
}
\providecommand{\href}[2]{#2}
\begin{thebibliography}{GHV72}

\bibitem[AS92]{AS}
Scott Axelrod and Isadore~M. Singer, \emph{Chern-{S}imons perturbation theory}, Proceedings of the {XX}th {I}nternational {C}onference on {D}ifferential {G}eometric {M}ethods in {T}heoretical {P}hysics, {V}ol. 1, 2 ({N}ew {Y}ork, 1991), World Sci. Publ., River Edge, NJ, 1992, pp.~3--45.

\bibitem[AS94]{AS2}
\bysame, \emph{Chern-{S}imons perturbation theory. {II}}, J. Differential Geom. \textbf{39} (1994), no.~1, 173--213.

\bibitem[BC98]{BC}
Raoul Bott and Alberto~S. Cattaneo, \emph{Integral invariants of {$3$}-manifolds}, J. Differential Geom. \textbf{48} (1998), no.~1, 91--133.

\bibitem[BC99]{BC2}
\bysame, \emph{Integral invariants of 3-manifolds. {II}}, J. Differential Geom. \textbf{53} (1999), no.~1, 1--13.

\bibitem[BN95]{BN}
Dror Bar-Natan, \emph{On the {V}assiliev knot invariants}, Topology \textbf{34} (1995), no.~2, 423--472.

\bibitem[BZ23]{MR4557281}
Valentina Bais and Daniele Zuddas, \emph{On {S}tiefel's parallelizability of 3-manifolds}, Expo. Math. \textbf{41} (2023), no.~1, 238--243.

\bibitem[CS21]{CS}
Alberto~S. Cattaneo and Tatsuro Shimizu, \emph{A note on the {$\Theta$}-invariant of 3-manifolds}, Quantum Topol. \textbf{12} (2021), no.~1, 111--127.

\bibitem[CV03]{MR2026331}
James Conant and Karen Vogtmann, \emph{On a theorem of {K}ontsevich}, Algebr. Geom. Topol. \textbf{3} (2003), 1167--1224.

\bibitem[DeB06]{DB}
Jason DeBlois, \emph{Totally geodesic surfaces and homology}, Algebr. Geom. Topol. \textbf{6} (2006), 1413--1428.

\bibitem[FH91]{FH}
William Fulton and Joe Harris, \emph{Representation theory}, Graduate Texts in Mathematics, vol. 129, Springer-Verlag, New York, 1991, A first course, Readings in Mathematics.

\bibitem[FM94]{FM}
William Fulton and Robert MacPherson, \emph{A compactification of configuration spaces}, Ann. of Math. (2) \textbf{139} (1994), no.~1, 183--225.

\bibitem[FS92]{FS}
Ronald Fintushel and Ronald~J. Stern, \emph{Integer graded instanton homology groups for homology three-spheres}, Topology \textbf{31} (1992), no.~3, 589--604.

\bibitem[GHV72]{GHSV}
Werner Greub, Stephen Halperin, and Ray Vanstone, \emph{Connections, curvature, and cohomology. {V}ol. {I}: {D}e {R}ham cohomology of manifolds and vector bundles}, Pure and Applied Mathematics, vol. Vol. 47, Academic Press, New York-London, 1972.

\bibitem[Joy12]{Joyce2012}
Dominic Joyce, \emph{On manifolds with corners}, Advances in geometric analysis, Adv. Lect. Math. (ALM), vol.~21, Int. Press, Somerville, MA, 2012, pp.~225--258.

\bibitem[Kna86]{MR0855239}
Anthony~W. Knapp, \emph{Representation theory of semisimple groups}, Princeton Mathematical Series, vol.~36, Princeton University Press, Princeton, NJ, 1986, An overview based on examples.

\bibitem[Kon93]{MR1247289}
Maxim Kontsevich, \emph{Formal (non)commutative symplectic geometry}, The {G}el'fand {M}athematical {S}eminars, 1990--1992, Birkh\"{a}user Boston, Boston, MA, 1993, pp.~173--187.

\bibitem[Kon94]{Ko}
\bysame, \emph{Feynman diagrams and low-dimensional topology}, First {E}uropean {C}ongress of {M}athematics, {V}ol. {II} ({P}aris, 1992), Progr. Math., vol. 120, Birkh\"{a}user, Basel, 1994, pp.~97--121.

\bibitem[KS23]{KS}
Teruaki Kitano and Tatsuro Shimizu, \emph{Gluing formula for an invariant related to the {C}hern-{S}imons perturbation theory}, preprint, 2023, arXiv:2306.09717.

\bibitem[Les04]{Les}
Christine Lescop, \emph{On the {K}ontsevich-{K}uperberg-{T}hurston construction of a configuration-space invariant for rational homology 3-spheres}, math.GT/0411088 Prepublication Institut Fourier \textbf{655} (2004), 71 pages.

\bibitem[Les10]{Les10}
Christine Lescop, \emph{On the cube of the equivariant linking pairing for knots and 3-manifolds of rank one}, preprint, 2010, arXiv:1008.5026.

\bibitem[Les20]{MR4521898}
Christine Lescop, \emph{Invariants of links and 3-manifolds that count graph configurations}, Winter Braids Lect. Notes \textbf{7} (2020), no.~Winter Braids X (Pisa, 2020), Exp. No. 1, 35.

\bibitem[Mel]{Melrose}
Richard Melrose, \emph{Differential analysis on manifolds with corners (in preparation)}, partially available at \url{http://www-math.mit.edu/~rbm/book.html}.

\bibitem[Mil85]{Mi}
John~J. Millson, \emph{A remark on {R}aghunathan's vanishing theorem}, Topology \textbf{24} (1985), no.~4, 495--498.

\bibitem[Por13]{P}
Joan Porti, \emph{Local and infinitesimal rigidity of representations of hyperbolic three manifolds}, 2013, p.~154–177.

\bibitem[Rag65]{R1}
Madabusi~S. Raghunathan, \emph{On the first cohomology of discrete subgroups of semisimple {L}ie groups}, Amer. J. Math. \textbf{87} (1965), 103--139.

\bibitem[Sav12]{Sa}
Nikolai Saveliev, \emph{Lectures on the topology of 3-manifolds}, revised ed., De Gruyter Textbook, Walter de Gruyter \& Co., Berlin, 2012, An introduction to the Casson invariant.

\bibitem[Shi21]{Shi3}
Tatsuro Shimizu, \emph{A geometric description of the {R}eidemeister-{T}uraev torsion of 3-manifolds}, RIMS preprint (2021), 34 pages.

\bibitem[Shi23]{Shi}
\bysame, \emph{Morse homotopy for the {$SU(2)$}-{C}hern-{S}imons perturbation theory}, J. Differential Geom. \textbf{123} (2023), no.~2, 363--390.

\bibitem[Sti35]{MR1509530}
Eduard Stiefel, \emph{Richtungsfelder und {F}ernparallelismus in n-dimensionalen {M}annigfaltigkeiten}, Comment. Math. Helv. \textbf{8} (1935), no.~1, 305--353.

\bibitem[Wit89]{MR990772}
Edward Witten, \emph{Quantum field theory and the {J}ones polynomial}, Comm. Math. Phys. \textbf{121} (1989), no.~3, 351--399.

\end{thebibliography}

\noindent
H.~Kodani\\
Institute of Mathematics for Industry,
Kyushu University, \\
744, Motooka, Nishi-ku, Fukuoka, 819-0395, Japan\\
\textit{Email address}: \texttt{kodani@imi.kyushu-u.ac.jp}

\ 

\noindent
B.~Liu\\
Universit\"at zu K\"oln, Department Mathematik/Informatik, \\ Weyertal  86-90, 50931 K\"oln, Germany\\
\textit{Email address}: \texttt{bingxiao.liu@uni-koeln.de}

\end{document}